\documentclass[12pt, eqno]{article}
\usepackage{amssymb,latexsym}
\usepackage{amsmath,latexsym}
\usepackage{amscd}

\newcommand{\K}{{\mathbb K}}

\newcommand{\N}{{\mathbb N}}

\newcommand{\R}{{\mathbb R}}

\newcommand{\Z}{{\mathbb Z}}

\newtheorem{theorem}{Theorem}[section]

\newtheorem{corollary}[theorem]{Corollary}

\newtheorem{remark}[theorem]{Remark}

\newtheorem{lemma}[theorem]{Lemma}

\newtheorem{claim}[theorem]{Claim}

\begin{document}


\title{Corrigendum:\\ The Conley conjecture for
Hamiltonian systems\\
 on the cotangent bundle\\
  and its analogue for
Lagrangian systems }

\author{\\Guangcun Lu
\thanks{Partially supported by the NNSF  10671017 and 10971014 of
China, and PCSIRT and Research Fund for the Doctoral Program Higher
Education of China (Grant No. 200800270003).}\\
{\normalsize School of Mathematical Sciences, Beijing Normal University},\\
{\normalsize Laboratory of Mathematics
 and Complex Systems,  Ministry of
  Education},\\
  {\normalsize Beijing 100875, The People's Republic
 of China}\\
{\normalsize (gclu@bnu.edu.cn)}}
\date{February 9,  2011} \maketitle \vspace{-0.1in}


\abstract{In lines 8-11 of \cite[pp. 2977]{Lu} we wrote: ``For
integer $m\ge 3$, if $M$ is $C^m$-smooth and $C^{m-1}$-smooth
$L:\R\times TM\to\R$ satisfies the assumptions (L1)-(L3), then the
functional ${\cal L}_\tau$ is $C^2$-smooth, bounded below, satisfies
the Palais-Smale condition, and all critical points of it have
finite Morse indexes and nullities (see \cite[Prop.4.1, 4.2]{AbF}
and \cite{Be}).'' However, as proved in \cite{AbSc1} the claim that
${\cal L}_\tau$ is $C^2$-smooth is true if and only if for every
$(t,q)$ the function $v\mapsto L(t,q,v)$ is a polynomial of degree
at most $2$. So the arguments in \cite{Lu} is only valid for the
physical Hamiltonian in (1.2) and corresponding Lagrangian therein.
 In this note we shall
correct our arguments in \cite{Lu} with a new splitting lemma
obtained in \cite{Lu2}. } \vspace{-0.1in}
\medskip\vspace{12mm}

 \vspace{2mm}

\section{A splitting lemma for $C^1$-functionals}\label{sec:1}
\setcounter{equation}{0}

In this section we shall give a special version of the splitting
lemma obtained by the author in \cite[Th. 2.1]{Lu2} recently. For
completeness we shall outline its proof because it is much simpler
than general case. The reader may refer to \cite{Lu2} for details.

Let $H$ be a Hilbert space with inner product $(\cdot,\cdot)_H$ and
the induced norm $\|\cdot\|$, and let $X$ be a Banach space with
norm $\|\cdot\|_X$, such that
\begin{description}
\item[(S)]  $X\subset H$ is dense in $H$ and
 $\|x\|\le \|x\|_X\;\forall x\in X$.
\end{description}
For an open neighborhood $V$ of the origin $\theta\in H$, $V\cap X$
is also an open neighborhood of $\theta$ in $X$, and we shall write
$V\cap X$ as $V_X$ when viewed as an open neighborhood of $\theta$
in $X$. For a $C^1$ functional ${\cal L}:V\to\mathbb{R}$ with
$\theta$ as an isolated critical point, suppose that there exist
maps $A\in C^1(V_X, X)$ and $B\in C(V_X, L_s(H))$ such that
\begin{eqnarray}
&&{\cal L}'(x)(u)=(A(x), u)_H\quad\forall x\in V_X\;\hbox{and}\;
u\in X,\label{e:1.1}\\
&&(A'(x)(u), v)_H=(B(x)u, v)_H\;\forall x\in V_X\;\hbox{and}\; u,
v\in X. \label{e:1.2}
\end{eqnarray}
(These imply: (a) ${\cal L}|_{V_X}\in C^2(V_X, \R)$,  (c)  $d^2{\cal
L}|_{V_X}(x)(u,v)=(B(x)u,v)_H$ for any $x\in V_X$ and $u, v\in X$,
(c) $B(x)(X)\subset X\;\forall x\in V_X$). Furthermore we also
assume $B$ to satisfy  the following properties:
\begin{description}
\item[(B1)]  If $u\in H$ such that $B(\theta)(u)=v$ for some
$v\in X$, then $u\in X$. Moreover,  all eigenfunctions of the
operator $B(\theta)$ that correspond to negative eigenvalues belong
to $X$.

\item[(B2)] The map $B:V_X\to
L_s(H,H)$  has a decomposition
$$
B(x)=P(x)+ Q(x)\quad\forall x\in V\cap X,
$$
where $P(x):H\to H$ is a positive definitive linear operator and
$Q(x):H\to H$ is a compact linear  operator with the following
properties:
\begin{description}
\item[(i)] For any sequence $\{x_k\}\subset
V\cap X$ with $\|x_k\|\to 0$ it holds that
$\|P(x_k)u-P(\theta)u\|\to 0$ for any $u\in H$;

\item[(ii)] The  map $Q:V\cap X\to
L(H,H)$ is continuous at $\theta$ with respect to the topology
induced from $H$ on $V\cap X$;

\item[(iii)] There exist positive constants $\eta_0>0$ and  $C_0>0$ such that
$$
(P(x)u, u)\ge C_0\|u\|^2\quad\forall u\in H,\;\forall x\in
B_H(\theta,\eta_0)\cap X.
$$
\end{description}
\end{description}
{\it Note}: since $B(\theta)\in L_s(H)$ is a self-adjoint Fredholm
operator, either $0\notin \sigma(B(\theta))$ or $0$ is an isolated
point in $\sigma(B(\theta))$ which is also an eigenvalue of finite
multiplicity. See Proposition B.2 in Appendix of \cite{Lu2}.

  Let $H^0:={\rm Ker}(B(\theta))$ and let $H^-$ (resp. $H^+$) be the positive
subspace (resp. negative definite) of $B(\theta)$. They are all
invariant subspaces of $B(\theta)$, and there exists an orthogonal
decomposition $H=H^0\oplus H^\pm=H^0\oplus H^-\oplus H^+$. Clearly,
\begin{equation}\label{e:1.3}
\left.\begin{array}{ll}
&(B(\theta)u, v)_H=0\;\forall u\in H^+\oplus H^-,\;v\in H^0,\\
&(B(\theta)u, v)_H=0\;\forall  u\in H^-\oplus H^0,\;v\in H^+,\\
&(B(\theta)u, v)_H=0\;\forall u\in H^+\oplus H^0,\;v\in H^-.
\end{array}\right\}
\end{equation}
Moreover, the conditions (B1) and (B2) imply that both $H^0$ and
$H^-$ are finitely dimensional subspaces contained in $X$, and that
there exists a small $a_0>0$ such that $[-2a_0,
2a_0]\cap\sigma(B(\theta))$ at most contains a point $0$. Hence
\begin{equation}\label{e:1.4}
\left.\begin{array}{ll}
  (B(\theta)u, u)_H\ge
2a_0\|u\|^2\quad\forall u\in H^+,\\
  (B(\theta)u, u)_H\le
-2a_0\|u\|^2\quad\forall u\in H^-.
\end{array}\right\}
\end{equation}
 Note that $H^\pm:=H^+ + H^-$ is the image of $B(\theta)$. Denote by $P^\ast$ the orthogonal
projections onto $H^\ast$, $\ast=+, -, 0$, and by $X^\ast=X\cap
H^\ast=P^\ast(X),\;\ast=+, -$. Then $X^+$ is dense in $H^+$, and
$(I-P^0)|_X=(P^++P^-)|_X: (X, \|\cdot\|_X)\to (X^\pm, \|\cdot\|)$ is
also continuous because all norms are equivalent on a linear space
of finite dimension, where $X^\pm:=X\cap (I-P^0)(H)=X\cap H^\pm=X^-+
P^+(X)=X^-+ H^+\cap X$.  These give the following topological direct
sum decomposition:
$$
X=H^0\oplus X^\pm=H^0\oplus X^+\oplus X^-.
$$
Let $m^0=\dim H^0$ and $m^-=\dim H^-$. They are called the {\it
nullity} and the {\it Morse index} of critical point $\theta$ of
${\cal L}$, respectively. The critical point $\theta$ is said to be
{\it nondegenerate} if $m^0=0$. For a normed vector space $(H,
\|\cdot\|)$ and $\delta>0$ let $B_H(\theta, \delta)=\{x\in
H\,|\,\|x\|<\delta\}$ and $\bar B_H(\theta, \delta)=\{x\in
H\,|\,\|x\|\le\delta\}$. Since the norms $\|\cdot\|$ and
$\|\cdot\|_X$ are equivalent on the finite dimension space $H^0$ we
shall not point out the norm used without occurring of confusions.

\begin{theorem}\label{th:1.1}
Under the above assumptions (S) and (B1)-(B2), there exist a
positive $\epsilon\in\R$,  a  $C^1$ map
$h:B_{H^0}(\theta,\epsilon)=B_{H}(\theta,\epsilon)\cap H^0\to X^\pm$
satisfying $h(\theta^0)=\theta^\pm$ and
\begin{equation}\label{e:1.5}
 (I-P^0)A(z+ h(z))=0\quad\forall z\in B_{H^0}(\theta,\epsilon),
 \end{equation}
an open neighborhood $W$ of $\theta$ in $H$ and an origin-preserving
homeomorphism
\begin{equation}\label{e:1.6}
\Phi: B_{H^0}(\theta,\epsilon)\times
\left(B_{H^+}(\theta^+,\epsilon) +
B_{H^-}(\theta^-,\epsilon)\right)\to W
\end{equation}
of form $\Phi(z, u^++ u^-)=z+ h(z)+\phi_z(u^++ u^-)$ with
$\phi_z(u^++ u^-)\in H^\pm$  such that
\begin{equation}\label{e:1.7}
{\cal L}\circ\Phi(z, u^++ u^-)=\|u^+\|^2-\|u^-\|^2+ {\cal L}(z+
h(z))
\end{equation}
for all $(z, u^+ + u^-)\in B_{H^0}(\theta,\epsilon)\times
\left(B_{H^+}(\theta^+,\epsilon) +
B_{H^-}(\theta^-,\epsilon)\right)$, and that
\begin{equation}\label{e:1.8}
\Phi\left(B_{H^0}(\theta,\epsilon)\times
\bigl(B_{H^+}(\theta^+,\epsilon)\cap X +
B_{H^-}(\theta^-,\epsilon)\bigr)\right)\subset X.
\end{equation}
Moreover, the maps $\Phi$, $h$ and the function
$B_{H^0}(\theta,\epsilon)\ni z\mapsto {\cal L}^\circ(z):={\cal L}(z+
h(z))$ also satisfy:
\begin{description}
\item[(i)]  For each $z\in B_{H^0}(\theta,\epsilon)$, $\Phi(z, \theta^\pm)=z+ h(z)$,
$\phi_z(u^++ u^-)\in H^-$ if and only if $u^+=\theta^+$;

\item[(ii)]  $h'(z)=-[(I-P^0)A'(z+
h(z))|_{X^\pm}]^{-1}(I-P^0)A'(z+ h(z))|_{H^0} \quad\forall z\in
B_{H^0}(\theta,\epsilon)$;

\item[(iii)]  ${\cal L}^\circ$ is $C^{2}$, $d^2{\cal
L}^\circ(\theta^0)=0$ and
 $$
d{\cal L}^\circ(z_0)(z)=(A(z_0+ h(z_0)), z)_H\quad\forall z_0\in
B_{H^0}(\theta, \epsilon),\; z\in H^0;
 $$

\item[(iv)]  $\theta^0$ is also an isolated critical point
of ${\cal L}^\circ$.
\end{description}
\end{theorem}

\begin{corollary}\label{cor:1.2}{\rm (Shifting)}
Under the assumptions of Theorem~\ref{th:1.1},  if $\theta$ is an
isolated critical point  of ${\cal L}$, for any Abelian group ${\bf
K}$ it holds that
$$
C_q({\cal L}, \theta;{\bf K})\cong C_{q-m^-}({\cal L}^{\circ},
\theta^0; {\bf K})\quad\forall q=0, 1,\cdots,
$$
where  ${\cal L}^{\circ}(z)={\cal L}(h(z)+z)$. Consequently,
$C_q({\cal L}, \theta;{\bf K})=0$ for $q\notin [m^-, m^-+ m^0]$, and
$C_q({\cal L}, \theta;{\bf K})$ is isomorphic to a finite direct sum
$r_1{\bf K}\oplus\cdots\oplus r_s{\bf K}$ for each $q\in [m^-, m^-+
m^0]$, where each $r_j\in\{0,1\}$. Moreover, if $C_q({\cal L},
\theta;{\bf K})\ne 0$ for some $q\in (m^-, m^-+ m^0]$, then $\theta$
must be a non-minimal saddle point of ${\cal L}$.
\end{corollary}

The proof of the final claim is as follows. If $m^0=0$ then the
conclusion is a consequence of (\ref{e:1.7}). If $m^0>0$  it follows
 from \cite[Ex.1, pp.33]{Ch} that $\theta^0$ is a non-minimal
saddle point of ${\cal L}^\circ$. This implies that $\theta$ is a
non-minimal saddle point of ${\cal L}$.

Our proof needs the following parameterized version of the
Morse-Palais lemma due to Duc-Hung-Khai (Theorem 1.1 in \cite{DHK}),
whose proof can be obtained by almost repeating the proof in
\cite{DHK} (cf. Appendix A of \cite{Lu2} for details).

\begin{theorem}\label{th:1.3}
Let $(H, \|\cdot\|)$ be a normed vector space and let $\Lambda$ be a
 compact topological space. Let $J:\Lambda\times B_H(\theta,
2\delta)\to\R$ be continuous, and for every $\lambda\in\Lambda$ the
function $J(\lambda, \cdot): B_H(\theta, 2\delta)\to\R$ is
continuously directional differentiable.
 Assume that there exist a closed
vector subspace $H^+$ and a finite-dimensional vector subspace $H^-$
of $H$ such that $H^+\oplus H^-$ is a direct sum decomposition of
$H$ and
\begin{description}
\item[(i)] $J(\lambda, \theta)=0$ and $D_2J(\lambda, \theta)=0$,
\item[(ii)] $[D_2J(\lambda, x+ y_2)-D_2J(\lambda, x+ y_1)](y_2-y_1)<0$ for any $(\lambda, x)\in\Lambda\times\bar
B_{H^+}(\theta^+,\delta)$, $y_1, y_2\in\bar
B_{H^-}(\theta^-,\delta)$ and $y_1\ne y_2$,

\item[(iii)] $D_2J(\lambda, x+y)(x-y)>0$ for any $(\lambda, x, y)\in\Lambda\times\bar B_{H^+}(\theta^+,
\delta)\times\bar B_{H^-}(\theta^-,\delta)$ and $(x,y)\ne (\theta^+,
\theta^-)$,

\item[(iv)] $D_2J(\lambda, x)x>p(\|x\|)$ for any $(\lambda, x)\in\Lambda\times\bar
B_{H^+}(\theta^+,\delta)\setminus\{\theta^+\}$, where $p:(0,
\delta]\to (0, \infty)$ is a non-decreasing function.
\end{description}
Then there exist a positive $\epsilon\in\R$, an open neighborhood
$U$ of $\Lambda\times\{\theta\}$ in $\Lambda\times H$ and a
homeomorphism
$$
\phi: \Lambda\times \bigl(B_{H^+}(\theta^+, \sqrt{p(\epsilon)/2})+
B_{H^-}(\theta^-, \sqrt{p(\epsilon)/2})\bigr)\to U
$$
such that
$$
J(\lambda, \phi(\lambda, x+ y))=\|x\|^2-\|y\|^2\quad\hbox{and}\quad
\phi(\lambda, x+ y)=(\lambda, \phi_\lambda(x+y))\in\Lambda\times H
$$
for all $(\lambda, x,y)\in \Lambda\times B_{H^+}(\theta^+,
\sqrt{p(\epsilon)/2})\times B_{H^-}(\theta^-,
\sqrt{p(\epsilon)/2})$. Moreover, for each $\lambda\in\Lambda$,
$\phi_\lambda(0)=0$, $\phi_\lambda(x+y)\in H^-$ if and only if
$x=0$.
\end{theorem}

\noindent{\bf Proof of Theorem~\ref{th:1.1}}.\quad \textsf{Step 1}.
As noted below the condition ({\bf B2}), either $0\notin
\sigma(B(\theta))$ or $0$ is an isolated point in
$\sigma(B(\theta))$. Using this, and $A'(\theta)=B(\theta)|_X$ and
the condition ({\bf B1}) it was proved in \cite{JM} that
$B(\theta)(X^\pm)\subset X^\pm$ and $B(\theta)|_{X^\pm}: X^\pm\to
X^\pm$ is a Banach space isomorphism. Since $A\in C^1(V\cap X, X)$,
we can directly apply the implicit function theorem
\cite[Th.3.7.2]{Schi} to $C^1$-map
$$
T:(H^0\cap V)\times (X^\pm\cap V)\to X^\pm,\; (z,x)\mapsto
(I-P^0)A(z+ x),
$$
and get $\delta>0$,  a (unique) $C^1$ map
$$
h:B_{H^0}(\theta, 2\delta)=B_{H}(\theta, 2\delta)\cap H^0\subset
V\cap X\to X^\pm
$$
satisfying $h(\theta^0)=\theta^\pm$ and (\ref{e:1.5}), i.e.
$(I-P^0)A(z+ h(z))=0$ for all $z\in B_{H^0}(\theta, 2\delta)$.
Moreover, the standard arguments show that the map $h$ and the
function ${\cal L}^\circ: B_{H^0}(\theta, 2\delta)\to\R$ given by
${\cal L}^\circ(z)={\cal L}(z+ h(z))$
 satisfy the conclusions (ii)-(iv) in Theorem~\ref{th:1.1}.
Let us shrink $\delta>0$ (if necessary) so that
\begin{equation}\label{e:1.9}
z+ h(z)+ u\in V\quad\forall (z,u)\in (\bar B_H(\theta, \delta)\cap
H^0)\times (\bar B_H(\theta, \delta)\cap H^\pm).
\end{equation}

 Define a $C^1$ functional $F:\bar B_{H^0}(\theta, \delta)\times B_{H^\pm}(\theta, \delta)\to\R$ as
\begin{equation}\label{e:1.10}
 F(z, u)={\cal L}(z+ h(z)+ u)-{\cal L}(z+ h(z)).
\end{equation}
 Then for each
$(z,u)\in \bar B_{H^0}(\theta, \delta)\times B_{H^\pm}(\theta,
\delta)$ and  $v\in H^\pm$ it holds that
\begin{eqnarray}\label{e:1.11}
D_2F(z,u)(v)&=&(A(z+ h(z)+ u), v)_H\nonumber\\
&=&((I-P^0)A(z+ h(z)+ u), v)_H.
\end{eqnarray}
 It follows from this and (\ref{e:1.5}) that
\begin{eqnarray}\label{e:1.12}
F(z, \theta^\pm)=0\quad\hbox{and}\quad D_2F(z,
\theta^\pm)(v)=0\;\forall v\in H^\pm.
\end{eqnarray}

In next step we shall show that Theorem~\ref{th:1.3} can be applied
to the functional $F$.

\textsf{Step 2}. {\it Claim 1}.
 There exists a function $\omega:V\cap X\to [0, \infty)$
 such that $\omega(x)\to 0$ as $x\in V\cap X$ and $\|x\|\to
0$, and that
$$
|(B(x)u, v)_H- (B(\theta)u, v)_H |\le \omega(x) \|u\|\cdot\|v\|
$$
for any $x\in V\cap X$,  $u\in H^0\oplus H^-$ and $v\in H$.

This is Lemma 2.15 in \cite{Lu2}.  Firstly, by a contradiction
argument the condition (i) of (B2) can be equivalently expressed as:
{\it For any $u\in H$ it holds that $\|P(x)u-P(\theta)u\|\to 0$ as
$x\in V\cap X$ and $\|x\|\to 0$.}

Next  let $e_1,\cdots, e_m$ be a basis of $H^0\oplus H^-$ with
$\|e_i\|=1$, $i=1,\cdots, m$. Then
 $$
\left(\sum^m_{i=1}|t_i|^2\right)^{1/2}\le C_1\|u\|
$$
for some constant $C_1>0$ and $u=\sum^m_{i=1}t_ie_i\in H^0\oplus
H^-$.  Since
\begin{eqnarray*}
|(B(x)e_i, v)_H- (B(\theta)e_i, v)_H | \le \|(P(x)e_i-
P(\theta)e_i\|\cdot\|v\|+ \|Q(x)-Q(\theta)\|\cdot\|v\|,
\end{eqnarray*}
 for any $u=\sum^m_{i=1}t_ie_i\in H^0\oplus H^-$ we have
\begin{eqnarray*}
&&|(B(x)u, v)_H- (B(\theta)u, v)_H | \\
&\le &\sum^m_{i=1}|t_i|\|P(x)e_i- P(\theta)e_i\|\cdot\|v\|+
\sum^m_{i=1}|t_i|\|Q(x)-Q(\theta)\|\cdot\|v\|\\
&\le &\left(\sum^m_{i=1}\|P(x)e_i-
P(\theta)e_i\|^2\right)^{1/2}\left(\sum^m_{i=1}|t_i|^2\right)^{1/2}\|v\|\\
&&\quad +
\sqrt{m}\left(\sum^m_{i=1}|t_i|^2\right)^{1/2}\|Q(x)-Q(\theta)\|\cdot\|v\|\\
&\le &\left[C_1
\biggl(\sum^m_{i=1}\|P(x)e_i-P(\theta)e_i\|^2\biggr)^{1/2}
+ C_1\sqrt{m}\|Q(x)-Q(\theta)\| \right]\|u\|\|v\|\\
&=&\omega(x)\|u\|\|v\|,
\end{eqnarray*}
 where
$$
\omega(x)=\left[C_1 \biggl(\sum^m_{i=1}\|P(x)e_i-
P(\theta)e_i\|^2\biggr)^{1/2} + C_1\sqrt{m}\|Q(x)-Q(\theta)\|
\right]\to 0
$$
as $x\in V\cap X$ and $\|x\|\to 0$ (because of the conditions (i)
and (ii) in ({\bf B2})). $\Box$\vspace{2mm}

As in the proof of Lemma 2 in \cite[page 201]{Skr} (see also Lemma
5.2 of \cite{Va}) we can prove:

{\it Claim 2}. There exists a  small neighborhood $U\subset V$ of
$\theta$ in $H$ and a number $a_1\in (0, 2a_0]$ such that for any
$x\in U\cap X$,
\begin{description}
\item[(i)] $(B(x)u, u)_H\ge a_1\|u\|^2\;\forall u\in H^+$;
\item[(ii)] $|(B(x)u,v)_H|\le\omega(x)\|u\|\cdot\|v\|\;\forall u\in H^+, \forall v\in
H^-\oplus H^0$;
\item[(iii)] $(B(x)u,u)_H\le-a_0\|u\|^2\;\forall u\in H^-$.
\end{description}

The reader may refer to Lemma~2.16 in \cite{Lu2} for a detailed
proof of it.

Since $h(\theta^0)=\theta^\pm$,  we may take $\varepsilon\in (0,
\delta)$ so small that
 \begin{equation}\label{e:1.13}
 z+ h(z)+
u^++ u^-\in U
\end{equation}
for all $z\in\bar B_{H^0}(\theta,\varepsilon)$, $u^+\in \bar
B_{H^+}(\theta,\varepsilon)$ and $u^-\in\bar
B_{H^-}(\theta,\varepsilon)$.

\textsf{Step 3}. The restriction of the function $F$ in
(\ref{e:1.10}) to $\bar B_{H^0}(\theta,\varepsilon)\times \bigl(\bar
B_{H^+}(\theta,\varepsilon)\oplus \bar
B_{H^-}(\theta,\varepsilon)\bigr)$ satisfies the conditions in
Theorem~\ref{th:1.3}.

This is Lemma~2.17 in \cite{Lu2}.
 By (\ref{e:1.12}) we only need to prove that
$F$ satisfies conditions (ii)-(iv) in Theorem~\ref{th:1.3}.

For $z\in\bar B_{H^0}(\theta,\varepsilon)$, $u^+\in \bar
B_H(\theta,\varepsilon)\cap X^+$ and $u^-_1, u^-_2\in\bar
B_{H^-}(\theta,\varepsilon)$,  (\ref{e:1.1}) gives
\begin{eqnarray*}
&&[D_2F(z, u^+ + u^-_2)-D_2F(z, u^++ u^-_1)](u^-_2-u^-_1)\\
&=&(A(z+ h(z)+ u^++u^-_2), u^-_2-u^-_1)_H - (A(z+ h(z)+ u^++u^-_1),
u^-_2-u^-_1)_H.
\end{eqnarray*}
By the mean value theorem we have $t\in (0, 1)$ such that
\begin{eqnarray*}
&&(A(z+ h(z)+ u^++u^-_2), u^-_2-u^-_1)_H - (A(z+ h(z)+ u^++u^-_1), u^-_2-u^-_1)_H\\
&=&\left(DA(z+ h(z)+ u^++ u^-_1+ t(u^-_2-u^-_1), u^-_2-u^-_1),
u^-_2-u^-_1\right)_H\\
&\stackrel{(\ref{e:1.2})}{=}&\left(B(z+ h(z)+ u^++ u^-_1+
t(u^-_2-u^-_1))(u^-_2-u^-_1),
u^-_2-u^-_1\right)_H\\
&\le& -a_0\|u^-_2-u^-_1\|^2,
\end{eqnarray*}
where the final inequality comes from (iii) of Claim 2 in Step 2.
Hence
\begin{eqnarray*}
[D_2F(z, u^+ + u^-_2)-D_2F(z, u^++ u^-_1)](u^-_2-u^-_1)\le
-a_0\|u^-_2-u^-_1\|^2.
\end{eqnarray*}
 Since $\bar
B_H(\theta,\varepsilon)\cap X^+$ is dense in $\bar
B_H(\theta,\varepsilon)\cap H^+$ we get
\begin{eqnarray}\label{e:1.14}
[D_2F(z, u^+ + u^-_2)-D_2F(z, u^++ u^-_1)](u^-_2-u^-_1)\le
-a_0\|u^-_2-u^-_1\|^2.
\end{eqnarray}
for all $z\in\bar B_{H^0}(\theta,\varepsilon)$, $u^+\in \bar
B_H(\theta,\varepsilon)\cap H^+$ and $u^-\in\bar
B_{H}(\theta,\varepsilon)\cap H^-$. It shows that the condition (ii)
in Theorem~\ref{th:1.3} holds for $F$.

Next, for $z\in\bar B_{H^0}(\theta,\varepsilon)$, $u^+\in \bar
B_H(\theta,\varepsilon)\cap X^+$ and $u^-\in\bar
B_{H^-}(\theta,\varepsilon)$, using (\ref{e:1.12}), the mean value
theorem and  (\ref{e:1.1})-(\ref{e:1.2}), we may find a $t\in (0,
1)$ such that
\begin{eqnarray*}
&&D_2F(z, u^++u^-)(u^+-u^-)\\
&=&D_2F(z, u^++u^-)(u^+-u^-)- D_2F(z, \theta^\pm)(u^+-u^-)\\
&=&(A(z+ h(z)+ u^++u^-), u^+-u^-)_H-(A(z+ h(z)+ \theta^\pm), u^+-u^-)_H\\
&=&\left(B(z+ h(z)+ t(u^++u^-))(u^++u^-), u^+-u^-\right)_H\\
&=&\left(B(z+ h(z)+ t(u^++u^-))u^+, u^+\right)_H-\left(B(z+ h(z)+
t(u^++u^-))u^-, u^-\right)_H\\
&\ge & a_1\|u^+\|^2+ a_0\|u^-\|^2.
\end{eqnarray*}
Here the final inequality is due to (i) and (iii) in Claim 2.  As
above this inequality also holds for all $u^+\in \bar
B_{H^+}(\theta,\varepsilon)$ because $\bar
B_H(\theta,\varepsilon)\cap X^+$ is dense in $\bar
B_H(\theta,\varepsilon)\cap H^+$. It is more than zero when $(u^+,
u^-)\ne (\theta^+, \theta^-)$. Hence the condition (iii) of
Theorem~\ref{th:1.3} is satisfied.

Finally, for $z\in\bar B_{H^0}(\theta,\varepsilon)$ and $u^+\in \bar
B_H(\theta,\varepsilon)\cap X^+$, as above we have $t\in (0, 1)$
such that
\begin{eqnarray*}
D_2F(z, u^+)u^+
&=&D_2F(z, u^+)u^+- D_2F(z, \theta^\pm)u^+\\
&=&(A(z+ h(z)+ u^+), u^+)_H-(A(z+ h(z)+ \theta^\pm), u^+)_H\\
&=&\left(B(z+ h(z)+ tu^+)u^+, u^+\right)_H\\
&\ge& a_1\|u^+\|^2
\end{eqnarray*}
because of Claim 2(i). So for the function $p:(0, \varepsilon]\to
(0, \infty),\;t\mapsto\frac{a_1}{2}t^2$ it holds that
$$
D_2F(z, u^+)u^+ \ge a_1\|u^+\|^2> p(\|u^+\|)\quad\forall u^+\in \bar
B_H(\theta,\varepsilon)\cap H^+\setminus\{\theta^+\}.
$$
This shows that $F$ satisfies the condition (iv) in
Theorem~\ref{th:1.3}.

\textsf{Step 4}.
 Applying Theorem~\ref{th:1.3} to $F$  we can
 get a positive number $\epsilon$, an open neighborhood
 ${\cal W}$ of $\bar B_{H^0}(\theta,\varepsilon)\times\{\theta^\pm\}$ in $\bar B_{H^0}(\theta,\varepsilon)\times H^\pm$,
  and  an origin-preserving homeomorphism
\begin{eqnarray}\label{e:1.15}
\phi: \bar B_{H^0}(\theta,\varepsilon)\times \left(B_{H^+}(\theta^+,
\epsilon)+ B_{H^-}(\theta^-, \epsilon)\right)\to {\cal W}
\end{eqnarray}
of form $\phi(z, u^+ + u^-)=(z, \phi_z(u^+ + u^-))\in (\bar
B_{H^0}(\theta,\varepsilon)\times H^\pm$ such that $\phi_z(\theta^+
+ \theta^-))=\theta^\pm$ and
\begin{eqnarray}\label{e:1.16}
 &&{\cal L}(z+ h(z)+ \phi_z(u^+, u^-))-{\cal L}(z+ h(z))\nonumber\\
 &=&F(\phi(z, u^+, u^-))=\|u^+\|^2-\|u^-\|^2
\end{eqnarray}
for all $(z, u^+, u^-)\in \bar B_{H^0}(\theta,\varepsilon)\times
B_{H^+}(\theta^+, \epsilon)\times B_{H^-}(\theta^-, \epsilon)$.
Moreover, $\phi_z(u^++u^-)\in H^-$ if and only if $u^+=\theta^+$.

Consider the continuous map
\begin{eqnarray}\label{e:1.17}
&&\Phi:B_{H^0}(\theta,\varepsilon)\times \left(B_{H^+}(\theta^+,
\epsilon)+ B_{H^-}(\theta^-, \epsilon)\right)\to H,\\
&&\hspace{10mm} (z, u^+ + u^-)\mapsto z+ h(z)+ \phi_z(u^+ +
u^-).\nonumber
\end{eqnarray}
Then (\ref{e:1.16}) gives (\ref{e:1.7}), i.e. ${\cal L}(\Phi(z, u^+,
u^-))=\|u^+\|^2-\|u^-\|^2+ {\cal L}(z+ h(z))$.

{\it Claim 3}.  $W:={\rm Im}(\Phi)$ is an open neighborhood of
$\theta$ in $H$ and $\Phi$  is an origin-preserving homeomorphism
onto $W$.

In fact, assume that $\Phi(z_1, u^+_1+ u^-_1)=\Phi(z_2, u^+_2+
u^-_2)$ for $(z_1, u^+_1+ u^-_1)$ and $(z_2, u^+_2+ u^-_2)$ in
$B_{H^0}(\theta,\varepsilon)\times (B_{H^+}(\theta^+, \epsilon)+
B_{H^-}(\theta^-, \epsilon))$. Then
$$
z_1=z_2\quad\hbox{and}\quad h(z_1)+ \phi_{z_1}(u^+_1+ u^-_1)=h(z_2)+
\phi_{z_2}(u^+_2+ u^-_2).
$$
It follows that $h(z_1)=h(z_2)$ and $\phi_{z_1}(u^+_1+ u^-_1)=
\phi_{z_2}(u^+_2+ u^-_2)$. They show that $\Phi(z_1, u^+_1+ u^-_1)=
\Phi(z_2, u^+_2+ u^-_2)$ and thus $(u^+_1, u^-_1)=(u^+_2, u^-_2)$.
So $\Phi$ is a bijection.

For a point $(z, u^++ u^-)$ and a sequence $\{(z_k, u^+_k+ u^-_k)\}$
 in $B_{H^0}(\theta,\varepsilon)\times (B_{H^+}(\theta^+,
\epsilon)+ B_{H^-}(\theta^-, \epsilon))$, suppose that $\Phi(z_k,
u^+_k+ u^-_k)\to\Phi(z, u^++ u^-)$. Then
\begin{eqnarray*}
&&P^0\Phi(z_k, u^+_k+ u^-_k)\to P^0\Phi(z, u^++
u^-)\quad\hbox{and}\\
&& (P^++P^-)\Phi(z_k, u^+_k+ u^-_k)\to(P^++P^-)\Phi(z, u^++ u^-).
\end{eqnarray*}
It follows that $z_k\to z$, and thus $h(z_k)\to h(z)$ and
$\phi_{z_k}(u^+_k+ u^-_k)\to \phi_z(u^++ u^-)$. These imply that
$\phi(z_k, u^+_k+ u^-_k)\to\phi(z, u^++ u^-)$ and hence $(z_k,
u^+_k+ u^-_k)\to (z, u^++ u^-)$ since $\phi$ is a homeomorphism.
That is, $\Phi^{-1}$ is continuous. Hence $\Phi$ is
 a homeomorphism onto $W$, and so $W$ is open in $H$.
The proof of Theorem~\ref{th:1.1} is completed. $\Box$\vspace{2mm}

Consider a tuple $(H, X, {\cal L}, A, B=P+ Q)$, where $H$ (resp.
$X$) is a Hilbert (resp. Banach) space satisfying the condition
({\bf S}) as in Section~\ref{sec:1}, the functional ${\cal
L}:H\to\R$ and maps $A:X\to H$ and $B:X\to L_s(H, H)$ satisfy, at
least near the origin $\theta\in H$, the conditions ({\bf B1})-({\bf
B2}) in Section~\ref{sec:1}. Let $(\widehat H, \widehat X,
\widehat{\cal L}, \widehat A, \widehat B=\widehat P+ \widehat Q)$ be
another such a tuple. Suppose that $J:H\to\widehat H$ is a linear
injection such that $J(X)\subset X$ and
\begin{equation}\label{e:1.18}
(Ju, Jv)_{\widehat H}=(u, v)_H\quad\hbox{and}\quad \|Jx\|_{\widehat
X}=\|x\|_X
\end{equation}
for all $u, v\in H$ and $x\in X$. Furthermore, we assume
\begin{equation}\label{e:1.19}
\widehat{\cal L}\circ J={\cal L}\quad\hbox{and}\quad \widehat
P(J(x))\circ J=J\circ P(x)\;\forall x\in X.
\end{equation}
Then we have
\begin{equation}\label{e:1.20}
\left.\begin{array}{ll}
 & \widehat A(J(x))=J\circ
A(x),\quad\widehat B(J(x))\circ J=J\circ B(x)\;\forall x\in X, \\
& \hbox{and thus}\quad \widehat Q(J(x))\circ J=J\circ Q(x)\;\forall
x\in X.
\end{array}\right\}
\end{equation}
Let $H=H^0\oplus H^+\oplus H^-$, $X=H^0\oplus X^+\oplus X^-$ and
$\widehat H=\widehat H^0\oplus \widehat H^+\oplus \widehat H^-$ and
$\widehat X=\widehat H^0\oplus \widehat X^+\oplus \widehat X^-$ be
the corresponding decompositions. Namely, $\widehat H^0={\rm
Ker}(\widehat B(\theta))$, and $\widehat H^+$ (resp. $\widehat H^-$)
is the positive (resp. negative) definite subspace of $\widehat
B(\theta)$. Denote by $P^\ast$ (resp. $\widehat P^\ast$) the
orthogonal projections from $H$ (resp. $\widehat H$) to $H^\ast$
(resp. $\widehat H^\ast$) for $\ast=+, -, 0$.
 We also assume that the Morse index and nullity of ${\cal L}$ at
$\theta\in H$ are equal to those of $\widehat{\cal L}$ at
$\theta\in\widehat H$, i.e.,
\begin{equation}\label{e:1.21}
m^-({\cal L}, \theta)=m^-(\widehat{\cal L},
\theta)\quad\hbox{and}\quad m^0({\cal L}, \theta)=m^0(\widehat{\cal
L}, \theta).
\end{equation}
Since $\widehat B(\theta)\circ J=J\circ B(\theta)$ by
(\ref{e:1.20}), (\ref{e:1.21}) implies
\begin{equation}\label{e:1.22}
\left.\begin{array}{ll}
 & JH^0=\widehat H^0, \quad  \widehat P^0\circ J=J\circ P^0, \\
& JH^-={\widehat H}^-, \quad  {\widehat P}^-\circ J=J\circ P^-, \\
& JH^+\subset {\widehat H}^+, \quad  \widehat P^+\circ J=J\circ P^+.
\end{array}\right\}
\end{equation}
 The following  functor property of the
splitting lemma Theorem~\ref{th:1.1}  is a special version of
Theorem~2.25 in \cite{Lu2}.

\begin{theorem}\label{th:1.4}
Under the assumptions above, for the $C^1$ maps
 $h:B_{H^0}(\theta,\epsilon)\to
X^\pm$ and $\hat h:B_{\widehat H^0}(\theta,\epsilon)\to \widehat
X^\pm$, and the origin-preserving homeomorphisms constructed in
Theorem~\ref{th:1.1},
\begin{eqnarray*}
&&\Phi: B_{H^0}(\theta,\epsilon)\times
\left(B_{H^+}(\theta,\epsilon) +
B_{H^-}(\theta,\epsilon)\right)\to W,\\
&&\widehat\Phi: B_{\widehat H^0}(\theta,\epsilon)\times \bigl(
B_{\widehat H^+}(\theta,\epsilon) + B_{\widehat
H^-}(\theta,\epsilon)\bigr)\to \widehat W,
\end{eqnarray*}
it holds that
$$
\hat h(Jz)=J\circ h(z)\quad\hbox{and}\quad \widehat\Phi(Jz, Ju^+ +
Ju^-)= J\circ\Phi(z, u^+ + u^-)
$$
 for all $(z, u^+, u^-)\in B_{H^0}(\theta,\epsilon)\times
B_{H^+}(\theta,\epsilon)\times B_{H^-}(\theta,\epsilon)$.
Consequently,
\begin{eqnarray*}
&&\widehat{\cal L}\circ\widehat\Phi(Jz, Ju^++ Ju^-)={\cal
L}\circ\Phi(z, u^++ u^-),\\
&&\widehat{\cal L}(Jz+ \hat h(Jz))={\cal L}(z+ h(z))
\end{eqnarray*}
for all $(z, u^+, u^-)\in B_{H^0}(\theta,\epsilon)\times
B_{H^+}(\theta,\epsilon) \times B_{H^-}(\theta,\epsilon)$.
\end{theorem}

By Step 1 of the proof of Theorem~\ref{th:1.1} and
(\ref{e:1.18})-(\ref{e:1.22}) one easily concludes $\hat
h(Jz)=J\circ h(z)$ for any $z\in B_{H^0}(\theta,\epsilon)$.
Carefully checking the proof of Theorem~\ref{th:1.1} it is not hard
to derive the other conclusions. See the proof of Theorem~2.25 in
\cite{Lu2}.

We actually need a variant of Theorem~\ref{th:1.4} above. For $1\le
r<\infty$ suppose that the first relations in (\ref{e:1.18}) and
(\ref{e:1.19}) are replaced by the following
\begin{equation}\label{e:1.23}
(Ju, Jv)_{\widehat H}=r(u, v)_H\quad\hbox{and}\quad \widehat{\cal
L}\circ J=r{\cal L}
\end{equation}
 for all $u, v\in H$ and $x\in X$, and other assumptions are not changed.
What are the corresponding conclusions? In order to understand this
question we define $\overline{H}$ to be the Hilbert space $\widehat
H$ equipped with an equivalent inner
$$
(u,v)_{\overline H}=\frac{1}{r}(u,v)_{\widehat H}.
$$
Note that we have still $\|u\|_{\overline H}=\|u\|_{\widehat
H}/\sqrt{r}\le\|u\|_{\widehat X}\;\forall u\in\widehat X$ since
$r\ge 1$. Namely, the condition ({\bf S}) is satisfied for the space
$\overline H$ and $\widehat X$. Set $\overline{\cal L}={\cal L}/r$.
It is easily checked that for the functional $\overline{\cal L}$ on
the Hilbert space $\overline H$ the corresponding maps $\overline A$
and $\overline B$ (given by (\ref{e:1.1})-(\ref{e:1.2})) are equal
to $\widehat A$ and $\widehat B$, respectively. Hence the conditions
of Theorem~\ref{th:1.4} hold for the tuples $(H, X, {\cal L}, A,
B=P+ Q)$ and $(\overline H, \widehat X, \overline {\cal L},
\overline A, \overline B=\overline P+ \overline Q)$. Obverse that
$B_{\overline H^\ast}(\theta,\epsilon)= B_{\widehat H^\ast}(\theta,
\sqrt{r}\epsilon)$ for $\ast=+,0,-$. By shrinking $\epsilon>0$ (if
necessary) Theorem~\ref{th:1.4} yields immediately:

\begin{corollary}\label{cor:1.5}
 Suppose for $1\le r<\infty$ that the first relations in the above
 assumptions (\ref{e:1.18}) and (\ref{e:1.19})  are changed into ones in
 (\ref{e:1.23}). Then there exist $\epsilon>0$,  the $C^1$ maps
 $h:B_{H^0}(\theta,\epsilon)\to
X^\pm$ and $\hat h:B_{\widehat H^0}(\theta, \sqrt{r}\epsilon)\to
\widehat X^\pm$, satisfying $h(\theta^0)=\theta^\pm$, $\hat
h(\theta^0)=\theta^\pm$ and
\begin{eqnarray*}
 &&(I-P^0)A(z+ h(z))=0\quad\forall z\in B_{H^0}(\theta,\epsilon),\\
&&(I-\widehat P^0)\widehat A(z+ \hat h(z))=0\quad\forall z\in
B_{\widehat H^0}(\theta,\sqrt{r}\epsilon),
 \end{eqnarray*}
and the origin-preserving homeomorphisms
\begin{eqnarray*}
&&\Phi: B_{H^0}(\theta,\epsilon)\times
\left(B_{H^+}(\theta,\epsilon) +
B_{H^-}(\theta,\epsilon)\right)\to W,\\
&&\widehat\Phi: B_{\widehat H^0}(\theta,\sqrt{r}\epsilon)\times
\bigl( B_{\widehat H^+}(\theta,\sqrt{r}\epsilon) + B_{\widehat
H^-}(\theta,\sqrt{r}\epsilon)\bigr)\to \widehat W
\end{eqnarray*}
satisfying (\ref{e:1.7}) for ${\cal L}$ and $\widehat {\cal L}$
respectively, such that
$$
\hat h(Jz)=J\circ h(z)\quad\hbox{and}\quad \widehat\Phi(Jz, Ju^+ +
Ju^-)=J\circ\Phi(z, u^+ + u^-)
$$
 for all $(z, u^+, u^-)\in B_{H^0}(\theta,\epsilon)\times
B_{H^+}(\theta,\epsilon)\times B_{H^-}(\theta,\epsilon)$.
Consequently,
\begin{eqnarray*}
&&\widehat{\cal L}\circ\widehat\Phi(Jz, Ju^++ Ju^-)= r{\cal
L}\circ\Phi(z, u^++ u^-),\\
&&\widehat{\cal L}(Jz+ \hat h(Jz))= r{\cal L}(z+ h(z))
\end{eqnarray*}
for all $(z, u^+, u^-)\in B_{H^0}(\theta,\epsilon)\times
B_{H^+}(\theta,\epsilon) \times B_{H^-}(\theta,\epsilon)$.
\end{corollary}

Write ${\cal L}\circ\Phi=\beta+ \alpha$, where $\alpha(z)={\cal
L}^\circ(z)={\cal L}(z+ h(z))$. Then $\beta$ and $\alpha$ are
$C^\infty$ and $C^2$, respectively, and the final two equalities in
Corollary~\ref{cor:1.5} imply
\begin{equation}\label{e:1.24}
\widehat\alpha\circ J=r\alpha\quad\hbox{and}\quad \widehat\beta\circ
J=r\beta.
\end{equation}

\section{An abstract theorem }\label{sec:2}
\setcounter{equation}{0}

Having the theory in Section~\ref{sec:1}, from the arguments on
critical modules under iteration maps in \cite{Lo, Lu} we may derive
the following abstract result.

\begin{theorem}\label{th:2.1}

 Let  tuples $(H_i, X_i, {\cal
L}_i, A_i, B_i=P_i+ Q_i)$ with open neighborhood $V_i=H_i$ of the
origin $\theta_i$ in $H_i$, $i=1,2$, satisfy the conditions ({\bf
S}) and ({\bf B1})-({\bf B2}) in Section~\ref{sec:1}. Suppose that
${\cal L}_i\in C^{2-0}(V_i,\R)$, satisfy the (PS) condition and
\begin{equation}\label{e:2.1}
m^-({\cal L}_1,\theta_1)=m^-({\cal
L}_2,\theta_2)\quad\hbox{and}\quad m^0({\cal
L}_1,\theta_1)=m^0({\cal L}_2,\theta_2).
\end{equation}
For some  constant $k>0$ let $J:H_1\to H_2$ be a linear injection
such that
\begin{eqnarray}
&& (Jx, J(y))_{H_2}=k\cdot (x, y)_{H_1}\;\forall x,y\in
V_1;\label{e:2.2} \\
&&J(X_1)\subset X_2\quad\hbox{and}\quad
\|Jx\|_{X_2}=\|x\|_{X_1}\;\forall x\in X;\label{e:2.3}\\
&& {\cal L}_2(Jx)=k\cdot {\cal L}_1(x)\;\forall x\in
V_1.\label{e:2.4}
\end{eqnarray}
(These imply:
\begin{equation}\label{e:2.5}
\left.\begin{array}{ll} & \nabla {\cal L}_2(Jx)=J\nabla {\cal
L}_1(x)\;\forall x\in V_1,\\
&A_2(Jx)=JA_1(x)\;\forall x\in V_1\cap X_1,\\
&B_2(Jx)\circ J=J\circ B_1(x)\;\forall x\in V_1\cap X_1
\end{array}\right\}
\end{equation}
 and thus
$J|_{H^0_1}:H^0_1\to H^0_2$ and $J|_{H^-_1}:H^-_1\to H^-_2$ are
linear isomorphisms.) Then for $c_i={\cal L}_i(\theta_i)$ and any
small $\varepsilon>0$ there exist  Gromoll-Meyer pairs of ${\cal
L}_i$ at $\theta_i\in H_1$ (with respect to the negative gradient
flows), $(W_i, W^-_i)$, which can be contained in $V_i$, such that
\begin{eqnarray}\label{e:2.6}
&&(W_1,  W^-_1)\subset\bigl({\cal L}_1^{-1}[c_1-\varepsilon, c_1+
\varepsilon], {\cal L}_1^{-1}(c_1-\varepsilon)\bigr),\nonumber
 \\
&&(W_{2}, W^-_{2})\subset\bigl({\cal L}_2^{-1}[c_2-k\varepsilon,
c_2+
k\varepsilon], {\cal L}_2^{-1}(c_2-k\varepsilon)\bigr),\nonumber\\
&& (J(W_1), J(W^-_1))\subset (W_{2}, W^-_{2}),
\end{eqnarray}
and  the induced homomorphisms
$$
J_\ast: H_\ast(W_1, W^-_1; \K)\to H_\ast(W_2, W^-_2; \K)
$$
are isomorphisms. (Actually these are true for any Gromoll-Meyer
pairs satisfying (\ref{e:2.6}), see Corollary~\ref{cor:2.8}).
\end{theorem}

\noindent{\bf Proof}.\quad Using Corollary~\ref{cor:1.5} one may
prove it as in \cite{Lo, Lu} directly. Here is a slightly different
proof with some proof ideas of \cite[Th.5.2]{Ch} partially.

By the construction of Gromoll-Meyer pairs (cf. \cite[page 49]{Ch})
we can require them to be contained in a given small neighborhood of
$\theta_i$. Hence we always assume $V_i=H_i$, $i=1,2$, below.

\noindent{\bf Step 1.} By the construction on page 49 of \cite{Ch},
we set
\begin{eqnarray*}
&&W_1:={\cal L}_1^{-1}[c_1-\varepsilon, c_1+
\varepsilon]\cap\bigl\{x\in H_1\,|\,
\lambda{\cal L}_1(x)+ \|x\|^2_{H_1}\le \mu\bigr\},\\
&&W_1^-:={\cal L}_1^{-1}(c_1-\varepsilon)\cap\bigl\{x\in H_1\,|\,
\lambda{\cal L}_1(x)+ \|x\|^2_{H_1}\le \mu\bigr\},\\
&&W_2:={\cal L}_2^{-1}[c_2- k\varepsilon, c_2+
k\varepsilon]\cap\bigl\{x\in H_2\,|\,
\lambda{\cal L}_2(x)+ \|x\|^2_{H_2}\le k\mu\bigr\},\\
&&W_2^-:={\cal L}_2^{-1}(c_2-k\varepsilon)\cap\bigl\{x\in H_2\,|\,
\lambda{\cal L}_2(x)+ \|x\|^2_{H_2}\le k\mu\bigr\},
\end{eqnarray*}
where positive numbers $\lambda, \mu, \varepsilon$ and $k\lambda,
k\mu, k\varepsilon$ are such that the conditions as in (5.13)-(5.15)
 on page 49 of \cite{Ch} hold. Then $(W_i,  W_i^-)$ are
 Gromoll-Meyer pairs of ${\cal L}_i$ at $\theta_i$, $i=1,2$, and
 satisfy (\ref{e:2.6}). We wish to prove

\begin{claim}\label{cl:2.2}
The map $J$ induces isomorphisms
$$
J_*: H_\ast\bigl(W_1, W_1^-; \K\bigr)\longrightarrow
         H_\ast\bigl(W_2, W_2^-; \K\bigr).
$$
\end{claim}

 Since the  Gromoll-Meyer pairs  $(W_i,  W_i^-)$ are with respect to the negative
gradient vector fields $-\nabla {\cal L}_i$, $i=1,2$, it follows
from the first equality in (\ref{e:2.5}) that
\begin{equation}\label{e:2.7}
J(\eta^{(1)}(t, x))= \eta^{(2)}(t, Jx)\quad\forall x\in H_1,
\end{equation}
where $\eta^{(j)}$ are the flows of  $-\nabla {\cal L}_{j}$,
$j=1,2$. Recall the proof of \cite[Th.5.2]{Ch}. Let
\begin{eqnarray*}
U_+^{(j)}= \cup_{0<t<\infty}\eta^{(j)}(t, W_j),\quad \tilde
U_+^{(j)}= \cup_{0<t<\infty}\eta^{(j)}(t, W_j^-)
\end{eqnarray*}
and $\mathfrak{F}_j$ be the continuous functions on $\tilde
U_+^{(j)}$ defined by the condition:
$$
\eta^{(j)}(\mathfrak{F}_j(x), x)\in ({\cal
L}_{j})_{c_j-j\varepsilon}\cap\tilde U_+^{(j)},\quad\hbox{but}\;
\eta^{(j)}(t, x)\notin ({\cal L}_{j})_{c_j-j\varepsilon}\cap\tilde
U_+^{(j)}\;\hbox{if}\; t<\mathfrak{F}_j(x).
$$
Then
\begin{equation}\label{e:2.8}
\sigma^{(j)}(t,x)=\eta^{(j)}(\mathfrak{F}_j(x), x)\quad t\in [0,1],
x\in \tilde U_+^{(j)}
\end{equation}
define strong deformation retracts
$$
\tilde U_+^{(j)}\to ({\cal L}_{j})_{c_j-j\varepsilon}\cap\tilde
U_+^{(j)}=({\cal L}_{j})_{c_j-j\varepsilon}\cap U_+^{(j)},
$$
and thus isomorphisms
$$
(\sigma_1^{(j)})_\ast: H_\ast\bigl(U_+^{(j)}, \tilde
U_+^{(j)};\K\bigr)\to H_\ast\bigl(({\cal L}_{j})_{c_j+
j\varepsilon}\cap U_+^{(j)}, ({\cal L}_{j})_{c_j-j\varepsilon}\cap
U_+^{(j)};\K\bigr),
$$
where $\sigma_1^{(j)}(\cdot)=\sigma^{(j)}(1,\cdot)$, $j=1,2$. By
(\ref{e:2.7}) and (\ref{e:2.8}) we have
$$
J(\sigma^{(1)}(t,x))=\sigma^{(2)}(t, Jx)\quad\forall x.
$$
This leads to the following commutative diagram:
\begin{equation}\label{e:2.9}
\begin{CD}
H_\ast\bigl(({\cal L}_{1})_{c_1+\varepsilon}\cap U_+^{(1)}, ({\cal
L}_{1})_{c_1-\varepsilon}\cap U_+^{(1)};\K\bigr)
@>(\sigma_1^{(1)})_\ast>> H_\ast\bigl(U_+^{(1)}, \tilde U_+^{(1)};\K\bigr) \\
@V J_\ast VV @VV J_\ast V \\
H_\ast\bigl(({\cal L}_{2})_{c_2+k\varepsilon}\cap U_+^{(2)}, ({\cal
L}_{2})_{c_2-k\varepsilon}\cap
U_+^{(2)};\K\bigr)@>(\sigma_1^{(2)})_\ast
>>H_\ast\bigl(U_+^{(2)}, \tilde U_+^{(2)};\K\bigr),
\end{CD}
\end{equation}
For $\delta>0$ let
\begin{eqnarray*}
 \tilde U_\delta^{(j)}= \cup_{\delta<t<\infty}\eta^{(j)}(t,
 W_j^-),\quad j=1,2.
\end{eqnarray*}
Then it follows from (\ref{e:2.7}) that $J(\tilde
U_\delta^{(1)})\subset \tilde U_\delta^{(2)}$ and
$$
J(U_+^{(1)}\setminus\tilde U_\delta^{(1)})\subset
U_+^{(2)}\setminus\tilde U_\delta^{(2)},\quad J(\tilde
U_+^{(1)}\setminus\tilde U_\delta^{(1)})\subset \tilde
U_+^{(2)}\setminus\tilde U_\delta^{(2)}.
$$
Hence we may get  the following commutative diagram:
\begin{equation}\label{e:2.10}
\begin{CD}
H_\ast\bigl( U_+^{(1)}, \tilde U_+^{(1)};\K\bigr)
@>\hbox{isomorphism}>> H_\ast\bigl(U_+^{(1)}\setminus
\tilde U_\delta^{(1)}, \tilde U_+^{(1)}\setminus \tilde U_\delta^{(1)};\K\bigr) \\
@V J_\ast VV @VV J_\ast V \\
H_\ast\bigl( U_+^{(2)}, \tilde U_+^{(2)};\K\bigr)
@>\hbox{isomorphism}
>>H_\ast\bigl(U_+^{(2)}\setminus \tilde U_\delta^{(2)}, \tilde
U_+^{(2)}\setminus \tilde U_\delta^{(2)};\K\bigr),
\end{CD}
\end{equation}
where two isomorphisms are given by the excision property. Moreover,
the reversed flows
$$
\eta^{(j)}(-t):\bigl(U_+^{(j)}\setminus \tilde U_\delta^{(j)},
\tilde U_+^{(j)}\setminus \tilde U_\delta^{(j)}\bigr)\to\bigl( W_j,
W_j^-\bigr),\;j=1,2,
$$
are also strong deformation retracts. As in (\ref{e:2.9}) we get the
following commutative diagram:
\begin{equation}\label{e:2.11}
\begin{CD}
H_\ast\bigl(U_+^{(1)}\setminus \tilde U_\delta^{(1)}, \tilde
U_+^{(1)}\setminus \tilde U_\delta^{(1)};\K\bigr)
@>\hbox{isomorphism}>> H_\ast\bigl(W_1, W_1^-;\K\bigr) \\
@V J_\ast VV @VV J_\ast V \\
H_\ast\bigl(U_+^{(2)}\setminus \tilde U_\delta^{(2)}, \tilde
U_+^{(2)}\setminus \tilde
U_\delta^{(2)};\K\bigr)@>\hbox{isomorphism}
>>H_\ast\bigl(W_2, W_2^-;\K\bigr).
\end{CD}
\end{equation}
Finally,  by the Deformation Theorem 3.2 in \cite{Ch} (with the
flows of $-\nabla {\cal L}_{j}/\|\nabla {\cal L}_{j}\|^2_{H_j}$,
$j=1,2$) we have also the commutative diagram: {\scriptsize
$$
\begin{CD}
H_\ast\bigl(({\cal L}_{1})_{c_1+\varepsilon}\cap U_+^{(1)}, ({\cal
L}_{1})_{c_1-\varepsilon}\cap U_+^{(1)};\K\bigr) @>\hbox{isom}>>
H_\ast\bigl(({\cal L}_{1})_{c_1}\cap U_+^{(1)}, \bigl(({\cal
L}_{1})_{c_1}\setminus\{0\}\bigr)\cap U_+^{(1)};\K\bigr)
 \\
@V J_\ast VV @VV J_\ast V \\
H_\ast\bigl(({\cal L}_{2})_{c_2+k\varepsilon}\cap U_+^{(2)}, ({\cal
L}_{2})_{c_2-k\varepsilon}\cap U_+^{(2)};\K\bigr)@>\hbox{isom}
>>H_\ast\bigl(({\cal L}_{2})_{c_2}\cap
U_+^{(2)}, \bigl(({\cal L}_{2})_{c_2}\setminus\{0\}\bigr)\cap
U_+^{(2)};\K\bigr).
\end{CD}
$$}
From this and the commutative diagrams (\ref{e:2.9})-(\ref{e:2.11})
it follows that Claim~\ref{cl:2.2} is equivalent to

\begin{claim}\label{cl:2.3}
The map $J$ induces isomorphisms {\scriptsize
$$
J_\ast: H_\ast\bigl(({\cal L}_{1})_{c_1}\cap U_+^{(1)}, \bigl(({\cal
L}_{1})_{c_1}\setminus\{0\}\bigr)\cap U_+^{(1)};\K\bigr)\to
H_\ast\bigl(({\cal L}_{2})_{c_2}\cap U_+^{(2)}, \bigl(({\cal
L}_{2})_{c_2}\setminus\{0\}\bigr)\cap U_+^{(2)};\K\bigr).
$$}
\end{claim}

As in (\ref{e:2.10}) we may have
   the following commutative diagram:{\tiny
$$
\begin{CD}
H_\ast\bigl(({\cal L}_{1})_{c_1}\cap U_+^{(1)}, \bigl(({\cal
L}_{1})_{c_1}\setminus\{0\}\bigr)\cap U_+^{(1)};\K\bigr)
@>\hbox{iso}>> H_\ast\bigl(({\cal L}_{1})_{c_1}\cap
(U_+^{(1)}\setminus
\tilde U_\delta^{(1)}), (({\cal L}_{1})_{c_1}\setminus\{0\})\cap (U_+^{(1)}\setminus \tilde U_\delta^{(1)});\K\bigr) \\
@V J_\ast VV @VV J_\ast V \\
H_\ast\bigl(({\cal L}_{2})_{c_2}\cap U_+^{(2)}, \bigl(({\cal
L}_{2})_{c_2}\setminus\{0\}\bigr)\cap U_+^{(2)};\K\bigr)
@>\hbox{iso}
>>H_\ast\bigl(({\cal L}_{2})_{c_2}\cap (U_+^{(2)}\setminus \tilde U_\delta^{(2)}),
(({\cal L}_{2})_{c_2}\setminus\{0\})\cap U_+^{(2)}\setminus \tilde
U_\delta^{(2)};\K\bigr),
\end{CD}
$$}

 So Claim~\ref{cl:2.3} is equivalent to

\begin{claim}\label{cl:2.4}
The map $J$ induces isomorphisms  {\tiny
$$
J_\ast:H_\ast\bigl(({\cal L}_{1})_{c_1}\cap (U_+^{(1)}\setminus
\tilde U_\delta^{(1)}), (({\cal L}_{1})_{c_1}\setminus\{0\})\cap
(U_+^{(1)}\setminus \tilde U_\delta^{(1)});\K\bigr) \to
H_\ast\bigl(({\cal L}_{2})_{c_2}\cap (U_+^{(2)}\setminus \tilde
U_\delta^{(2)}), (({\cal L}_{2})_{c_2}\setminus\{0\})\cap
U_+^{(2)}\setminus \tilde U_\delta^{(2)};\K\bigr).
$$}
\end{claim}

Note that $U_+^{(j)}\setminus \tilde U_\delta^{(j)}$ are
neighborhoods $\theta_j\in H_j$, $j=1,2$. By the construction of
Gromoll-Meyer pairs (cf. \cite[page 49]{Ch}), for $\delta>0$
sufficiently small we can require that no other critical points of
${\cal L}_j$ is contained in them. Hence the excision property of
singular homology  implies that Claim~\ref{cl:2.4} is equivalent to

\begin{claim}\label{cl:2.5}
There exist small open neighborhoods $V^{(j)}$ of $\theta_j\in H_j$
 with $J(V^{(1)})\subset V^{(2)}$, such that $J$ induces isomorphisms {\footnotesize
$$
J_\ast: H_\ast\bigl(({\cal L}_{1})_{c_1}\cap V^{(1)}, (({\cal
L}_{1})_{c_1}\setminus\{0\})\cap V^{(1)};\K\bigr)\to
H_\ast\bigl({\cal L}_{2})_{c_2}\cap V^{(2)}, (({\cal
L}_{c_2})_{c_2}\setminus\{0\})\cap V^{(2)};\K\bigr).
$$}
\end{claim}

\noindent{\bf Step 2}. Consider the orthogonal decompositions
$$
H_j=H^0_{j}\oplus H^-_{j}\oplus H^+_{j}=H^0_{j}\oplus H^\bot_{j},
$$
where $H^0_j$, $H^-_j$ and $H^+_j$ are the null, negative, and
positive definite spaces of $B_{j}(\theta_j)$, $j=1,2$,
respectively. By Corollary~\ref{cor:1.5}
 there exist $\epsilon>0$,  the $C^1$ maps
 $$
 h_1:B_{H^0_1}(\theta,\epsilon)\to
X_1^\pm\quad\hbox{and}\quad
h_2:B_{H^0_2}(\theta,\sqrt{k}\epsilon)\to X^\pm_2
$$
satisfying $h_j(\theta^0_j)=\theta_j^\pm$, $j=1,2$, and the
origin-preserving homeomorphisms
\begin{eqnarray*}
&&\Phi_1: B_{H^0_1}(\theta,\epsilon)\oplus
B_{H^+_1}(\theta,\epsilon) \oplus
B_{H^-_1}(\theta,\epsilon)\to W_1,\\
&&\Phi_2: B_{H^0_2}(\theta,\sqrt{k}\epsilon)\oplus
B_{H^+_2}(\theta,\sqrt{k}\epsilon) \oplus
B_{H^-_2}(\theta,\sqrt{k}\epsilon)\to W_2,
\end{eqnarray*}
such that $h_2(Jz)=J\circ h_1(z)$ and
\begin{eqnarray}\label{e:2.12}
&& \Phi_2(Jz+ Ju^+ +
Ju^-)=J\circ\Phi_1(z+ u^+ + u^-),\\
&&{\cal L}_1\circ\Phi_1(z+ u^++
u^-)=\|u^+\|^2_{H_1}-\|u^-\|^2_{H_1}+ {\cal L}_1(z+
h_1(z))\nonumber\\
&&\hspace{37mm}\equiv\beta_1(u^++u^-)+ \alpha_1(z)\label{e:2.13}
\end{eqnarray}
 for all $(z, u^+, u^-)\in B_{H^0_1}(\theta,\epsilon)\times
B_{H^+_1}(\theta,\epsilon)\times B_{H^-_1}(\theta,\epsilon)$, and
that
\begin{eqnarray}\label{e:2.14}
&&{\cal L}_2\circ\Phi_2(z+ u^++
u^-)=\|u^+\|^2_{H_2}-\|u^-\|^2_{H_2}+ {\cal L}_2(z+
h_2(z))\nonumber\\
&&\hspace{37mm}\equiv\beta_2(u^++u^-)+ \alpha_2(z)
\end{eqnarray}
for all $(z, u^+, u^-)\in B_{H^0_2}(\theta,\sqrt{k}\epsilon)\times
B_{H^+_2}(\theta,\sqrt{k}\epsilon)\times
B_{H^-_2}(\theta,\sqrt{k}\epsilon)$.  Consequently,
\begin{equation}\label{e:2.15}
\alpha_2\circ J=k\alpha_1\quad\hbox{and}\quad \beta_2\circ
J=k\beta_1.
\end{equation}

Take open convex neighborhoods of the origin $\theta$ in $H^0_1,
H^-_1, H^+_1$, ${\cal U}^0_1, {\cal U}^-_1, {\cal U}^+_1$, and that
of $\theta$ in $H^+_{2}$, ${\cal U}^+_{2}$,
 such
that
\begin{eqnarray*}
&&{\cal U}_1:={\cal U}^0_1\oplus{\cal U}^-_1\oplus{\cal
U}^+_1\subset B_{H^0_1}(\theta,\epsilon)\oplus
B_{H^+_1}(\theta,\epsilon) \oplus
B_{H^-_1}(\theta,\epsilon),\\
&&\Phi_1({\cal U}_1)\subset B_{H^0_1}(\theta,\epsilon)\oplus
B_{H^+_1}(\theta,\epsilon) \oplus
B_{H^-_1}(\theta,\epsilon),\\
 &&J({\cal U}_1)\subset
B_{H^0_2}(\theta,\sqrt{k}\epsilon)\oplus
B_{H^+_2}(\theta,\sqrt{k}\epsilon) \oplus
B_{H^-_2}(\theta,\sqrt{k}\epsilon)
\end{eqnarray*}
and that $J({\cal U}^+_1)\subset{\cal U}^+_{2}$ and
\begin{eqnarray*}\label{e:2.19}
&&{\cal U}_{2}:= J({\cal U}^0_1)\oplus J({\cal U}^-_1)\oplus{\cal
U}_{2}^+\subset B_{H^0_2}(\theta,\sqrt{k}\epsilon)\oplus
B_{H^+_2}(\theta,\sqrt{k}\epsilon) \oplus
B_{H^-_2}(\theta,\sqrt{k}\epsilon),\\
&&\Phi_2({\cal U}_{2})\subset
B_{H^0_2}(\theta,\sqrt{k}\epsilon)\oplus
B_{H^+_2}(\theta,\sqrt{k}\epsilon) \oplus
B_{H^-_2}(\theta,\sqrt{k}\epsilon).
\end{eqnarray*}
By (\ref{e:2.12}) we have the commutative diagrams
$$
\begin{CD}
{\cal U}_1 @> \Phi_1 >> \Phi_1({\cal U}_1)
 \\
@V J VV @VV J V \\
{\cal U}_{2}@> \Phi_2
>>\Phi_2({\cal U}_{2}).
\end{CD}
$$
and thus {\scriptsize
$$
\begin{CD}
H_\ast\bigl(({\cal L}_{1}\circ\Phi_1)_{c_1}\cap {\cal U}_1, (({\cal
L}_{1}\circ\Phi_1)_{c_1}\setminus\{0\})\cap {\cal U}_1;\K\bigr)
@>(\Phi_1)_\ast
>> H_\ast\bigl(({\cal
L}_{1})_{c_1}\cap\Phi_1({\cal U}_1), (({\cal
L}_{1})_{c_1}\setminus\{0\})\cap\Phi_1({\cal U}_1);\K\bigr)
 \\
@V J_\ast VV @VV J_\ast V \\
H_\ast\bigl(({\cal L}_{2}\circ\Phi_2)_{c_2}\cap {\cal U}_{2},
(({\cal L}_{2}\circ\Phi_{2})_{c_2}\setminus\{0\})\cap {\cal
U}_{2};\K\bigr)@>(\Phi_2)_\ast
>>H_\ast\bigl(({\cal
L}_{2})_{c_2}\cap\Phi_{2}({\cal U}_{2}), (({\cal
L}_{2})_{c_2}\setminus\{0\})\cap\Phi_{2}({\cal U}_{2});\K\bigr).
\end{CD}
$$}
By this,  (\ref{e:2.13})-(\ref{e:2.14}), and  ${\cal
L}_{j}\circ\Phi_{j}=\beta_{j}+ \alpha_{j}$, $j=1,2$, and the fact
that $(\Phi_{j})_\ast$ are isomorphisms, $j=1,2$, taking
$V^{(1)}=\Phi_{1}({\cal U}_{1})$ and $V^{(2)}=\Phi_{2}({\cal
U}_{2})$,  Claim~\ref{cl:2.5} is equivalent to

\begin{claim}\label{cl:2.6}
$J$ induces  isomorphisms {\scriptsize
$$
H_\ast\bigl((\beta_{1}+ \alpha_{1})_{c_1}\cap {\cal U}_1,
((\beta_{1}+ \alpha_{1})_{c_1}\setminus\{0\})\cap {\cal
U}_1;\K\bigr)\to H_\ast\bigl((\beta_{2}+ \alpha_{2})_{c_2}\cap {\cal
U}_{2}, ((\beta_{2}+ \alpha_{2})_{c_2}\setminus\{0\})\cap {\cal
U}_{2};\K\bigr) .
$$}
\end{claim}
Since the deformation  retracts
\begin{eqnarray*}
&&H^0_1\oplus H^-_1\oplus H^+_1\times [0, 1]\to H^0_1\oplus
H^-_1\oplus H^+_1,\\
&&\qquad (x^0+ x^-+ x^+, t)\mapsto
x^0+ x^-+ tx^+,\\
&&H^0_{2}\oplus H^-_{2}\oplus H^+_{2}\times [0, 1]\to H^0_{2}\oplus
H^-_{2}\oplus H^+_{2},\\
&&\qquad (x^0+ x^-+ x^+, t)\mapsto x^0+ x^-+ tx^+
\end{eqnarray*}
commute with $J$, Claim~\ref{cl:2.6} is equivalent to

\begin{claim}\label{cl:2.7}
$J$ induces  isomorphisms from
\begin{eqnarray*}
H_\ast\bigl((\beta_{1}+ \alpha_{1})_{c_1}\cap ({\cal
U}^0_1\oplus{\cal U}_1^-\oplus\{0\}), ((\beta_{1}+
\alpha_{1})_{c_1}\setminus\{0\})\cap ({\cal U}^0_1\oplus{\cal
U}_1^-\oplus\{0\});\K\bigr)
\end{eqnarray*}
to {\scriptsize
\begin{eqnarray*}
H_\ast\bigl((\beta_{2}+ \alpha_{2})_{c_2}\cap (J({\cal
U}^0_{1})\oplus J({\cal U}_{1}^-)\oplus\{0\}),  ((\beta_{2}+
\alpha_{2})_{c_2}\setminus\{0\})\cap (J({\cal U}^0_{1})\oplus
J({\cal U}_{1}^-)\oplus\{0\});\K\bigr) .
\end{eqnarray*}}
\end{claim}

But $J:{\cal U}^0_1\oplus{\cal U}_1^-\oplus\{0\}\to J({\cal
U}^0_{\tau})\oplus J({\cal U}_{\tau}^-)\oplus\{0\}$ is a linear
diffeomorphism and
$$
(\beta_{2}+ \alpha_{2})( Jx^0+ Jx^-)=k(\beta_{1}+ \alpha_{1})(x^0+
x^-)\quad\forall x^0+ x^-\in{\cal U}^0_1\oplus{\cal U}^-_1
$$
because of (\ref{e:2.15}). Claim~\ref{cl:2.7} follows immediately.
Hence the homomorphisms in (\ref{e:2.3}) are isomorphisms.
 Theorem~\ref{th:2.1} is
proved. $\Box$.\vspace{2mm}

\begin{corollary}\label{cor:2.8}
Under the assumptions of Theorem~\ref{th:2.1} one has:\\
\noindent{\rm (i)} For any neighborhoods $\widetilde V_i$ of
$\theta_i\in H_1$ with $J(\widetilde V_1)\subset\widetilde V_2$ the
map $J$ induces isomorphisms
$$
J_\ast: H_\ast\bigl(({\cal L}_{1})_{c_1}\cap \widetilde V_1, (({\cal
L}_{1})_{c_1}\setminus\{0\})\cap \widetilde V_1;\K\bigr)\to
H_\ast\bigl({\cal L}_{2})_{c_2}\cap \widetilde V_2, (({\cal
L}_{c_2})_{c_2}\setminus\{0\})\cap \widetilde V_2;\K\bigr).
$$
\noindent{\rm (ii)} For any  Gromoll-Meyer pairs of ${\cal L}_i$ at
$\theta_i\in H_1$ (with respect to the negative gradient flows),
$(\widehat W_1,  \widehat W^-_1)$ with $(J(\widehat W_1), J(\widehat
W^-_1))\subset (\widehat W_{2}, \widehat W^-_{2})$, the map $J$
induces isomorphisms
$$
J_\ast: H_\ast(\widehat W_1, \widehat W^-_1; \K)\to H_\ast(\widehat
W_2, \widehat W^-_2; \K)
$$
\end{corollary}

\noindent{\bf Proof}. (i) For the neighborhoods $V^{(i)}$ in
Claim~\ref{cl:2.5} let us take open neighborhoods $\widehat V_i$ of
$\theta_i\in H_1$ with $J(\widehat V_1)\subset\widehat V_2$, such
that $Cl(\widehat V_i)\subset {\rm Int}(\widetilde V_i)\cap {\rm
Int}(V^{(i)})$, $i=1,2$. Then we have the commutative diagrams:
{\scriptsize
$$
\begin{CD}
H_\ast\bigl(({\cal L}_{1})_{c_1}\cap V^{(1)}, (({\cal
L}_{1})_{c_1}\setminus\{0\})\cap V^{(1)};\K\bigr) @> {\rm Isom}
>> H_\ast\bigl(({\cal
L}_{1})_{c_1}\cap \widehat V_1, (({\cal
L}_{1})_{c_1}\setminus\{0\})\cap \widehat V_1;\K\bigr)
 \\
@V J_\ast VV @VV J_\ast V \\
H_\ast\bigl(({\cal L}_{2})_{c_2}\cap V^{(2)}, (({\cal
L}_{2})_{c_2}\setminus\{0\})\cap V^{(2)};\K\bigr)@>{\rm Isom}
>>H_\ast\bigl(({\cal
L}_{2})_{c_2}\cap \widehat V_2, (({\cal
L}_{2})_{c_2}\setminus\{0\})\cap\widehat V_2;\K\bigr)
\end{CD}
$$}
and {\scriptsize
$$
\begin{CD}
H_\ast\bigl(({\cal L}_{1})_{c_1}\cap \widetilde V_1, (({\cal
L}_{1})_{c_1}\setminus\{0\})\cap \widetilde V_1;\K\bigr) @> {\rm
Isom}
>> H_\ast\bigl(({\cal
L}_{1})_{c_1}\cap \widehat V_1, (({\cal
L}_{1})_{c_1}\setminus\{0\})\cap \widehat V_1;\K\bigr)
 \\
@V J_\ast VV @VV J_\ast V \\
H_\ast\bigl(({\cal L}_{2})_{c_2}\cap \widetilde V_2, (({\cal
L}_{2})_{c_2}\setminus\{0\})\cap \widetilde V_2;\K\bigr)@>{\rm Isom}
>>H_\ast\bigl(({\cal
L}_{2})_{c_2}\cap \widehat V_2, (({\cal
L}_{2})_{c_2}\setminus\{0\})\cap\widehat V_2;\K\bigr)
\end{CD}
$$}
Here four ``Isom'' come from the excision property.
Claim~\ref{cl:2.5} gives the desired conclusion.\

\noindent{(ii)} By the proof of Theorem~\ref{th:2.1} the conclusion
required is equivalent to the corresponding result of
Claim~\ref{cl:2.3}, that is, \textsf{the map $J$ induces
isomorphisms} {\scriptsize
$$
J_\ast: H_\ast\bigl(({\cal L}_{1})_{c_1}\cap \widehat U_+^{(1)},
\bigl(({\cal L}_{1})_{c_1}\setminus\{0\}\bigr)\cap \widehat
U_+^{(1)};\K\bigr)\to H_\ast\bigl(({\cal L}_{2})_{c_2}\cap \widehat
U_+^{(2)}, \bigl(({\cal L}_{2})_{c_2}\setminus\{0\}\bigr)\cap
\widehat U_+^{(2)};\K\bigr),
$$}
where
\begin{eqnarray*}
\widehat U_+^{(j)}= \cup_{0<t<\infty}\eta^{(j)}(t, \widehat
W_j),\quad j=1,2.
\end{eqnarray*}
Since $\widehat U_+^{(j)}$ are neighborhoods of $\theta_j\in H_j$
and $J(\widehat U_+^{(1)})\subset \widehat U_+^{(2)}$ by
(\ref{e:2.7}), the desired conclusion follows from (i).
$\Box$\vspace{2mm}

\section{Variational setup}\label{sec:3}
\setcounter{equation}{0}

For integers $m\ge 3$ and $k\in\N$,  a compact $C^m$-smooth manifold
$M$ without boundary and $C^{m-1}$-smooth $L:\R\times TM\to\R$
satisfying the assumptions (L1)-(L3) in \cite[Section 1]{Lu}, on the
$C^{m-1}$-smooth Hilbert manifold $H_{k\tau}=W^{1,2}(S_{k\tau}, M)$,
where $S_{k\tau}:=\R/k\tau\Z=\{[s]_{k\tau}\,|\,[s]_{k\tau}=s+
k\tau\Z,\, s\in\R\}$,
\begin{equation}\label{e:3.1}
{\cal L}_{k\tau}(\gamma)=\int^{k\tau}_0L(t, \gamma(t),\dot
\gamma(t))dt \quad\forall \gamma\in H_{k\tau},
\end{equation}
defines a functional ${\cal L}_{k\tau}$ is $C^{2-0}$-smooth, bounded
below, satisfies the Palais-Smale condition (cf.
\cite[Prop.2.2]{AbSc1}). By \cite[Th.3.7.2]{Fa}, all critical points
of ${\cal L}_{k\tau}$ are all of class $C^{m-1}$ and therefore
correspond to all $k\tau$-periodic solutions of  the Lagrangian
system on $M$:
\begin{equation}\label{e:3.2}
\frac{d}{dt}\Big(\frac{\partial
L}{\partial\dot{q_i}}\Big)-\frac{\partial L}{\partial q_i}=0
\end{equation}
in any local coordinates $(q_1,\cdots, q_n)$. However, by the
condition (L2) and \cite[pp.175]{Du}, the functional ${\cal
L}_{k\tau}$ is $C^2$-smooth on the $C^2$-Banach manifold
$$
X_{k\tau}:=C^1(S_{k\tau}, M)
$$
with the usual topology of uniform convergence of the curves and
their derivatives. So ${\cal L}_{k\tau}$ has the same critical point
set on $H_{k\tau}$ and $X_{k\tau}$. Denote by  ${\cal L}_{k\tau}^X$
the restriction of ${\cal L}_{k\tau}|_{X_{k\tau}}$ to $X_{k\tau}$.

For a critical point $\gamma_0$ of ${\cal L}_{k\tau}$, which
actually sits in $C^2(S_{k\tau}, M)\subset X_{k\tau}$ due to our
assumptions, by the proof of Theorem 3.1 in \cite{Lu}, near
$\gamma_0$  we can pullback $L$ to $\tilde L:\R\times
B^n_\rho(0)\times\R^n$ by \cite[(3.15)]{Lu}. Denote by
$$
\tilde V_{k\tau}:=W^{1,2}(S_{k\tau}, B^n_\rho(0)),\qquad \tilde
X_{k\tau}:=C^1(S_{k\tau}, \R^n),\qquad \tilde
H_{k\tau}=W^{1,2}(S_{k\tau}, \R^n).
$$
Let $\tilde\gamma\in \tilde V_{k\tau}\cap \tilde X_{k\tau}$ (or
$\tilde V_{k\tau}$) be the pullback of $\gamma\in  X_{k\tau}$ (or
$H_{k\tau}$) near $\gamma_0$ by $\phi_{k\tau}$ as in
\cite[(3.8)]{Lu}. Then $\tilde\gamma_0=0$. Define
\begin{eqnarray*}
&&\tilde{\cal L}_{k\tau}(\tilde\alpha)=\int^{k\tau}_0 \tilde
L\left(t, \tilde\alpha(t),
\dot{\tilde\alpha}(t)\right)dt\quad\forall\tilde\alpha\in
\tilde V_{k\tau},\\
&&{\tilde{\cal L}}^X_{k\tau}(\tilde\alpha)=\int^{k\tau}_0 \tilde
L\left(t, \tilde\alpha(t),
\dot{\tilde\alpha}(t)\right)dt\quad\forall\tilde\alpha\in \tilde
V_{k\tau}\cap \tilde X_{k\tau}.
\end{eqnarray*}
 Then $\tilde{\cal L}_{k\tau}$ is $C^{2-0}$ in
$\tilde V_{k\tau}\subset \tilde H_{k\tau}$, and ${\tilde{\cal
L}}^X_{k\tau}$ is $C^2$ in $\tilde V_{k\tau}\cap \tilde
X_{k\tau}\subset\tilde X_{k\tau}$. Moreover, the zero  is the
critical point of both.

We shall prove that the functional $\tilde{\cal L}_{k\tau}$, and
spaces $\tilde X_{k\tau}$, $\tilde H_{k\tau}$ and $\tilde V_{k\tau}$
satisfy the conditions of Theorem~\ref{th:1.1}. Recall that
\begin{eqnarray*}
d\tilde{\cal L}_{k\tau}(\tilde\gamma) (\tilde \xi)
 \!\!\!&=&\!\!\!\int_0^{k\tau}  \left( D_{\tilde q} \tilde
L\left(t,\tilde\gamma(t),\dot{\tilde\gamma}(t)\right)\cdot\tilde\xi(t)
+ D_{\tilde v} \tilde
L\left(t,\tilde\gamma(t),\dot{\tilde\gamma}(t)\right)\cdot\dot{\tilde\xi}(t)
\right)  \, dt
\end{eqnarray*}
for any $\tilde\gamma\in \tilde V_{k\tau},\;\tilde\xi\in \tilde
H_{k\tau}$ and $k\in\N$, and that
\begin{eqnarray}\label{e:3.3}
 d^2 {\tilde{\cal L}}^X_{k\tau}
  (\tilde\gamma)(\tilde\xi,\tilde\eta)
   = \int_0^{k\tau} \Bigl(\!\! \!\!\!&&\!\!\!\!\!D_{\tilde v\tilde v}
    \tilde L\left(t,\tilde\gamma(t),\dot{\tilde\gamma}(t)\right)
\bigl(\dot{\tilde\xi}(t), \dot{\tilde \eta}(t)\bigr) \nonumber\\
&&+ D_{\tilde q\tilde v} \tilde
  L\left(t,\tilde\gamma(t), \dot{\tilde\gamma}(t)\right)
\bigl(\tilde\xi(t), \dot{\tilde\eta}(t)\bigr)\nonumber \\
&& + D_{\tilde v\tilde q} \tilde
  L\left(t,\tilde\gamma(t),\dot{\tilde\gamma}(t)\right)
\bigl(\dot{\tilde\xi}(t), \tilde\eta(t)\bigr) \nonumber\\
&&+  D_{\tilde q\tilde q} \tilde L\left(t,\tilde\gamma(t),
\dot{\tilde\gamma}(t)\right) \bigl(\tilde\xi(t),\tilde
\eta(t)\bigr)\Bigr) \, dt
\end{eqnarray}
for any $\tilde\gamma\in \tilde V_{k\tau}\cap \tilde X_{k\tau},\;
\tilde\xi,\tilde\eta\in \tilde X_{k\tau}$ and $k\in\N$.
 Let  $\nabla\tilde{\cal
L}_{k\tau}(\tilde\gamma)\in \tilde H_{k\tau}$ be the gradient of
$\tilde{\cal L}_{k\tau}$ at $\tilde\gamma\in \tilde V_{k\tau}$. If
$\tilde\gamma\in \tilde V_{k\tau}\cap \tilde X_{k\tau}$ and
$\tilde\xi\in \tilde X_{k\tau}$, then
\begin{equation}\label{e:3.4}
d{\tilde{\cal L}}^X_{k\tau}(\tilde\gamma) (\tilde \xi)=d\tilde{\cal
L}_{k\tau}(\tilde\gamma) (\tilde \xi)=\bigl(\nabla\tilde{\cal
L}_{k\tau}(\tilde\gamma), \tilde \xi\bigr)_{W^{1,2}}.
\end{equation}
We need to compute $\nabla\tilde{\cal
L}_{k\tau}(\tilde\gamma)\in \tilde H_{k\tau}$. Note that the
function $s\mapsto G_{k\tau}(\tilde\gamma)(s)$ given by
$$
G_{k\tau}(\tilde\gamma)(s):=\int^s_0\Bigl[D_{\tilde v} \tilde
L\left(t,\tilde\gamma(t),\dot{\tilde\gamma}(t)\right)-\frac{1}{k\tau}\int^{k\tau}_0D_{\tilde
v} \tilde
L\left(t,\tilde\gamma(t),\dot{\tilde\gamma}(t)\right)dt\Bigr]dt
$$
 is a $k\tau$-periodic primitive function of the
function
$$
s\mapsto D_{\tilde v} \tilde
L\left(s,\tilde\gamma(s),\dot{\tilde\gamma}(s)\right)-\frac{1}{k\tau}\int^{k\tau}_0D_{\tilde
v} \tilde L\left(s,\tilde\gamma(s),\dot{\tilde\gamma}(s)\right)ds
$$
and that
\begin{eqnarray*}
&&\quad\int_0^{k\tau}   D_{\tilde v} \tilde
L\left(t,\tilde\gamma(t),\dot{\tilde\gamma}(t)\right)\cdot\dot{\tilde\xi}(t)
  \, dt\\
  &&=\int_0^{k\tau}\Bigl[D_{\tilde v} \tilde
L\left(t,\tilde\gamma(t),\dot{\tilde\gamma}(t)\right)-\frac{1}{k\tau}\int^{k\tau}_0D_{\tilde
v} \tilde
L\left(t,\tilde\gamma(t),\dot{\tilde\gamma}(t)\right)dt\Bigr]\cdot\dot{\tilde\xi}(t)
  \, dt.
\end{eqnarray*}
Hence
\begin{eqnarray*}
&&\quad\int_0^{k\tau}  \left( D_{\tilde q} \tilde
L\left(t,\tilde\gamma(t),\dot{\tilde\gamma}(t)\right)\cdot\tilde\xi(t)
+ D_{\tilde v} \tilde
L\left(t,\tilde\gamma(t),\dot{\tilde\gamma}(t)\right)\cdot\dot{\tilde\xi}(t)
\right)  \, dt\\
&&=\int_0^{k\tau}  \left( D_{\tilde q} \tilde
L\left(t,\tilde\gamma(t),\dot{\tilde\gamma}(t)\right)-G_{k\tau}(\tilde\gamma)(t)\right)\cdot\tilde\xi(t)\,dt
+(G_{k\tau}(\tilde\gamma),\xi)_{W^{1,2}}.
\end{eqnarray*}

\begin{lemma}\label{lem:3.1}
If $f\in L^1(S_T, \R^n)$ is bounded, then the equation
$$
x''(t)- x(t)=f(t)
$$
has an unique $T$-periodic solution
$$
x(t)=\frac{1}{2}\int^\infty_t e^{t-s}f(s)\,ds +
\frac{1}{2}\int_{-\infty}^t e^{s-t}f(s)\,ds.
$$
\end{lemma}

Since
$$
\tilde\xi\mapsto \int_0^{k\tau}  \left( D_{\tilde q} \tilde
L\left(t,\tilde\gamma(t),\dot{\tilde\gamma}(t)\right)-G_{k\tau}(\tilde\gamma)(t)\right)\cdot\tilde\xi(t)\,dt
$$
is a bounded linear functional on $\tilde H_{k\tau}$, the Riesz
theorem yields an unique $F(\tilde\gamma)\in \tilde H_{k\tau}$ such
that
\begin{equation}\label{e:3.5}
\int_0^{k\tau}  \left( D_{\tilde q} \tilde
L\left(t,\tilde\gamma(t),\dot{\tilde\gamma}(t)\right)-
G_{k\tau}(\tilde\gamma)(t)\right)\cdot\tilde\xi(t)\,dt=\bigl(F(\tilde\gamma),\tilde\xi\bigr)_{W^{1,2}}
\end{equation}
for any $\tilde\xi\in \tilde H_{k\tau}$. It follows that
\begin{equation}\label{e:3.6}
\nabla\tilde{\cal L}_{k\tau}(\tilde\gamma)= G_{k\tau}(\tilde\gamma)+
F(\tilde\gamma).
\end{equation}
By Lemma~\ref{lem:3.1} and a direct computation we get
\begin{eqnarray*}
F(\tilde\gamma)(t)&=&\frac{e^t}{2}\int^\infty_te^{-s}\left(
G_{k\tau}(\tilde\gamma)(s)- D_{\tilde q} \tilde
L\left(s,\tilde\gamma(s),\dot{\tilde\gamma}(s)\right)\right)\, ds\\
&+&\frac{e^{-t}}{2}\int^t_{-\infty}e^{s}\left(G_{k\tau}(\tilde\gamma)(s)-D_{\tilde
q} \tilde L\left(s,\tilde\gamma(s),\dot{\tilde\gamma}(s)\right)
\right)\, ds
\end{eqnarray*}
 for any $t\in\R$. This and (\ref{e:3.6}) lead to
\begin{eqnarray*}
\nabla\tilde{\cal
L}_{k\tau}(\tilde\gamma)(t)&=&\frac{e^t}{2}\int^\infty_te^{-s}\left(
G_{k\tau}(\tilde\gamma)(s) -D_{\tilde q} \tilde
L\left(s,\tilde\gamma(s),\dot{\tilde\gamma}(s)\right)\right)\, ds\\
&+&\frac{e^{-t}}{2}\int^t_{-\infty}e^{s}\left(G_{k\tau}(\tilde\gamma)(s)-D_{\tilde
q} \tilde L\left(s,\tilde\gamma(s),\dot{\tilde\gamma}(s)\right)
\right)\, ds + G_{k\tau}(\tilde\gamma)(t),\\
\frac{d}{dt}\nabla\tilde{\cal
L}_{k\tau}(\tilde\gamma)(t)&=&\frac{e^t}{2}\int^\infty_te^{-s}\left(
G_{k\tau}(\tilde\gamma)(s) D_{\tilde q}- \tilde
L\left(s,\tilde\gamma(s),\dot{\tilde\gamma}(s)\right)\right)\, ds\\
&-&\frac{e^{-t}}{2}\int^t_{-\infty}e^{s}\left(G_{k\tau}(\tilde\gamma)(s)-D_{\tilde
q} \tilde L\left(s,\tilde\gamma(s),\dot{\tilde\gamma}(s)\right)
\right)\, ds \\
&+& D_{\tilde v} \tilde
L\left(t,\tilde\gamma(t),\dot{\tilde\gamma}(t)\right)-\frac{1}{k\tau}\int^{k\tau}_0D_{\tilde
v} \tilde L\left(s,\tilde\gamma(s),\dot{\tilde\gamma}(s)\right)\,ds.
\end{eqnarray*}
From these it easily follows that
\begin{eqnarray*}
&& \tilde\gamma\in \tilde V_{k\tau}\cap \tilde
X_{k\tau}\Longrightarrow\nabla\tilde{\cal
L}_{k\tau}(\tilde\gamma)\in
\tilde X_{k\tau}\;\hbox{and}\\
&&\tilde V_{k\tau}\cap \tilde X_{k\tau}\ni\tilde\gamma\mapsto
\nabla\tilde{\cal L}_{k\tau}(\tilde\gamma)\in \tilde
X_{k\tau}\quad\hbox{is continuous}.
\end{eqnarray*}

\begin{lemma}\label{lem:}\label{lem:3.2}
With the topology on $\tilde V_{k\tau}\cap \tilde X_{k\tau}$ induced
from $\tilde X_{k\tau}$ the map
$$
A_{k\tau}: \tilde V_{k\tau}\cap \tilde X_{k\tau}\to \tilde X_{k\tau}
$$
defined by $A_{k\tau}(\tilde\gamma)=\nabla\tilde{\cal
L}_{k\tau}(\tilde\gamma)$
 is continuously differentiable.
\end{lemma}

\noindent{\bf Proof}. For $\tilde\gamma\in \tilde V_{k\tau}\cap
\tilde X_{k\tau}$ and $\tilde\xi\in \tilde X_{k\tau}$, a direct
computation gives
\begin{eqnarray*}
&&\hspace{-5mm}G'_{k\tau}(\tilde\gamma)(\tilde\xi)(t)=\int^t_0\Bigl[D_{\tilde
v\tilde q} \tilde
L\left(s,\tilde\gamma(s),\dot{\tilde\gamma}(s)\right)\cdot{\tilde\xi}(s)+
D_{\tilde v\tilde v} \tilde
L\left(s,\tilde\gamma(s),\dot{\tilde\gamma}(s)\right)\cdot\dot{\tilde\xi}(s)
-\\
&&\hspace{15mm}\frac{1}{k\tau}\int^{k\tau}_0\left(D_{\tilde v\tilde
q} \tilde
L\left(s,\tilde\gamma(s),\dot{\tilde\gamma}(s)\right)\cdot{\tilde\xi}(s)+
D_{\tilde v\tilde v} \tilde
L\left(s,\tilde\gamma(s),\dot{\tilde\gamma}(s)\right)\cdot\dot{\tilde\xi}(s)\right)ds\Bigr]ds,\\
&&\hspace{-5mm} \frac{d}{dt}G'_{k\tau}(\tilde\gamma)(\tilde\xi)(t)=
D_{\tilde v\tilde q} \tilde
L\left(t,\tilde\gamma(t),\dot{\tilde\gamma}(t)\right)\cdot{\tilde\xi}(t)+
D_{\tilde v\tilde v} \tilde
L\left(t,\tilde\gamma(t),\dot{\tilde\gamma}(t)\right)\cdot\dot{\tilde\xi}(t)
-\\
&&\hspace{15mm}\frac{1}{k\tau}\int^{k\tau}_0\left(D_{\tilde v\tilde
q} \tilde
L\left(s,\tilde\gamma(s),\dot{\tilde\gamma}(s)\right)\cdot{\tilde\xi}(s)+
D_{\tilde v\tilde v} \tilde
L\left(s,\tilde\gamma(s),\dot{\tilde\gamma}(s)\right)\cdot\dot{\tilde\xi}(s)\right)ds.
\end{eqnarray*}
It follows that for any $\varepsilon>0$ there exists
$\delta=\delta(\tilde\gamma)>0$ such that
$$
\|(G'_{k\tau}(\tilde\gamma+ \tilde
h)-G'_{k\tau}(\tilde\gamma))(\tilde\xi)\|_{C^1}\le\varepsilon\|\tilde\xi\|_{C^1}
$$
for any $\tilde h\in \tilde V_{k\tau}\cap \tilde X_{k\tau}$ with
$\|\tilde h\|_{C^1}<\delta$ and $\tilde \xi\in \tilde X_{k\tau}$.
Namely
$$
\tilde V_{k\tau}\cap \tilde X_{k\tau}\ni\tilde\gamma\mapsto
G_{k\tau}(\tilde\gamma)\in \tilde
X_{k\tau}\quad\hbox{is}\;C^1-\hbox{smooth}.
$$
Similarly, we have
\begin{eqnarray*}
A'_{k\tau}(\tilde\gamma)(\tilde\xi)(t)&=&\frac{e^t}{2}\int^\infty_te^{-s}\Bigl(
G'_{k\tau}(\tilde\gamma)(\tilde\xi)(s) -D_{\tilde q\tilde q} \tilde
L\bigl(s,\tilde\gamma(s),\dot{\tilde\gamma}(s)\bigr)\cdot\tilde{\xi}(s)\\
&&\hspace{15mm} -D_{\tilde q\tilde v} \tilde
L\bigl(s,\tilde\gamma(s),\dot{\tilde\gamma}(s)\bigr)\cdot\dot{\tilde\xi}(s)
\Bigr)\, ds\\
&+&\frac{e^{-t}}{2}\int^t_{-\infty}e^{s}\Bigl(
G'_{k\tau}(\tilde\gamma)(\tilde\xi)(s) -D_{\tilde q\tilde q} \tilde
L\bigl(s,\tilde\gamma(s),\dot{\tilde\gamma}(s)\bigr)\cdot\tilde{\xi}(s)\\
&&\hspace{15mm} -D_{\tilde q\tilde v} \tilde
L\bigl(s,\tilde\gamma(s),\dot{\tilde\gamma}(s)\bigr)\cdot\dot{\tilde\xi}(s)
\Bigr)\, ds + G'_{k\tau}(\tilde\gamma)(\tilde\xi)(t),\\
\frac{d}{dt}A'_{k\tau}(\tilde\gamma)(\tilde\xi)(t)&=&\frac{e^t}{2}\int^\infty_te^{-s}\Bigl(
G'_{k\tau}(\tilde\gamma)(\tilde\xi)(s) -D_{\tilde q\tilde q} \tilde
L\bigl(s,\tilde\gamma(s),\dot{\tilde\gamma}(s)\bigr)\cdot\tilde{\xi}(s)\\
&&\hspace{15mm} -D_{\tilde q\tilde v} \tilde
L\bigl(s,\tilde\gamma(s),\dot{\tilde\gamma}(s)\bigr)\cdot\dot{\tilde\xi}(s)
\Bigr)\, ds\\
&-&\frac{e^{-t}}{2}\int^t_{-\infty}e^{s}\Bigl(
G'_{k\tau}(\tilde\gamma)(\tilde\xi)(s) -D_{\tilde q\tilde q} \tilde
L\bigl(s,\tilde\gamma(s),\dot{\tilde\gamma}(s)\bigr)\cdot\tilde{\xi}(s)\\
&&\hspace{15mm} -D_{\tilde q\tilde v} \tilde
L\bigl(s,\tilde\gamma(s),\dot{\tilde\gamma}(s)\bigr)\cdot\dot{\tilde\xi}(s)
\Bigr)\, ds + \frac{d}{dt}G'_{k\tau}(\tilde\gamma)(\tilde\xi)(t).
\end{eqnarray*}
From them it easily follows that for any $\varepsilon>0$ there
exists $\delta=\delta(\tilde\gamma)>0$ such that
$$
\|(A'_{k\tau}(\tilde\gamma+ \tilde
h)-A'_{k\tau}(\tilde\gamma))(\tilde\xi)\|_{C^1}\le\varepsilon\|\tilde\xi\|_{C^1}
$$
for any $\tilde h\in \tilde V_{k\tau}\cap \tilde X_{k\tau}$ with
$\|\tilde h\|_{C^1}<\delta$ and  $\tilde \xi\in \tilde X_{k\tau}$,
and thus
$$
\|A'_{k\tau}(\tilde\gamma+ \tilde
h)-A'_{k\tau}(\tilde\gamma)\|_{C^1}\le\varepsilon\|\tilde\xi\|
$$
for any $\tilde h\in \tilde V_{k\tau}\cap \tilde X_{k\tau}$ with
$\|\tilde h\|_{C^1}<\delta$. This proves that $A_{k\tau}$ is $C^1$.
$\Box$\vspace{2mm}

In summary,  the functional $\tilde{\cal L}_{k\tau}$ satisfies the
condition (\ref{e:1.1}) in Section 1 for $\tilde X_{k\tau}$, $\tilde
V_{k\tau}$ and $A=\nabla\tilde{\cal L}_{k\tau}|_{\tilde
V_{k\tau}\cap \tilde X_{k\tau}}$.

\begin{remark}\label{rm:3.3}
{\rm For $\tilde\gamma\in \tilde V_\tau\cap \tilde X_\tau$ and
$\tilde\gamma^k\in \tilde V_{k\tau}\cap \tilde X_{k\tau}$, $k\in\N$,
we have
\begin{eqnarray*}
\nabla\tilde{\cal
L}_{k\tau}(\tilde\gamma^k)(t)&=&\frac{e^t}{2}\int^\infty_te^{-s}\left(
G_{\tau}(\tilde\gamma)(s)- D_{\tilde q}\tilde
L\left(s,\tilde\gamma(s),\dot{\tilde\gamma}(s)\right)\right)\, ds\\
&+&\frac{e^{-t}}{2}\int^t_{-\infty}e^{s}\left(G_{\tau}(\tilde\gamma)(s)-D_{\tilde
q} \tilde L\left(s,\tilde\gamma(s),\dot{\tilde\gamma}(s)\right)
\right)\, ds + G_{\tau}(\tilde\gamma)(t),
\end{eqnarray*}
which is also $\tau$-periodic. }
\end{remark}

By (\ref{e:3.3}) it is easily checked:
\begin{description}
\item[(i)] For any $\tilde\gamma\in \tilde V_{k\tau}\cap \tilde X_{k\tau}$ there
exists a constant $C(\tilde\gamma)$ such that
$$
|d^2 {\tilde{\cal L}}^X_{k\tau}
  (\tilde\gamma)(\tilde\xi,\tilde\eta)|\le
  C(\tilde\gamma)\|\tilde\xi\|_{W^{1,2}}\cdot\|\tilde\eta\|_{W^{1,2}}
\quad\forall \tilde\xi,\tilde\eta\in \tilde X_{k\tau};
$$
\item[(ii)] $\forall\varepsilon>0,\;\exists\;\delta>0$, such that
for any $\tilde\gamma_1, \tilde\gamma_2\in \tilde V_{k\tau}\cap
\tilde X_{k\tau}$ with
$\|\tilde\gamma_1-\tilde\gamma_2\|_{C^1}<\delta$,
$$
|d^2 {\tilde{\cal L}}^X_{k\tau}
  (\tilde\gamma_1)(\tilde\xi,\tilde\eta)- d^2{\tilde{\cal L}}^X_{k\tau}
  (\tilde\gamma_2)(\tilde\xi,\tilde\eta)|\le
  \varepsilon \|\tilde\xi\|_{W^{1,2}}\cdot\|\tilde\eta\|_{W^{1,2}}
\quad\forall \tilde\xi,\tilde\eta\in \Tilde X_{k\tau}.
$$
\end{description}
It follows (cf. \cite[Prop.2.1]{JM}) that there exists a  map
$$
B_{k\tau}: \tilde V_{k\tau}\cap \tilde X_{k\tau}\to L(\tilde
H_{k\tau}),
$$
which is uniformly continuous with respect to the induced topology
on $\tilde V_{k\tau}\cap \tilde X_{k\tau}$ from $\tilde X_{k\tau}$,
such that for any $\tilde\gamma\in \tilde V_{k\tau}\cap X_{k\tau}$
and $\tilde\xi, \tilde\eta\in
  \tilde X_{k\tau}$ one has
\begin{equation}\label{e:3.7}
d^2 {\tilde{\cal L}}^X_{k\tau}
  (\tilde\gamma)(\tilde\xi,\tilde\eta)=\bigl(B_{k\tau}(\tilde\gamma)\tilde\xi,
  \tilde\eta\bigr)_{W^{1,2}}.
\end{equation}
By (i) the right side of (\ref{e:3.3}) is also a bounded symmetric
bilinear form on $\tilde H_{k\tau}$, each $B_{k\tau}(\tilde\gamma)$
is a bounded linear self-adjoint operator on $\tilde H_{k\tau}$.
From (\ref{e:3.3}), (\ref{e:3.4}) and Lemma~\ref{lem:3.2} one easily
derive
\begin{equation}\label{e:3.8}
\bigl(dA_{k\tau}(\tilde\gamma)\tilde\xi,
\tilde\eta\bigr)_{W^{1,2}}=\bigl(B_{k\tau}(\tilde\gamma)\tilde\xi,
  \tilde\eta\bigr)_{W^{1,2}}\quad\forall\tilde\xi,
  \tilde\eta\in\tilde X_{k\tau}.
\end{equation}
That is, (\ref{e:1.2}) is satisfied.

Moreover, if $\tilde\gamma\in \tilde V_{k\tau}^E\cap E\tilde
X_{k\tau}$, where $\tilde V_{k\tau}^E=\tilde V_{k\tau}\cap E\tilde
H_{k\tau}$ and
$$
E\tilde H_{k\tau}:=\bigl\{\tilde\gamma\in\tilde H_{k\tau}\,|\,
\tilde\gamma(-t)=\tilde\gamma(t)\;\forall t\bigr\},\quad E\tilde
X_{k\tau}:=\bigl\{\tilde\gamma\in\tilde X_{k\tau}\,|\,
\tilde\gamma(-t)=\tilde\gamma(t)\;\forall t\bigr\},
$$
then it is not difficult to check that $A_{k\tau}\bigl(\tilde
V^E_{k\tau}\cap  E\tilde X_{k\tau} \bigr)\subset E\tilde X_{k\tau}$
and
$$
B_{k\tau}(\tilde\gamma)(E\tilde H_{k\tau})\subset E\tilde
H_{k\tau}\quad\forall\tilde\gamma\in \tilde V^E_{k\tau}\cap \tilde
EX_{k\tau}.
$$
Hence $A_{k\tau}$ and $B_{k\tau}$ restrict to a $C^1$ map
\begin{equation}\label{e:3.9}
A_{k\tau}^E: \tilde V^E_{k\tau}\cap \tilde EX_{k\tau} \to E\tilde
X_{k\tau}
\end{equation}
and a continuous map
\begin{equation}\label{e:3.10}
B_{k\tau}^E:
\tilde V^E_{k\tau}\cap \tilde EX_{k\tau} \to L_s(E\tilde X_{k\tau})
\end{equation}
respectively.  Let $\tilde{\cal L}^E_{k\tau}$ (resp. $\tilde{\cal
L}^{EX}_{k\tau}$) is the restriction of $\tilde{\cal L}_{k\tau}$
(resp. $\tilde{\cal L}^X_{k\tau}$) to $\tilde V^E_{k\tau}$ (resp.
$\tilde V^E_{k\tau}\cap E\tilde X_{k\tau}$). Then (\ref{e:3.7}) and
(\ref{e:3.8}) imply
\begin{eqnarray}
&&d\tilde{\cal L}^{EX}_{k\tau}(\tilde\gamma) (\tilde
\xi)=d\tilde{\cal L}^E_{k\tau}(\tilde\gamma) (\tilde
\xi)=\bigl(A^E_{k\tau}(\tilde\gamma), \tilde
\xi\bigr)_{W^{1,2}},\label{e:3.11}\\
&&  \bigl(dA^E_{k\tau}(\tilde\gamma)\tilde\xi,
\tilde\eta\bigr)_{W^{1,2}}=\bigl(B^E_{k\tau}(\tilde\gamma)\tilde\xi,
  \tilde\eta\bigr)_{W^{1,2}}\label{e:3.12}
\end{eqnarray}
for any $\tilde\gamma\in \tilde V^E_{k\tau}\cap E\tilde X_{k\tau}$
and $\tilde\xi, \tilde\eta\in
  E\tilde X_{k\tau}$.

For any $\tilde\gamma\in \tilde V_{k\tau}\cap X_{k\tau}$ set
\begin{eqnarray*}
&&\hat P_\gamma(t)=D_{\tilde v\tilde v}
    \tilde L\left(t,\tilde\gamma(t),\dot{\tilde\gamma}(t)\right),\\
&&\hat Q_\gamma(t)=D_{\tilde q\tilde v} \tilde
  L\left(t,\tilde\gamma(t), \dot{\tilde\gamma}(t)\right),\\
&&\hat R_\gamma(t)=D_{\tilde q\tilde q} \tilde
L\left(t,\tilde\gamma(t),
\dot{\tilde\gamma}(t)\right),\\
&&\hat L_\gamma(t, \tilde y, \tilde v)=\frac{1}{2}\hat
P_\gamma(t)\tilde v\cdot\tilde v + \hat Q_\gamma(t)\tilde
y\cdot\tilde v + \frac{1}{2}\hat
R_\gamma(t)\tilde y\cdot\tilde y
\end{eqnarray*}
and
$$
\hat f_{k\tau,\gamma}(\tilde y)=\int^{k\tau}_0\hat L_\gamma\left(t,
\tilde y(t),\dot {\tilde y}(t)\right)dt\quad\forall\tilde y\in
\tilde H_{k\tau}.
$$
Then $\hat f_{k\tau,\gamma}$ is
$C^2$-smooth on $\tilde H_{k\tau}$ and $\tilde X_{k\tau}$, and
$\tilde y=0\in \tilde H_{k\tau}$ is a critical point of $\hat
f_{k\tau,\gamma}$.
 It is also easily checked that
\begin{eqnarray}\label{e:3.13}
&& d^2 \hat{f}_{k\tau,\gamma}
  (0)(\tilde\xi,\tilde\eta)=\int^{k\tau}_0\left[\bigl(\hat P_\gamma\dot{\tilde\xi}+
\hat Q_\gamma\tilde\xi\bigr)\cdot\dot{\tilde\eta}+
Q_\gamma^T\dot{\tilde\xi}\cdot \tilde\eta+
 \hat R_\gamma\tilde\xi\cdot \tilde\eta\right]dt\nonumber\\
&&\hspace{4mm} =\int^{k\tau}_0  \Bigl(D_{\tilde v\tilde v}
    \tilde L\left(t,\tilde\gamma(t),\dot{\tilde\gamma}(t)\right)
\bigl(\dot{\tilde\xi}(t), \dot{\tilde \eta}(t)\bigr)  + D_{\tilde
q\tilde v} \tilde
  L\left(t,\tilde\gamma(t), \dot{\tilde\gamma}(t)\right)
\bigl(\tilde\xi(t), \dot{\tilde\eta}(t)\bigr) \nonumber\\
&&\hspace{7mm}  + D_{\tilde v\tilde q} \tilde
  L\left(t,\tilde\gamma(t), \dot{\tilde\gamma}(t)\right)
\bigl(\dot{\tilde\xi}(t), \tilde\eta(t)\bigr) + D_{\tilde q\tilde q}
\tilde L\left(t,\tilde\gamma(t),
\dot{\tilde\gamma}(t)\right)\bigl(\tilde\xi(t),\tilde
\eta(t)\bigr)\Bigr) \, dt\nonumber\\
&&\hspace{4mm} =\bigl(B_{k\tau}(\tilde\gamma)\tilde\xi,
  \tilde\eta\bigr)_{W^{1,2}}\quad\forall \tilde\xi,\tilde\eta\in \tilde H_{k\tau}.
\end{eqnarray}

\begin{lemma}\label{lem:3.4}
$B_{k\tau}(0)$ satisfies the condition ({\bf B1}).
\end{lemma}

\noindent{\bf Proof}. Firstly, by Lemma~\ref{lem:3.1}  there exists
a linear symmetric compact operator $\Xi_{k\tau}:\tilde H_{k\tau}\to
\tilde H_{k\tau}$ such that
$$
\int^{k\tau}_0\tilde\xi(t)\cdot\tilde\eta(t)
dt=(\Xi_{k\tau}\tilde\xi,
\tilde\eta)_{W^{1,2}}\quad\forall\tilde\xi, \tilde\eta\in \tilde
H_{k\tau}.
$$
Note that for sufficiently large $M>0$ there exists a $\delta>0$
such that
$$
M\int^{k\tau}_0\tilde\xi(t)\cdot\tilde\xi(t) dt + d^2
\hat{f}_{k\tau,0}
  (0)(\tilde\xi,\tilde\xi)\ge
  \delta\|\tilde\xi\|^2_{W^{1,2}}\quad\forall\tilde\xi\in \tilde H_{k\tau}.
$$
Hence $R_{k\tau}:=M\Xi_{k\tau} + B_{k\tau}(0): \tilde H_{k\tau}\to
\tilde H_{k\tau}$ is a bounded linear positive definite operator.
Since $C_{k\tau}:=M(R_{k\tau})^{-1}\Xi_{k\tau}$ is compact,
$$
B_{k\tau}(0)=R_{k\tau}- M\Xi_{k\tau}=R_{k\tau}(I-
C_{k\tau})\quad\hbox{implies}
$$
\begin{description}
\item[a)] $0$ is an isolated spectrum point of $B_{k\tau}(0)$,
\item[b)] the maximal negative subspace of $B_{k\tau}(0))$ in
$\tilde H_{k\tau}$ is finitely dimensional and is contained in
$\tilde X_{k\tau}$.
\end{description}
(See  the arguments in \cite[pp.176-177]{Du}).

Next, suppose that $\tilde\xi\in \tilde H_{k\tau}$ satisfies:
$B_{k\tau}(0)\tilde\xi=\tilde\zeta\in \tilde
X_{k\tau}=C^{1}(S_{k\tau}, \R^n)$. We want to prove $\tilde\xi\in
\tilde X_{k\tau}$. To this goal let
$$
J_{k\tau}(s):=\int^s_0\Bigl[\bigl(\hat P_0(t)\dot{\tilde\xi}(t)+
\hat
Q_0(t)\tilde\xi(t)\bigr)-\frac{1}{k\tau}\int^{k\tau}_0\bigl(\hat
P_0(t)\dot{\tilde\xi}(t) + \hat Q_0(t)\tilde\xi(t)\bigr) dt\Bigr]dt
$$
for $s\in\R$.  Since
$$
\int^{k\tau}_0\left[(\hat P_0\dot{\tilde\xi}+ \hat
Q_0\tilde\xi)\cdot\dot{\tilde\eta}+ Q_0^T\dot{\tilde\xi}\cdot
\tilde\eta+ \hat R_0\tilde\xi\cdot
\tilde\eta\right]dt=(\tilde\zeta,\tilde\eta)_{W^{1,2}}
$$
for any $\tilde\eta\in \tilde H_{k\tau}$,
 we have
$$
\int^{k\tau}_0(\hat P_0\dot{\tilde\xi}+ \hat
Q_0\tilde\xi)\cdot\dot{\tilde\eta}\,dt=\int^{k\tau}_0\dot{J}_{k\tau}\cdot{\tilde\eta}\,dt
$$
and thus
$$
\int^{k\tau}_0\left[(Q_0^T\dot{\tilde\xi} + \hat R_0\tilde\xi-
J_{k\tau})\cdot \tilde\eta\right]dt=(\tilde\zeta-
J_{k\tau},\tilde\eta)_{W^{1,2}}.
$$
for any $\tilde\eta\in \tilde H_{k\tau}$. As in the arguments below
(\ref{e:3.5}) we get
\begin{eqnarray*}
\tilde\zeta(t)-
J_{k\tau}(t)&=&\frac{e^t}{2}\int^\infty_te^{-s}\left( J_{k\tau}(s)-
Q_0^T(s)\dot{\tilde\xi}(s)- \hat R_0(s)\tilde\xi(s)\right)\, ds\\
&+&\frac{e^{-t}}{2}\int^t_{-\infty}e^{s}\left(J_{k\tau}(s)-
Q_0^T(s)\dot{\tilde\xi}(s)- \hat R_0(s)\tilde\xi(s) \right)\, ds
\end{eqnarray*}
 for any $t\in\R$. This leads to
\begin{eqnarray*}
 \dot{J}_{k\tau}(t)&=&\dot{\tilde\zeta}(t)-
\frac{e^t}{2}\int^\infty_te^{-s}\left( J_{k\tau}(s)-
Q_0^T(s)\dot{\tilde\xi}(s)- \hat R_0(s)\tilde\xi(s)\right)\, ds\\
&+&\frac{e^{-t}}{2}\int^t_{-\infty}e^{s}\left(J_{k\tau}(s)-
Q_0^T(s)\dot{\tilde\xi}(s)- \hat R_0(s)\tilde\xi(s) \right)\, ds
\end{eqnarray*}
 is continuous in $t\in\R$. Note that $\hat P_0(t)$
 is invertible and
$$
\hat P_0(t)\dot{\tilde\xi}(t)+ \hat
Q_0(t)\tilde\xi(t)=\dot{J}_{k\tau}(t)+
\frac{1}{k\tau}\int^{k\tau}_0(\hat P_0(t)\dot{\tilde\xi}(t) + \hat
Q_0(t)\tilde\xi(t)) dt.
$$
This shows that $\dot{\tilde\xi}(t)$ is continuous in $t$. Hence
$\tilde\xi\in C^1(S_{k\tau},\R^n)$. Lemma~\ref{lem:3.4} is proved.
$\Box$

\begin{lemma}\label{lem:3.5}
$B_{k\tau}$ satisfies the condition ({\bf B2}).
\end{lemma}

\noindent{\bf Proof}. Recall that under the assumptions (L1)-(L3) in
\cite{Lu}  there exist constants $0<c<C$ such that for all $(t, q,
v)\in\R\times B^n_\rho(0)\times\R^n$ the following inequalities
hold.
\begin{eqnarray}
&&|\tilde L(t,q,v)|\le C(1+ |v|^2),\nonumber\\
&&\Bigl| \frac{\partial\tilde L}{\partial q_i}(t,q,v)\Bigr|\le C(1+
|v|^2),\quad\Bigl| \frac{\partial\tilde L}{\partial
v_i}(t,q,v)\Bigr|\le C(1+ |v|),\nonumber\\
&&\Bigl| \frac{\partial^2\tilde L}{\partial q_i\partial
q_j}(t,q,v)\Bigr|\le C(1+ |v|^2),\quad \Bigl| \frac{\partial^2\tilde
L}{\partial q_i\partial v_j}(t,q,v)\Bigr|\le C(1+ |v|),\label{e:3.14}\\
&&\Bigl| \frac{\partial^2\tilde L}{\partial v_i\partial
v_j}(t,q,v)\Bigr|\le C\quad\hbox{and}\quad
\sum_{ij}\frac{\partial^2\tilde L}{\partial v_i\partial
v_j}(t,q,v)u_iu_j\ge c|{\bf u}|^2\label{e:3.15}
\end{eqnarray}
for all ${\bf u}=(u_1,\cdots,u_n)\in\R^n$.

For any $\tilde\gamma\in \tilde V_{k\tau}\cap\tilde X_{k\tau}$, by
(\ref{e:3.3}) and (\ref{e:3.7}) we have
\begin{eqnarray*}
(B_{k\tau}(\tilde\gamma)\tilde\xi,
  \tilde\eta)_{W^{1,2}}
   = \int_0^{k\tau} \Bigl(\!\! \!\!\!&&\!\!\!\!\!D_{\tilde v\tilde v}
    \tilde L\left(t,\tilde\gamma(t),\dot{\tilde\gamma}(t)\right)
\bigl(\dot{\tilde\xi}(t), \dot{\tilde \eta}(t)\bigr) \nonumber\\
&&+ D_{\tilde q\tilde v} \tilde
  L\left(t,\tilde\gamma(t), \dot{\tilde\gamma}(t)\right)
\bigl(\tilde\xi(t), \dot{\tilde\eta}(t)\bigr)\nonumber \\
&& + D_{\tilde v\tilde q} \tilde
  L\left(t,\tilde\gamma(t),\dot{\tilde\gamma}(t)\right)
\bigl(\dot{\tilde\xi}(t), \tilde\eta(t)\bigr) \nonumber\\
&&+  D_{\tilde q\tilde q} \tilde L\left(t,\tilde\gamma(t),
\dot{\tilde\gamma}(t)\right) \bigl(\tilde\xi(t),\tilde
\eta(t)\bigr)\Bigr) \, dt
\end{eqnarray*}
for any $\tilde\xi,\tilde\eta\in C^1(S_{k\tau}, \R^n)$ and $k\in\N$.
Set
\begin{eqnarray*}
&&(P(\tilde\gamma)\tilde\xi,
  \tilde\eta)_{W^{1,2}}
   = \int_0^{k\tau} \Bigl(D_{\tilde v\tilde v}
    \tilde L\left(t,\tilde\gamma(t),\dot{\tilde\gamma}(t)\right)
\bigl(\dot{\tilde\xi}(t), \dot{\tilde \eta}(t)\bigr)+
\bigl(\tilde\xi(t),\tilde \eta(t)\bigr)\Bigr) \, dt,\\
&&(Q_1(\tilde\gamma)\tilde\xi,
  \tilde\eta)_{W^{1,2}}
   = \int_0^{k\tau}D_{\tilde v\tilde q} \tilde
  L\left(t,\tilde\gamma(t),\dot{\tilde\gamma}(t)\right)
\bigl(\dot{\tilde\xi}(t), \tilde\eta(t)\bigr)dt,\\
 &&(Q_2(\tilde\gamma)\tilde\xi,
  \tilde\eta)_{W^{1,2}}
   = \int_0^{k\tau}D_{\tilde q\tilde v} \tilde
  L\left(t,\tilde\gamma(t), \dot{\tilde\gamma}(t)\right)
\bigl(\tilde\xi(t), \dot{\tilde\eta}(t)\bigr)dt, \\
&&(Q_3(\tilde\gamma)\tilde\xi,
  \tilde\eta)_{W^{1,2}}
   = \int_0^{k\tau} \Bigl(D_{\tilde q\tilde q} \tilde L\left(t,\tilde\gamma(t),
\dot{\tilde\gamma}(t)\right) \bigl(\tilde\xi(t),\tilde
\eta(t)\bigr)-\bigl(\tilde\xi(t),\tilde \eta(t)\bigr)\Bigr) \, dt.
\end{eqnarray*}
We shall prove that the operator $P(\tilde\gamma)$ and
$Q(\tilde\gamma):=Q_1(\tilde\gamma)+ Q_2(\tilde\gamma)+
Q_3(\tilde\gamma)$ satisfy the conditions in ({\bf B2}).  It is
clear that (\ref{e:3.15}) implies
$$
\min\{c,1\}\|\tilde\xi\|^2_{W^{1,2}}\le (P(\tilde\gamma)\tilde\xi,
  \tilde\xi)_{W^{1,2}}\le \max\{C, 1\}\|\tilde\xi\|^2_{W^{1,2}}
$$
for any $\tilde\gamma\in\tilde V_{k\tau}\cap X_{k\tau}$ and
$\tilde\xi\in \tilde H_{k\tau}$, and in particular (iii) of ({\bf
B2}) for $P$. It remains to prove that the conditions (i)-(ii) in
({\bf B2}) are satisfied.

For any $\tilde\gamma, \tilde\alpha\in \tilde V_{k\tau}\cap\tilde
X_{k\tau}$ and $\tilde\xi, \tilde\eta\in\tilde H_{k\tau}$,
\begin{eqnarray*}
&&(P(\tilde\gamma)\tilde\xi- P(\tilde\alpha)\tilde\xi,
  \tilde\eta)_{W^{1,2}}\\
&=&\int_0^{k\tau} \Bigl(D_{\tilde v\tilde v}
    \tilde L\left(t,\tilde\gamma(t),\dot{\tilde\gamma}(t)\right)
\bigl(\dot{\tilde\xi}(t), \dot{\tilde \eta}(t)\bigr)- D_{\tilde
v\tilde v}
    \tilde L\left(t,\tilde\alpha(t),\dot{\tilde\alpha}(t)\right)
\bigl(\dot{\tilde\xi}(t), \dot{\tilde \eta}(t)\bigr)\Bigr) \, dt,\\
&=&\int_0^{k\tau}\sum^n_{j=1}\left[\sum^n_{i=1}\left(\frac{\partial^2\tilde
L}{\partial v_i\partial
v_j}\left(t,\tilde\gamma(t),\dot{\tilde\gamma}(t)\right)-
\frac{\partial^2\tilde L}{\partial v_i\partial
v_j}\left(t,\tilde\alpha(t),\dot{\tilde\alpha}(t)\right)\right)
\dot{\tilde\xi}_i(t)\right]\cdot\dot{\tilde \eta}_j(t) \, dt.
\end{eqnarray*}
It follows that
\begin{eqnarray*}
&&\|P(\tilde\gamma)\tilde\xi- P(\tilde\alpha)\tilde\xi\|_{W^{1,2}}\\
&\le
&\left(\int_0^{k\tau}\sum^n_{j=1}\left|\sum^n_{i=1}\left(\frac{\partial^2\tilde
L}{\partial v_i\partial
v_j}\left(t,\tilde\gamma(t),\dot{\tilde\gamma}(t)\right)-
\frac{\partial^2\tilde L}{\partial v_i\partial
v_j}\left(t,\tilde\alpha(t),\dot{\tilde\alpha}(t)\right)\right)
\dot{\tilde\xi}_i(t)\right|^2\, dt\right)^{1/2}.
\end{eqnarray*}
Thus we only need to prove
$$
\int_0^{k\tau}\sum^n_{j=1}\left|\sum^n_{i=1}\left(\frac{\partial^2\tilde
L}{\partial v_i\partial
v_j}\left(t,\tilde\gamma(t),\dot{\tilde\gamma}(t)\right)-
\frac{\partial^2\tilde L}{\partial v_i\partial
v_j}\left(t,\tilde\alpha(t),\dot{\tilde\alpha}(t)\right)\right)
\dot{\tilde\xi}_i(t)\right|^2\, dt\to 0
$$
as $\tilde\gamma\in\tilde V_{k\tau}\cap\tilde X_{k\tau}$ and
$\|\tilde\gamma-\tilde\alpha\|_{W^{1,2}}\to 0$.  By a contradiction
suppose that there exist $c_0>0$ and a sequence
$\{\tilde\gamma_m\}\subset\tilde V_{k\tau}\cap\tilde X_{k\tau}$ with
$\|\tilde\gamma_m-\tilde\alpha\|_{W^{1,2}}\to 0$ as $m\to\infty$,
such that
$$
\int_0^{k\tau}\sum^n_{j=1}\left|\sum^n_{i=1}\left(\frac{\partial^2\tilde
L}{\partial v_i\partial
v_j}\left(t,\tilde\gamma_m(t),\dot{\tilde\gamma}_m(t)\right)-
\frac{\partial^2\tilde L}{\partial v_i\partial
v_j}\left(t,\tilde\alpha(t),\dot{\tilde\alpha}(t)\right)\right)
\dot{\tilde\xi}_i(t)\right|^2\, dt\ge c_0
$$
for all $m=1,2,\cdots$. Using Proposition~1.2 in \cite{MW} and
Corollary~2.17 in \cite{AdFo} we have a subsequence
$\{\tilde\gamma_{m_l}\}$ such that $\{\tilde\gamma_{m_l}\}$
converges uniformly to $\tilde\alpha$ on $[0, k\tau]$ and that
$\{\dot{\tilde\gamma}_{m_l}\}$ converges pointwise almost everywhere
 to $\dot{\tilde\alpha}$ on $[0, k\tau]$. Since  (\ref{e:3.15})
 implies
$$
\left|\frac{\partial^2\tilde L}{\partial v_i\partial
v_j}\left(t,\tilde\gamma_{m_l}(t),\dot{\tilde\gamma}_{m_l}(t)\right)-
\frac{\partial^2\tilde L}{\partial v_i\partial
v_j}\left(t,\tilde\alpha(t),\dot{\tilde\alpha}(t)\right)\right|\le
2C\quad\forall l=1,2,\cdots,
$$
 the dominated convergence theorem leads to
$$
\int_0^{k\tau}\sum^n_{j=1}\left|\sum^n_{i=1}\left(\frac{\partial^2\tilde
L}{\partial v_i\partial
v_j}\left(t,\tilde\gamma_{m_l}(t),\dot{\tilde\gamma}_{m_l}(t)\right)-
\frac{\partial^2\tilde L}{\partial v_i\partial
v_j}\left(t,\tilde\alpha(t),\dot{\tilde\alpha}(t)\right)\right)
\dot{\tilde\xi}_i(t)\right|^2\, dt\to 0
$$
as $l\to\infty$. This contradiction affirms (i) of ({\bf B2}) for
$P$.

In order to prove (ii) of ({\bf B2}) for $Q$, we only need to prove
that each one of the operators $Q_1, Q_2$ and $Q_3$ satisfies (ii)
of ({\bf B2}). Viewing $D_{\tilde v\tilde q} \tilde
  L\left(t,\tilde\gamma(t),\dot{\tilde\gamma}(t)\right)$ as a matrix
  of order $n\times n$, we can write
$$
(Q_1(\tilde\gamma)\tilde\xi,
  \tilde\eta)_{W^{1,2}}
   = \int_0^{k\tau}\bigl(D_{\tilde v\tilde q} \tilde
  L\left(t,\tilde\gamma(t),\dot{\tilde\gamma}(t)\right)
\dot{\tilde\xi}(t), \tilde\eta(t)\bigr)_{\R^n}dt.
$$
As the arguments below Lemma~\ref{lem:3.1} we have
\begin{eqnarray*}
(Q_1(\tilde\gamma)\tilde\xi)(t)&=&-\frac{e^t}{2}\int^\infty_te^{-s}\left(
D_{\tilde v\tilde q} \tilde
  L\left(s,\tilde\gamma(s),\dot{\tilde\gamma}(s)\right)
\dot{\tilde\xi}(s)\right)\, ds\\
&-&\frac{e^{-t}}{2}\int^t_{-\infty}e^{s}\left(D_{\tilde v\tilde q}
\tilde
  L\left(s,\tilde\gamma(s),\dot{\tilde\gamma}(s)\right)
\dot{\tilde\xi}(s)\right)\, ds ,\\
\frac{d}{dt}(Q_1(\tilde\gamma)\tilde\xi)(t)&=&-\frac{e^t}{2}\int^\infty_te^{-s}\left(
D_{\tilde v\tilde q} \tilde
  L\left(s,\tilde\gamma(s),\dot{\tilde\gamma}(s)\right)
\dot{\tilde\xi}(s)\right)\, ds\\
&+&\frac{e^{-t}}{2}\int^t_{-\infty}e^{s}\left(D_{\tilde v\tilde q}
\tilde
  L\left(s,\tilde\gamma(s),\dot{\tilde\gamma}(s)\right)
\dot{\tilde\xi}(s)\right)\, ds.
\end{eqnarray*}
As in the proof of Lemma~\ref{lem:3.4} it easily follows that
$Q_1(\tilde\gamma)$ is a completely continuous operator from $\tilde
H_{k\tau}$ to $\tilde H_{k\tau}$. Now  for $0\le t\le k\tau$ it is
not hard to check that
\begin{eqnarray*}
&&|(Q_1(\tilde\gamma)\tilde\xi)(t)-
(Q_1(\tilde\alpha)\tilde\xi)(t)|\\
&\le & \frac{e^{k\tau}}{2}\int^\infty_0e^{-s}\left|D_{\tilde v\tilde
q} \tilde
  L(s,\tilde\gamma(s),\dot{\tilde\gamma}(s))
- D_{\tilde v\tilde q} \tilde
  L(s,\tilde\alpha(s),\dot{\tilde\alpha}(s)\right||\dot{\tilde\xi}(s)|\, ds\\
&+&\frac{1}{2}\int^{k\tau}_{-\infty}e^{s}\left|D_{\tilde v\tilde q}
\tilde
  L(s,\tilde\gamma(s),\dot{\tilde\gamma}(s))-D_{\tilde v\tilde q}
\tilde
  L(s,\tilde\alpha(s),\dot{\tilde\alpha}(s))\right|
|\dot{\tilde\xi}(s)|\, ds ,\\
&=&\frac{e^{k\tau}}{2}\sum^\infty_{i=0}\int^{(i+1)k\tau}_{ik\tau}e^{-s}\left|D_{\tilde
v\tilde q} \tilde
  L(s,\tilde\gamma(s),\dot{\tilde\gamma}(s))
- D_{\tilde v\tilde q} \tilde
  L(s,\tilde\alpha(s),\dot{\tilde\gamma}_0(s)\right||\dot{\tilde\xi}(s)|\, ds\\
&+&\frac{1}{2}\sum^1_{i=-\infty}\int^{ik\tau}_{(i-1)k\tau}e^{s}\left|D_{\tilde
v\tilde q} \tilde
  L(s,\tilde\gamma(s),\dot{\tilde\gamma}(s))-D_{\tilde v\tilde q}
\tilde
  L(s,\tilde\alpha(s),\dot{\tilde\alpha}(s))\right|
|\dot{\tilde\xi}(s)|\, ds \\
&\le&\frac{e^{k\tau}}{2}\sum^\infty_{i=0}e^{-ik\tau}\int^{k\tau}_{0}\left|D_{\tilde
v\tilde q} \tilde
  L(s,\tilde\gamma(s),\dot{\tilde\gamma}(s))
- D_{\tilde v\tilde q} \tilde
  L(s,\tilde\alpha(s),\dot{\tilde\alpha}(s)\right||\dot{\tilde\xi}(s)|\, ds\\
&+&\frac{1}{2}\sum^1_{i=-\infty}e^{ik\tau}\int^{k\tau}_{0}\left|D_{\tilde
v\tilde q} \tilde
  L(s,\tilde\gamma(s),\dot{\tilde\gamma}(s))-D_{\tilde v\tilde q}
\tilde
  L(s,\tilde\alpha(s),\dot{\tilde\alpha}(s))\right|
|\dot{\tilde\xi}(s)|\, ds\\
&\le & (e^{k\tau}+ 1)\int^{k\tau}_{0}\left|D_{\tilde v\tilde q}
\tilde
  L(s,\tilde\gamma(s),\dot{\tilde\gamma}(s))-D_{\tilde v\tilde q}
\tilde
  L(s,\tilde\alpha(s),\dot{\tilde\alpha}(s))\right|
|\dot{\tilde\xi}(s)|\, ds \\
&\le&(e^{k\tau}+ 1)\left(\int^{k\tau}_{0}\left|D_{\tilde v\tilde q}
\tilde
  L(s,\tilde\gamma(s),\dot{\tilde\gamma}(s))-D_{\tilde v\tilde q}
\tilde
  L(s,\tilde\alpha(s),\dot{\tilde\alpha}(s))\right|^2
  ds\right)^{1/2}\|\tilde\xi\|_{W^{1,2}}.
\end{eqnarray*}
Similarly, we have
\begin{eqnarray*}
&&\left|\frac{d}{dt}(Q_1(\tilde\gamma)\tilde\xi)(t)-
\frac{d}{dt}(Q_1(\tilde\gamma_0)\tilde\xi)(t)\right|\\
&\le&(e^{k\tau}+ 1)\left(\int^{k\tau}_{0}\left|D_{\tilde v\tilde q}
\tilde
  L(s,\tilde\gamma(s),\dot{\tilde\gamma}(s))-D_{\tilde v\tilde q}
\tilde
  L(s,\tilde\alpha(s),\dot{\tilde\alpha}(s))\right|^2
  ds\right)^{1/2}\|\tilde\xi\|_{W^{1,2}}.
\end{eqnarray*}
It follows that
\begin{eqnarray*}
&&\|Q_1(\tilde\gamma)-
Q_1(\tilde\gamma_0)\|_{L(\tilde H_{k\tau})}\\
&\le&2(e^{k\tau}+ 1)\left(\int^{k\tau}_{0}\left|D_{\tilde v\tilde q}
\tilde
  L(s,\tilde\gamma(s),\dot{\tilde\gamma}(s))-D_{\tilde v\tilde q}
\tilde
  L(s,\tilde\alpha(s),\dot{\tilde\alpha}(s))\right|^2
  ds\right)^{1/2}.
\end{eqnarray*}
By a theorem of Krasnosel'skii, the second inequality in
(\ref{e:3.14}) implies that the map
$$
W^{1,2}(S_{k\tau}, \R^n)\to
L^2(S_{k\tau},\R^n),\;\tilde\gamma\mapsto D_{\tilde v\tilde q}
\tilde
  L(\cdot,\tilde\gamma(\cdot),\dot{\tilde\gamma}(\cdot))
$$
is continuous. Hence we have proved that $Q_1$ satisfies the
condition (ii) of ({\bf B2}).

 Viewing $D_{\tilde q\tilde v} \tilde
  L\left(t,\tilde\gamma(t), \dot{\tilde\gamma}(t)\right)$ as a matrix
  of order $n\times n$, we can write
$$
(Q_2(\tilde\gamma)\tilde\xi,
  \tilde\eta)_{W^{1,2}}
   = \int_0^{k\tau}\bigl(\tilde\xi(t), [D_{\tilde q\tilde v} \tilde
  L\left(t,\tilde\gamma(t), \dot{\tilde\gamma}(t)\right)]^T\dot{\tilde\eta}(t)\bigr)_{\R^n}dt=
  \bigr(\tilde\xi,  (Q_2(\tilde\gamma))^\ast\tilde\eta\bigl)_{W^{1,2}}.
$$
Hence $(Q_2(\tilde\gamma))^\ast$ and thus $Q_2(\tilde\gamma)$
satisfies (ii) of ({\bf B2}).

Finally, for $Q_3$ we view $D_{\tilde q\tilde q} \tilde
L\left(t,\tilde\gamma(t), \dot{\tilde\gamma}(t)\right)$ as a matrix
  of order $n\times n$, we can write
$$
 (Q_3(\tilde\gamma)\tilde\xi,
  \tilde\eta)_{W^{1,2}}
   = \int_0^{k\tau} \bigl(D_{\tilde q\tilde q} \tilde L\bigl(t,\tilde\gamma(t),
\dot{\tilde\gamma}(t)\bigr)
\tilde\xi(t)-\tilde\xi(t),\eta(t)\bigr)_{\R^n} \, dt.
$$
As before we can get
\begin{eqnarray*}
(Q_3(\tilde\gamma)\tilde\xi)(t)&=&-\frac{e^t}{2}\int^\infty_te^{-s}\left(
D_{\tilde q\tilde q} \tilde L\bigl(s,\tilde\gamma(s),
\dot{\tilde\gamma}(s)\bigr) \tilde\xi(s)-\tilde\xi(s)\right)\, ds\\
&-&\frac{e^{-t}}{2}\int^t_{-\infty}e^{s}\left( D_{\tilde q\tilde q}
\tilde L\bigl(s,\tilde\gamma(s), \dot{\tilde\gamma}(s)\bigr)
\tilde\xi(s)-\tilde\xi(s)\right)\, ds ,\\
\frac{d}{dt}(Q_3(\tilde\gamma)\tilde\xi)(t)&=&-\frac{e^t}{2}\int^\infty_te^{-s}\left(
D_{\tilde q\tilde q} \tilde L\bigl(s,\tilde\gamma(s),
\dot{\tilde\gamma}(s)\bigr) \tilde\xi(s)-\tilde\xi(s)
\right)\, ds\\
&+&\frac{e^{-t}}{2}\int^t_{-\infty}e^{s}\left( D_{\tilde q\tilde q}
\tilde L\bigl(s,\tilde\gamma(s), \dot{\tilde\gamma}(s)\bigr)
\tilde\xi(s)-\tilde\xi(s) \right) \, ds.
\end{eqnarray*}
As for $Q_1$ we can use this to prove that the condition  (ii) of
({\bf B2}) holds for $Q_3$. Hence $Q$ satisfies (ii) of ({\bf B2}).
Lemma~\ref{lem:3.5} is proved. $\Box$\vspace{2mm}

Similar to the proofs of Lemmas~\ref{lem:3.4},~\ref{lem:3.5} we may
prove that $B^E_{k\tau}$ satisfies the conditions ({\bf B1}) and
({\bf B2}).

According to the definition above Theorem~\ref{th:1.1},
  the \verb"nullity" $m^0_{k\tau}(\tilde{\cal L}_{k\tau},
0)$ and \verb"Morse index" $m^-_{k\tau}(\tilde{\cal L}_{k\tau}, 0)$
of $\tilde{\cal L}_{k\tau}$ at $0\in \tilde V_{k\tau}$ are
respectively given by $\dim{\rm Ker}(B_{k\tau}(0))$ and  the
dimension of the negative definite space of $B_{k\tau}(0))$ in
$\tilde H_{k\tau}$. (\ref{e:3.13}) implies
\begin{equation}\label{e:3.16}
  m^-_{k\tau}(\tilde{\cal L}_{k\tau}, 0)=
 m^-_{k\tau}(\hat f_{k\tau}, 0) \quad{\rm and}\quad
 m^0_{k\tau}(\tilde{\cal L}_{k\tau}, 0)=m^0_{k\tau}(\hat f_{k\tau},
0).
\end{equation}
We define the  Morse index and nullity of ${\cal L}_{k\tau}$ at
$\gamma_0\in H_{k\tau}$ by
\begin{equation}\label{e:3.17}
  m^-_{k\tau}(\gamma_0):= m^-_{k\tau}(\tilde{\cal L}_{k\tau}, 0) \quad{\rm and}\quad
 m^0_{k\tau}(\gamma_0):=m^0_{k\tau}(\tilde{\cal L}_{k\tau}, 0).
\end{equation}
It is easily checked that they are the maximal value of dimensions
of linear subspaces $L\subset T_{\gamma_0} X_{k\tau}$ on which
 the second-order differential $d^2{\cal L}^X_{k\tau}(\gamma_0)<0$ and $d^2{\cal
L}^X_{k\tau}(\gamma_0)=0$, respectively.
 It follows from (\ref{e:3.16}) that
\begin{equation}\label{e:3.18}
 0\le m^0_{k\tau}(\gamma_0)=
m^0_{k\tau}(\hat f_{k\tau}, 0)\le 2n.
\end{equation}

Let $\Psi:[0, +\infty)\to{\rm Sp}(2n,\R)$ be the fundamental
solution of the problem $\dot{\bf u}(t)=J_0\hat S(t){\bf u}$ with
$\Psi(0)=I_{2n}$, and let $i_{k\tau}(\Psi)$ and $\nu_{k\tau}(\Psi)$
be the Maslov-type index of $\Psi$ on $[0, k\tau]$. Here
$$
\hat S(t)=\left(\begin{array}{cc}
\hat P_0(t)^{-1}& -\hat P_0(t)^{-1}\hat Q_0(t)\\
-\hat Q_0(t)^T\hat P_0(t)^{-1}& \hat Q_0(t)^T\hat P_0(t)^{-1}\hat
Q_0(t)-\hat R_0(t)\end{array}\right),
$$
and $\hat P_0$, $\hat Q_0$ and $\hat R_0$ are defined as below
(\ref{e:2.8}).  By \cite[(2.16)]{Lu}, for any $k\in\N$ we have
\begin{equation}\label{e:3.19}
m^-_{k\tau}(\hat f_{k\tau}, 0)=i_{k\tau}(\Psi)\quad{\rm and}\quad
m^0_{k\tau}(\hat f_{k\tau}, 0)=\nu_{k\tau}(\Psi)
\end{equation}
According to the definitions of the Morse index and nullity for  a
critical point $\gamma$ of ${\cal L}_\tau$ on $H_\tau$ in
(\ref{e:3.17}), from (\ref{e:3.19}) and Lemma~2.1 of \cite{Lu} we
deduce that  Theorem 3.1 in \cite{Lu} is still true, i.e.

\begin{theorem}\label{th:3.6}
 For  a critical point $\gamma$ of ${\cal
L}_\tau$ on $H_\tau$, assume that $\gamma^\ast TM\to S_\tau$ is
trivial. Then the mean Morse index
$$
\hat{m}^-_{\tau}(\gamma):=\lim_{k\to\infty}\frac{{m}^-_{k\tau}(\gamma^k)}{k}
$$
always exists, and it holds that
\begin{equation}\label{e:3.20}
{\rm max}\left\{0, k\hat {m}^-_{\tau}(\gamma)-n\right\}\le
{m}^-_{k\tau}(\gamma^k)\le k\hat {m}_{\tau}(\gamma_0)+n-
{m}^0_{k\tau}(\gamma^k) \quad\forall k\in\N.
\end{equation}
Consequently, for any critical point $\gamma$ of ${\cal L}_\tau$ on
$H_\tau$, $\hat {m}^-_{2\tau}(\gamma^2)$ exists and
\begin{equation}\label{e:3.21}
{\rm max}\left\{0, k\hat {m}^-_{2\tau}(\gamma^2)-n\right\}\le
{m}^-_{2k\tau}(\gamma^{2k})\le k\hat {m}_{2\tau}(\gamma)+ n-
{m}^0_{2k\tau}(\gamma^{2k}) \quad\forall k\in\N
\end{equation}
because $(\gamma^2)^\ast TM\to S_{2\tau}$ is always trivial.
\end{theorem}

Similarly, for a critical point $\gamma_0$ of ${\cal L}^E_{k\tau}$
on $EH_{k\tau}$ we may define the Morse index and nullity of it,
\begin{equation}\label{e:3.22}
 m^-_{1,k\tau}(\gamma)\quad\hbox{and}\quad
m^0_{1,k\tau}(\gamma).
\end{equation}
They are equal to the maximal value of dimensions of linear
subspaces $S\subset T_{\gamma_0} EX_{k\tau}$ on which $d^2{\cal
L}^{EX}_{k\tau}(\gamma_0)<0$ and $d^2{\cal
L}^{EX}_{k\tau}(\gamma_0)=0$, respectively. Here
$$
EX_{k\tau}=\bigl\{\gamma\in  X_{k\tau}\,|\,
\gamma(-t)=\gamma(t)\;\forall t\in\R\bigr\}
$$
and ${\cal L}^{EX}_{k\tau}$ is the restriction of ${\cal
L}^E_{k\tau}$ to $EX_{k\tau}$. Then $0\le m^0_{1,k\tau}(\gamma)\le
2n$ for any $k\in\N$. According to the definitions of the Morse
index  and nullity  in (\ref{e:3.22}) Theorem 3.3 in \cite{Lu} also
 holds, i.e.

\begin{theorem}\label{th:3.7}
 Let $L$ satisfy the conditions (L1)-(L4). Then for any critical point $\gamma$ of ${\cal
L}^E_\tau$ on $EH_\tau$, the mean Morse index
\begin{equation}\label{e:3.23}
\hat{m}^-_{1,\tau}(\gamma):=\lim_{k\to\infty}\frac{{m}^-_{1,
k\tau}(\gamma^k)}{k}
\end{equation}
 exists, and it holds that
\begin{equation}\label{e:3.24}
{m}^-_{1, k\tau}(\gamma^k)+ {m}^0_{1, k\tau}(\gamma^k)\le n
\quad\forall k\in\N \quad{\rm if}\quad \hat
{m}^-_{1,\tau}(\gamma)=0.
\end{equation}
\end{theorem}

\section{The corrections of Sections~4.1, 4.2 in
\cite{Lu}}\label{sec:4} \setcounter{equation}{0}

Though Theorems~4.4, 4.7 in \cite{Lu} are now direct consequences of
Theorem~\ref{th:2.1} we are also to use Theorem~\ref{th:1.1} to
revise their proofs in \cite{Lu}. (So far Lemma~4.2 of \cite{Lu} is
not needed). We shall use the same count for equations as that of
\cite{Lu}.

\subsection{The corrections of  Section~4.1 in \cite{Lu}}\label{sec:4.1}

 We only need to change  the part before \cite[Lemma
4.5]{Lu}  in the second passage  on the page 2994 of \cite{Lu} into:

 For later conveniences we outline the arguments therein. By the proof of Theorem~3.1 in \cite{Lu}
for $l\in \N$ we may choose the chart $\phi_{l\tau}$ therein so that
$$
\tilde\gamma^l=(\phi_{l\tau})^{-1}(\gamma^l)=0\in
W^{1,2}(S_{l\tau},\R^n).
$$
Sometimes, for clearness we write $(\phi_{l\tau})^{-1}(\gamma^l)$ as
$\tilde\gamma^l$ rather than $0$.   Let
\begin{eqnarray*}
W^{1,2}(S_{l\tau}, \R^n)=M^0(\tilde\gamma^l)\oplus
M(\tilde\gamma^l)^-\oplus M(\tilde\gamma^l)^+
=M^0(\tilde\gamma^l)\oplus M(\tilde\gamma^l)^\bot
\end{eqnarray*}
be the orthogonal decomposition of  the space $W^{1,2}(S_{l\tau},
\R^n)$  according to the null, negative, and positive definiteness
of the operator $B_{l\tau}(0)$, where
$B_{l\tau}(0)=B_{l\tau}(\tilde\gamma^l)$ is given as above
(\ref{e:3.7}).  By Theorem~\ref{th:1.4} we have homeomorphisms
$\tilde\Theta_{l\tau}$ from some open neighborhoods $\tilde
U_{l\tau}$ of $0$ in $W^{1,2}(S_{l\tau}, \R^n)$ to
$\tilde\Theta_{l\tau}(\tilde U_{l\tau})\subset W^{1,2}(S_{l\tau},
\R^n)$ with $\tilde\Theta_{l\tau}(0)=\tilde\gamma^l=0$, and  $C^1$
maps
$$
\tilde h_{l\tau}: \tilde U_{l\tau}\cap M(\tilde\gamma^l)^0\to
M(\tilde\gamma^l)^\bot\cap\tilde X_{l\tau}
$$
such that
$$
\qquad\tilde {\cal L}_{l\tau}(\tilde\Theta_{l\tau}(\eta+\xi))=\tilde
{\cal L}_{l\tau}\bigl(\eta+ \tilde h_{l\tau}(\eta)\bigr) +
\|\xi^+\|_{\tilde H_{l\tau}}^2-\|\xi^-\|_{\tilde H_{l\tau}}^2\equiv
\tilde\alpha_{l\tau}(\eta)+
\tilde\beta_{l\tau}(\xi)\hspace{6mm}\hbox{(4.12)}
$$
for any $\eta+\xi\in \tilde U_{l\tau}\cap
\left(M(\tilde\gamma^l)^0\oplus M(\tilde\gamma^l)^\bot\right)$ and
$l\in\N$.  By the constructions of the maps $A_{l\tau}$ in
Lemma~\ref{lem:3.2} and $B_{l\tau}$ above (\ref{e:3.7}) it is easily
checked that
$$
\hspace{12mm} \tilde\psi^l\bigl(A_{l\tau}(x)\bigr)=
A_{l\tau}\bigl(\tilde\psi^l(x)\bigr)\quad{\rm
and}\quad\tilde\psi^l\bigl(B_{l\tau}(x)\xi\bigr)=
B_{l\tau}\bigl(\tilde\psi^l(x)\bigr)\tilde\psi^l(\xi)\hspace{10mm}\hbox{(4.13)}
$$
for any $\tau, l\in \N$, $x\in W^{1,2}\bigl(S_\tau,
B^n_\rho(0)\bigr)\cap \tilde X_\tau$ and $\xi\in W^{1,2}(S_\tau,
\R^n)$. Now for $l=1,k$, using Corollary~\ref{cor:1.5} and shrinking
$\tilde U_\tau$ (if necessary) we may get
$$
\tilde\Theta_{k\tau}\circ\tilde\psi^k(\eta+\xi)=\tilde\Theta_{\tau}(\eta+\xi)
$$
and
$$
\hspace{22mm} \tilde\alpha_{k\tau}(\tilde\psi^k(\eta))=
k\tilde\alpha(\eta)\quad{\rm
and}\quad\tilde\beta_{k\tau}(\tilde\psi^k(\xi))=
k\tilde\beta_\tau(\xi)\hspace{28mm}\hbox{(4.14)}
$$
for any $\eta\in\tilde U_\tau\cap M^0(\tilde\gamma)$ and
$\xi\in\tilde U_\tau\cap M^\bot(\tilde\gamma)$.\

\subsection{The corrections of  Section~4.2 in \cite{Lu}}\label{sec:4.2}

We only need to change  the part between line 4 from below on the
page 2998 of \cite{Lu} and  \cite[(4.43)]{Lu} into:
\begin{eqnarray*}
EW^{1,2}(S_{k\tau}, \R^n)=M^0(\tilde\gamma^k)_E\oplus
M(\tilde\gamma^k)^-_E\oplus M(\tilde\gamma^k)^+_E =
M^0(\tilde\gamma^k)_E\oplus M(\tilde\gamma^k)^\bot_E
\end{eqnarray*}
be the orthogonal decomposition of  the space $EW^{1,2}(S_{k\tau},
\R^n)$  according to the null, negative, and positive definiteness
of the operator $B^E_{k\tau}(\tilde\gamma^k)=B^E_{k\tau}(0)$ in
(\ref{e:3.10}). As above Theorem~\ref{th:1.1} yields a homeomorphism
$\tilde\Theta^E_{k\tau}$ from some open neighborhood $\tilde
U^E_{k\tau}$ of $0$ in $EW^{1,2}(S_{k\tau}, \R^n)$ to
$\tilde\Theta^E_{k\tau}(\tilde U^E_{k\tau})\subset
EW^{1,2}(S_{k\tau}, \R^n)$ with
$\tilde\Theta^E_{k\tau}(0)=\tilde\gamma^k=0$, and a $C^1$ map
$\tilde h^E_{k\tau}: \tilde U^E_{k\tau}\cap M(\tilde\gamma^k)^0_E\to
M(\tilde\gamma^k)_E^\bot\cap E\tilde X_{k\tau}$ such that
$$
\begin{array}{rcl}
\tilde {\cal
L}^E_{k\tau}\bigl(\tilde\Theta^E_{k\tau}(\eta+\xi)\bigr)=\tilde
{\cal L}^E_{k\tau}\bigl(\eta+ \tilde h^E_{k\tau}(\eta)\bigr) +
\|\xi^+\|_{E\tilde H_{k\tau}}-\|\xi^-\|_{E\tilde H_{k\tau}}\equiv
\tilde\alpha^E_{k\tau}(\eta)+ \tilde\beta^E_{k\tau}(\xi)
\end{array}
$$
for any $\eta+\xi\in \tilde U^E_{k\tau}\cap
(M(\tilde\gamma^k)_E^0\oplus M(\tilde\gamma^k)_E^\bot)$, where
$\tilde\beta^E_{k\tau}$ and $\tilde\alpha^E_{k\tau}$ are
 $C^\infty$ and $C^2$ respectively.
For a fixed $k\in\N$,  by Corollary~\ref{cor:1.5} (and shrinking
$\tilde U_\tau$ if necessary) we may require
\begin{eqnarray*}
&&\tilde\Theta^E_{k\tau}\circ\tilde\psi^k(\eta+\xi)=\tilde\Theta^E_{\tau}(\eta+\xi)
\quad\hbox{and}\\
&& \tilde\alpha^E_{k\tau}(\tilde\psi^k(\eta))=
k\tilde\alpha^E(\eta)\quad{\rm
and}\quad\tilde\beta^E_{k\tau}(\tilde\psi^k(\xi))=
k\tilde\beta^E_\tau(\xi)
\end{eqnarray*}
for any $\eta\in\tilde U^E_\tau\cap M^0(\tilde\gamma)_E$ and
$\xi\in\tilde U^E_\tau\cap M^\bot(\tilde\gamma)_E$. Clearly,
$\tilde\Theta^E_{k\tau}$ induces isomorphisms on critical modules,
$$
\hspace{30mm}(\tilde\Theta^E_{k\tau})_\ast:C_\ast(\tilde\alpha^E_{k\tau}+
\tilde\beta^E_{k\tau}, 0;\K)\cong C_\ast(\tilde{\cal L}^E_{k\tau},
\tilde\gamma^k;\K).\hspace{20mm}\hbox{(4.41)}
$$
Note that
$$
\hspace{20mm}(W(\gamma^k)_E,
W^-(\gamma^k)_E):=\left(\phi^E_{k\tau}\bigl(\tilde
W(\tilde\gamma^k)_E\bigr), \phi^E_{k\tau}\bigl(\tilde
W^-(\gamma^k)_E\bigr)\right) \hspace{10mm}\hbox{(4.42)}
$$
 is  a Gromoll-Meyer pair of ${\cal L}^E_{k\tau}$ at $\gamma^k$.
 Define the  critical modules
$$
\hspace{30mm} C_\ast({\cal L}^E_{k\tau}, \gamma^k;
\K):=H_\ast\bigl(W(\gamma^k)_E, W^-(\gamma^k)_E;
 \K\bigr).\hspace{20mm}\hbox{(4.43)}
$$

\section{The corrections of  Section~4.3 in \cite{Lu}}\label{sec:5}

In this section we shall rewrite Section~4.3 of \cite{Lu} with some
corrections. One direct method is to follow the original line with
the theory developed in \cite[\S 3]{Lu2} (as in Section 3 or
\cite{LoLu}). We here choose another way for which
Theorem~\ref{th:1.1} is sufficient.

We always assume: $M$ is $C^5$-smooth, $L$ is $C^4$-smooth and
satisfies (L1)-(L3) in \cite{Lu}. The goal is to generalize
\cite[Th.2.5]{LoLu} to the present general case. However, unlike the
last two cases we cannot choose a local coordinate chart around a
critical orbit. For $\tau>0$, let $
S_\tau:=\R/\tau\Z=\{[s]_{\tau}\,|\,[s]_{\tau}=s+\tau\Z,\, s\in\R\}$.
 By Section~2.2 of Chapter 2 in \cite{Kl},  there
exist  equivariant and also isometric operations of
$S_{\tau}$-action on $H_{\tau}(\alpha)$ and $TH_\tau(\alpha)$:
\begin{equation}\label{e:5.1}
\left.\begin{array}{ll} &[s]_{\tau}\cdot \gamma(t)=\gamma(s+t),
\quad \forall [s]_{\tau}\in S_{\tau},
         \; \gamma\in H_{\tau}(\alpha),\\
 &[s]_{\tau}\cdot \xi(t)=\xi(s+t),
\quad \forall [s]_{\tau}\in S_{\tau},
         \; \xi\in T_{\gamma}H_{\tau}(\alpha)
 \end{array}\right\}
         \end{equation}
which are  continuous, but not differentiable. Clearly, ${\cal
L}_{\tau}$ is invariant under this action. Since under our
assumptions each critical point $\gamma$ of ${\cal L}_\tau$ is
$C^4$-smooth,  the orbit $S_{\tau}\cdot \gamma$ is a
$C^3$-submanifold in $H_{\tau}(\alpha)$ by \cite[page 499]{GM2}.
 It is easily checked that $S_{\tau}\cdot
\gamma$ is a $C^3$-smooth critical submanifold of ${\cal L}_\tau$.
Seemingly, the theory of \cite{Wa} cannot be applied to this case
because the action of $S_\tau$ is only continuous. However, as
pointed out in the second paragraph of \cite[page 500]{GM2} this
theory still hold since critical orbits are smooth and $S_\tau$ acts
by isometries.

 For any
$k\in\N$, there is a natural $k$-fold cover $\varphi_k$ from
$S_{k\tau}$ to $S_{\tau}$ defined by
\begin{equation}\label{e:5.2}
\varphi_k: [s]_{k\tau}\mapsto [s]_{\tau}.
 \end{equation}
 It is easy to check that the $S_{\tau}$-action
on $H_{\tau}(\alpha)$, the $S_{k\tau}$-action on
$H_{k\tau}(\alpha^k)$, and the $k$-th iteration map $\psi^k$ defined
above \cite[(3.9)]{Lu} satisfy:
\begin{equation}\label{e:5.3}
\left.\begin{array}{rcl}
  ([s]_{\tau}\cdot \gamma)^k &=& [s]_{k\tau}\cdot \gamma^k,  \\
  {\cal L}_{k\tau}([s]_{k\tau}\cdot \gamma^k) &=& k{\cal L}_{\tau}([s]_{\tau}\cdot \gamma)
     = k{\cal L}_{\tau}(\gamma)
\end{array}\right\}
 \end{equation}
 for all $\gamma\in H_{\tau}(\alpha)$, $k\in\N$,  and $s\in \R$.

 Let $\gamma_0\in H_{\tau}(\alpha)$ be a  critical point of ${\cal L}_{\tau}$.  Denote by ${\cal
O}=S_{\tau}\cdot \gamma_0$. If $\gamma_0$ is nonconstant and
 has minimal period $\tau/m$ for some $m\in\N$, then
${\cal O}=S_{\tau/m}\cdot \gamma_0$ is a $1$-dimensional
$C^3$-submanifold diffeomorphic to the circle. We define
\begin{equation}\label{e:5.4}
\bigl(m^-_{\tau}({\cal O}), m^0_{\tau}({\cal O})\bigr): =
\bigl(m^-_{\tau}(\gamma_0), m^0_{\tau}(\gamma_0)\bigr).
\end{equation}
if ${\cal O}$ is  a single point critical orbit ${\cal
O}=\{\gamma_0\}$, i.e., $\gamma_0$ is constant, and
\begin{equation}\label{e:5.5}
 \bigl(m^-_{\tau}({\cal O}), m^0_{\tau}({\cal O})\bigr) =
\bigl(m^-_{\tau}(x), m^0_{\tau}(x)-1\bigr)
       \quad \forall x\in{\cal O}
\end{equation}
if $\gamma_0$ is nonconstant, where $m^-_{\tau}(x)$ and
$m^0_{\tau}(x)$ are defined as in (\ref{e:3.17}).

 Let $c={\cal L}_\tau|_{\cal O}$. \textsf{Assume that ${\cal O}$ is isolated and nonconstant}. (The case
 of constant orbits has been included in Section 4.1).  We
 may  take a neighborhood $U$ of ${\cal O}$ such that
 ${\cal K}({\cal L}_\tau)\cap U={\cal
 O}$. By \cite[(4.1)]{Lu} we have critical group
 $C_\ast({\cal L}_\tau, {\cal O}; \K)$
of ${\cal L}_\tau$ at ${\cal O}$.   For every $s\in [0,\tau/m]$ the
tangent space $T_{s\cdot \gamma_0}(S_{\tau}\cdot \gamma_0)$ is
$\R(s\cdot \gamma_0)^\cdot$, and the fiber $N({\cal O})_{s\cdot
\gamma_0}$ at $s\cdot \gamma_0$ of the normal bundle $N({\cal O})$
of ${\cal O}$ is a subspace of codimension $1$ which is orthogonal
to $(s\cdot \gamma_0)^\cdot$ in $T_{s\cdot
\gamma_0}H_{\tau}(\alpha)$, i.e.
\begin{equation}\label{e:5.6}
N({\cal O})_{s\cdot \gamma_0}=\left\{ \xi\in T_{s\cdot
\gamma_0}H_{\tau}(\alpha)\,|\, \langle\!\langle\xi, (s\cdot
\gamma_0)^\cdot\rangle\!\rangle_1=0\,\right\}.
\end{equation}
Since $H_\tau(\alpha)$ is $C^4$-smooth and ${\cal O}$ is a
$C^3$-smooth submanifold, $N({\cal O})$ is $C^2$-smooth manifold.
\footnote{This is the reason that we require higher smoothness of
$M$ and $L$.} Notice that $N({\cal O})$ is invariant under the
$S_{\tau}$-actions in (\ref{e:5.3}) and each $[s]_{\tau}$ gives an
isometric bundle map
\begin{equation}\label{e:5.7}
N({\cal O})\to N({\cal O}),\,(z, v)\mapsto ([s]_{\tau}\cdot z,
[s]_{\tau}\cdot v).
\end{equation}

Note that for $j=1,k$ and sufficiently small $\delta>0$ the set
$$
N(\psi^j({\cal O}))(j\delta):=\bigl\{(y,v)\in N(\psi^j({\cal
O}))\,|\,y\in
  \psi^j({\cal O}),\,\|v\|_1<j\delta\bigr\}
   $$
is contained in an open neighborhood of the zero section of the
tangent bundle $TH_{j\tau}(\alpha)$. By Theorem~1.3.7 on the page 20
of \cite{Kl} we have a $C^2$-embedding  from $N(\psi^j({\cal
O}))(j\delta)$ to an open neighborhood of  the diagonal of
$H_{j\tau}(\alpha)\times H_{j\tau}(\alpha)$,
$$
N(\psi^j({\cal O}))(j\delta)\to H_{j\tau}(\alpha)\times
H_{j\tau}(\alpha),\; (y, v)\mapsto (y, \exp_yv),
$$
where $\exp$ is the exponential map of the chosen Riemannian metric
on $M$ and $(\exp_yv)(t)=\exp_{y(t)}{v(t)}\;\forall t\in\R$. This
yields  a $C^2$-diffeomorphism from $N(\psi^j({\cal O}))(j\delta)$
to an open neighborhood $U_{j\delta}(\psi^j({\cal O}))$ of
$\psi^j({\cal O})$,
\begin{equation}\label{e:5.8}
\Psi_{j\tau}: N(\psi^j({\cal O}))(j\delta)\to
U_{j\delta}(\psi^j({\cal O}))
\end{equation}
given by $\Psi_{j\tau}(y, v)(t)=\exp_{y(t)}{v(t)}\;\forall t\in\R$.
 Clearly,
\begin{equation}\label{e:5.9}
\Psi_{j\tau}(y,0)=y\;\forall y\in\psi^j({\cal O})\quad{\rm and}\quad
\Psi_{j\tau}([s]_{j\tau}\cdot y,\, [s]_{j\tau}\cdot
v)=[s]_{j\tau}\cdot\Psi_\tau(y,  v)
\end{equation}
for any $(y, v)\in N(\psi^j({\cal O}))(j\delta)$ and $[s]_{j\tau}\in
S_{j\tau}$. It follows that $U_{j\delta}(\psi^j({\cal O}))$ is an
$S_{j\tau}$-invariant neighborhood of $\psi^j({\cal O})$, and that
$\Psi_{j\tau}$ is $S_{j\tau}$-equivariant. We also require
$\delta>0$ so small that $U_{j\delta}(\psi^j({\cal O}))$ contains no
other critical orbit besides $\psi^j({\cal O})$, and that
$\Psi_{j\tau}\bigr(\{y\}\times N(\psi^j({\cal O}))_y(j\delta)\bigl)$
and $\psi^j({\cal O})$ have a unique intersection point $y$ (after
identifying $\psi^j({\cal O})$ with the zero section of
$N(\psi^j({\cal O}))$), where
$$
N(\psi^j({\cal O}))_y(j\delta):=N(\psi^j({\cal O}))(j\delta)\cap
N(\psi^j({\cal O}))_y.
$$
Define
\begin{equation}\label{e:5.10}
{\cal F}_{j\tau}: N(\psi^j({\cal O}))(j\delta)\to \R,\;(y,v)\mapsto
{\cal L}_{j\tau}\circ\Psi_{j\tau}(y,v).
\end{equation}
It is $C^{2-0}$, and  satisfies the (PS) condition and
$$
{\cal F}_{j\tau}([s]_{j\tau}\cdot y,\, [s]_{j\tau}\cdot v)={\cal
F}_{j\tau}(y,\,  v)
$$
for any $(y, v)\in N(\psi^j({\cal O}))(j\delta)$ and $[s]_{j\tau}\in
S_{j\tau}$.

Following \cite[Th.2.3]{Wa} we may construct Gromoll-Meyer pairs of
$\psi^j({\cal O})$ as  critical submanifolds of ${\cal F}_{j\tau}$
on $N(\psi^j({\cal O}))(j\delta)$,
\begin{equation}\label{e:5.11}
\bigl(W(\psi^j({\cal O})), W(\psi^j({\cal O}))^-\bigr),\;j=1,k.
\end{equation}
Precisely, set $h_{j\tau}(y,v)=\lambda{\cal F}_{j\tau}(y,v)+
\|v\|_1^2$, and
\begin{equation}\label{e:5.12}
\left.\begin{array}{ll} &W(\psi^j({\cal O}))=({\cal
F}_{j\tau})^{-1}[jc-j\epsilon_1,
jc+ j\epsilon_1]\cap(h_{j\tau})_{j\epsilon_2},\\
&W(\psi^j({\cal O}))^-=({\cal F}_{j\tau})^{-1}(jc-j\epsilon_1)\cap
W(\psi^j({\cal O})).
\end{array}\right\}
\end{equation}
Here positive constants $\lambda$, $\epsilon_1$ and $\epsilon_2$ are
determined by the following conditions.
\begin{eqnarray*}
&&\hbox{$\bullet$ ${\cal F}_{j\tau}$ has a unique
 critical
value $jc$ in $[jc-j\varepsilon, jc+ j\varepsilon]$};\\
 &&\hbox{$\bullet$ $N(\psi^j({\cal
O}))(\frac{j\delta}{2})\subset W(\psi^j({\cal O}))\subset
N(\psi^j({\cal O}))(j\delta)\cap({\cal
F}_{j\tau})^{-1}[jc-j\varepsilon, jc+ j\varepsilon]$};\\
&&\hbox{$\bullet$ $({\cal F}_{j\tau})^{-1}[jc-j\epsilon_1, jc+
j\epsilon_1]\cap(h_{j\tau})_{j\epsilon_2}\subset N(\psi^j({\cal
O}))(j\delta)\setminus N(\psi^j({\cal O}))(\frac{j\delta}{2})$};\\
&&\hbox{$\bullet$ $(dh_{j\tau}(y,v), d{\cal F}_{j\tau}(y,v))>0$ for
any $(y,v)\in N(\psi^j({\cal O}))(j\delta)\setminus N(\psi^j({\cal
O}))(\frac{j\delta}{2})$}.
\end{eqnarray*}
 (Note that different from \cite{Wa}
the present $S_{j\tau}$-action on $N(\psi^j({\cal O}))(j\delta)$ is
only continuous; but the arguments there can still be carried out
due to the special property of our $S_{j\tau}$-action  and the
definition of ${\cal F}_{j\tau}$.) For any $y\in \psi^j({\cal O})$,
the restriction
$$
{\cal F}_{j\tau}|_{N(\psi^j({\cal O}))_y(j\delta)}
$$
has a unique critical point $0=(y, 0)$ in $N(\psi^j({\cal
O}))_y(j\delta)$ (the fiber of disk bundle $N(\psi^j({\cal
O}))(j\delta)$ at $y$), and
\begin{eqnarray}\label{e:5.13}
&&\bigl(W(\psi^j({\cal O}))_y, W(\psi^j({\cal
O}))_y^-\bigr)\\
&&:=\bigl(W(\psi^j({\cal O}))\cap N(\psi^j({\cal O}))_y(j\delta),\,
W(\psi^j({\cal O}))^-\cap
  N(\psi^j({\cal O}))_y(j\delta)\bigr)\nonumber
   \end{eqnarray}
is a  Gromoll-Meyer pair of  ${\cal F}_{j\tau}|_{N(\psi^j({\cal
O}))_y(j\delta)}$ at its isolated critical point $0=(y,0)$
satisfying
\begin{equation}\label{e:5.14}
\left.\begin{array}{ll}
 &\bigl(W(\psi^j({\cal
O}))_{[s]_{j\tau}\cdot y},\, W(\psi^j({\cal O}))^-_{[s]_{j\tau}\cdot
y}\bigr)\\
&= \bigl([s]_{j\tau}\cdot W(\psi^j({\cal O}))_y,\, [s]_{j\tau}\cdot
W(\psi^j({\cal O}))^-_y\bigr)\end{array}\right\}
\end{equation}
for any $[s]_{j\tau}\in S_{j\tau}$ and $y\in \psi^j({\cal O})$
(\cite[Th.2.3]{Wa}). By (\ref{e:5.12}) it is easily checked that
\begin{equation}\label{e:5.15}
\psi^k\bigl(W({\cal O})_y\bigr)\subset W(\psi^k({\cal
O}))_{\psi^k(y)}\quad\hbox{and}\quad \psi^k\bigl(W({\cal
O})_y^-\bigr)\subset W(\psi^k({\cal O}))^-_{\psi^k(y)}
\end{equation}
for each $y\in {\cal O}$. Clearly, for $j=1,k$,
\begin{equation}\label{e:5.16}
\bigl(\widehat W(\psi^j({\cal O})), \widehat W(\psi^j({\cal
O}))^-\bigr):=\bigl(\Psi_{j\tau}(W(\psi^j({\cal O}))),
\Psi_{j\tau}(W(\psi^j({\cal O}))^-)\bigr)
\end{equation}
 are  Gromoll-Meyer pairs of ${\cal L}_{j\tau}$ at $\psi^j({\cal O})$, which is also
$S_{j\tau}$-invariant.

\begin{theorem}\label{th:5.1}
{\rm (\cite[Theorem 4.11]{Lu})}  For an isolated critical
submanifold ${\cal O}={\rm S}_{\tau}\cdot \gamma_0$ of ${\cal
L}_{\tau}$ in $H_{\tau}(\alpha)$, suppose that for some $k\in\N$ the
critical submanifold $\psi^k({\cal O})={\rm S}_{k\tau}\cdot
\gamma_0^k$ of ${\cal L}_{k\tau}$ in $H_{k\tau}(\alpha^k)$ is also
isolated, and that
\begin{equation}\label{e:5.17}
m^-_{k\tau}(\psi^k({\cal O}))= m^-_{\tau}({\cal O})\quad{\rm
and}\quad m^0_{k\tau}(\psi^k({\cal O}))= m^0_{\tau}({\cal O}).
\end{equation}
 Then for $c={\cal L}_\tau|_{\cal O}$ and small
   $\epsilon>0$ there exist
   Gromoll-Meyer pairs  of ${\cal L}_{\tau}$ at ${\cal O}\subset H_\tau(\alpha)$ and
 of ${\cal L}_{k\tau}$ at $\psi^k({\cal O})\subset H_{k\tau}(\alpha^k)$
\begin{eqnarray*}
&&(\widehat W({\cal O}), \widehat W({\cal
  O})^-)\subset\bigl(({\cal L}_\tau)^{-1}[c-\epsilon,
c+ \epsilon], ({\cal L}_\tau)^{-1}(c-\epsilon)\bigr)
 \quad\hbox{and}\\
&&\bigl(\widehat W(\psi^k({\cal O})), \widehat W(\psi^k({\cal
  O}))^-\bigr)\subset\bigl(({\cal
L}_{k\tau})^{-1}[kc-k\epsilon, kc+ k\epsilon], ({\cal
L}_{k\tau})^{-1}(kc-k\epsilon)\bigr),
\end{eqnarray*}
such that
$$
\bigl(\psi^k(\widehat W({\cal O})), \psi^k(\widehat W({\cal
  O})^-)\bigr)\subset \bigl(\widehat W(\psi^k({\cal O})), \widehat W(\psi^k({\cal
  O}))^-\bigr)
$$
and that the iteration map $\psi^k: H_{\tau}(\alpha)\to
H_{k\tau}(\alpha^k)$ induces an isomorphism:
\begin{eqnarray*}
 &&\psi^k_*:  C_*({\cal L}_{\tau}, {\cal O}; \K):=H_\ast\bigl(\widehat W({\cal O}), \widehat W({\cal
  O})^-; \K\bigr)\\
 &&\qquad \longrightarrow
         C_*({\cal L}_{k\tau}, \psi^k({\cal O}); \K):
         =H_\ast\bigl(\widehat W(\psi^k({\cal O})), \widehat W(\psi^k({\cal
  O}))^-; \K\bigr).
\end{eqnarray*}
\end{theorem}

\noindent{\bf Proof}. The following commutative diagram
$$
\begin{CD}
\bigl(W({\cal O}), W({\cal O})^-\bigr) @>\psi^k
>> \bigl(W(\psi^k({\cal O})), W(\psi^k({\cal
O}))^-\bigr)
 \\
@V \Psi_{\tau} VV @VV \Psi_{k\tau} V \\
\bigl(\widehat W({\cal O}), \widehat W({\cal O})^-\bigr)@>\psi^k
>>\bigl(\widehat W(\psi^k({\cal O})), \widehat W(\psi^k({\cal
O}))^-\bigr)
\end{CD}
$$
implies that we only need to prove

\begin{claim}\label{cl:5.2}  $\psi^k$ induces an isomorphism
$$
(\psi^k)_\ast:H_\ast(W({\cal O}), W({\cal O})^-; \K)\to
H_\ast(W(\psi^k({\cal O})), W(\psi^k({\cal O}))^-; \K).
$$
\end{claim}

Recall that for $j=1,2$ the fibers $N(\psi^j({\cal
O}))_{[s]_{j\tau}\cdot\gamma_0^j}$ are subspaces which are
orthogonal to $([s]_{j\tau}\cdot\gamma_0^j)^{.}$ in
$T_{[s]_{j\tau}\cdot\gamma_0^j}H_{j\tau}(\alpha^j)$. We have natural
bundle trivializations
\begin{eqnarray*}
&&\Gamma_j: N(\psi^j({\cal O}))\to S_{j\tau}\cdot\gamma_0^j\times
N(\psi^j({\cal O}))_{\gamma^j_0},\\
&& \hspace{10mm} ([s]_{j\tau}\cdot\gamma_0^j,\, v)\mapsto
([s]_{j\tau}\cdot\gamma_0^j,\, [-s]_{j\tau}\cdot v).\nonumber
\end{eqnarray*}
From (\ref{e:5.14})-(\ref{e:5.15}) we get the commutative diagram
$$
\begin{CD}
\bigl(W({\cal O}), W({\cal O})^-\bigr) @>\Gamma_1
>> \bigl(S_\tau\cdot\gamma_0\times W({\cal O})_{\gamma_0}, S_\tau\cdot\gamma_0\times W({\cal
O})_{\gamma_0}^-\bigr)
 \\
@V \psi^k VV @VV \psi^k V \\
\bigl(W(\psi^k({\cal O})), W(\psi^k({\cal O}))^-\bigr)@>\Gamma_k
>>\bigl(S_{k\tau}\cdot\gamma_0^k\times W(\psi^k({\cal O}))_{\gamma_0^k}, S_{k\tau}\cdot\gamma_0^k\times W(\psi^k({\cal
O}))_{\gamma_0^k}^-\bigr)
\end{CD}
$$
So Claim~\ref{cl:5.2} is equivalent to

\begin{claim}\label{cl:5.3}
  $\psi^k$ induces an isomorphism
\begin{eqnarray*}
&&(\psi^k)_\ast:H_\ast\bigl(S_\tau\cdot\gamma_0\times W({\cal
O})_{\gamma_0}, S_\tau\cdot\gamma_0\times W({\cal O})_{\gamma_0}^-;
\K\bigr)\to\\
&&\hspace{30mm} H_\ast\bigl(S_{k\tau}\cdot\gamma^k\times
W(\psi^k({\cal O}))_{\gamma_0^k}, S_{k\tau}\cdot\gamma^k\times
W(\psi^k({\cal O}))_{\gamma_0^k}^-; \K\bigr).
\end{eqnarray*}
\end{claim}

Since $(\psi^k)_\ast:H_\ast\bigl(S_\tau\cdot\gamma_0; \K\bigr)\to
H_\ast\bigl(S_{k\tau}\cdot\gamma_0^k; \K\bigr)$ is an isomorphism,
and
\begin{eqnarray*}
&&H_q\bigl(S_\tau\cdot\gamma_0\times W({\cal O})_{\gamma_0},
S_\tau\cdot\gamma_0\times W({\cal O})_{\gamma_0}^-;
\K\bigr)\\
&&\hspace{20mm}=\bigoplus^q_{j=0}H_j\bigl(S_\tau\cdot\gamma_0;
\K\bigr)\otimes H_{q-j}\bigl(W({\cal O})_{\gamma_0},
 W({\cal O})_{\gamma_0}^-;
\K\bigr),\\
&&H_q\bigl(S_{k\tau}\cdot\gamma_0^k\times W(\psi^k({\cal
O}))_{\gamma_0^k}, S_{k\tau}\cdot\gamma_0^k\times W(\psi^k({\cal
O}))_{\gamma_0^k}^-; \K\bigr)\\
&&\hspace{20mm}=\bigoplus^q_{j=0}H_j\bigl(S_{k\tau}\cdot\gamma_0^k;
\K)\otimes H_{q-j}\bigl(W(\psi^k({\cal O}))_{\gamma_0^k},
 W(\psi^k({\cal O}))_{\gamma_0^k}^-;\K\bigr)
\end{eqnarray*}
by the Kunneth formula, Claim~\ref{cl:5.3} is equivalent to

\begin{claim}\label{cl:5.4}
 $\psi^k$ induces an isomorphism
\begin{eqnarray*}
(\psi^k)_\ast:H_\ast\bigl(W({\cal O})_{\gamma_0},
 W({\cal O})_{\gamma_0}^-;
\K\bigr)\to H_\ast\bigl(W(\psi^k({\cal O}))_{\gamma_0^k},
 W(\psi^k({\cal O}))_{\gamma_0^k}^-;\K\bigr).
\end{eqnarray*}
\end{claim}

For conveniences we write
\begin{equation}\label{e:5.18}
{\cal F}^N_{j\tau}:={\cal F}_{j\tau}|_{N(\psi^j({\cal
O}))_{\gamma_0^j}(j\delta)}.
\end{equation}
Since $\bigl(W(\psi^j({\cal O}))_{\gamma_0^j}, W(\psi^j({\cal
O}))_{\gamma_0^j}^-\bigr)$ are  Gromoll-Meyer pairs of  ${\cal
F}^N_{j\tau}$ at isolated critical points $0=(\gamma_0^j,0)$ with
respect to the flow of negative gradients, $j=1, k$, and
\begin{eqnarray}
\psi^k\Bigl(N({\cal O})_{\gamma_0}\Bigr)\subset N(\psi^k({\cal
O}))_{\gamma^k_0},\label{e:5.19}
\end{eqnarray}
by (\ref{e:5.6}) and the relations $\langle\!\langle\psi^k(\xi),
\psi^k(\eta)\rangle\!\rangle_1=k \langle\!\langle\xi,
\eta\rangle\!\rangle_1\;\forall\xi,\eta\in
T_{\gamma_0}H_{\tau}(\alpha)$, the proofs from Claim~\ref{cl:2.2} to
Claim~\ref{cl:2.3} show that Claim~\ref{cl:5.4} is equivalent to

\begin{claim}\label{cl:5.5}
 There exist small open neighborhoods
$V^{(j)}$ of $0\in N(\psi^j({\cal O}))_{\gamma_0^j}(j\delta)$
 with $\psi^k(V^{(1)})\subset V^{(k)}$, such that $\psi^k$ induces isomorphisms {\footnotesize
$$
(\psi^k)_\ast: H_\ast\bigl(({\cal F}^N_{\tau})_{c}\cap V^{(1)},
(({\cal F}^N_{\tau})_{c}\setminus\{0\})\cap V^{(1)};\K\bigr)\to
H_\ast\bigl({\cal F}^N_{k\tau})_{kc}\cap V^{(k)}, (({\cal
F}^N_{k\tau})_{kc}\setminus\{0\})\cap V^{(k)};\K\bigr).
$$}
\end{claim}

For $j=1,k$ let $(\Psi_{j\tau})_{\gamma^j_0}$ be the restrictions of
the maps $\Psi_{j\tau}$ in (\ref{e:5.8}) to the fibers
$N(\psi^j({\cal O}))(j\delta)_{\gamma^j_0}$. By shrinking $\delta>0$
we may assume that their images are contained in the images of the
charts $\phi_{j\tau}$ in \cite[(3.8)]{Lu}. Let
\begin{equation}\label{e:5.20}
\Upsilon_{\gamma^j_0}:T_{\gamma_0^j}H_{j\tau}(\alpha^j)=W^{1,2}((\gamma_0^j)^\ast
TM)\to \tilde H_{j\tau}
\end{equation}
be the inverses of the tangent maps $d\phi_{j\tau}(0): \tilde
H_{j\tau}=T_0\tilde H_{j\tau}\to T_{\gamma_0^j}H_{j\tau}(\alpha^j)$
 given by
$$
d\phi_{j\tau}(0)(\tilde\alpha)(t)=\frac{d}{ds}\Bigm|_{s=0}\phi_{j\tau}(0)(s\tilde\alpha)(t)=\Phi(t)\tilde\alpha(t)\quad
\forall t.
$$
Clearly,  one has
\begin{equation}\label{e:5.21}
\psi^k\circ\Upsilon_{\gamma_0}=\Upsilon_{\gamma_0^k}\circ\psi^k.
\end{equation}
 By (\ref{e:5.8}) it is easily seen that the
compositions
$$
(\phi_{j\tau})^{-1}\circ(\Psi_{j\tau})_{\gamma^j_0}=\Upsilon_{\gamma^j_0}\bigm|_{
N(\psi^j({\cal O}))(j\delta)_{\gamma^j_0}}.
$$
 For $j=1,k$ let us define
\begin{eqnarray*}
&&XN(\psi^j({\cal O}))_{\gamma^j_0}:=N(\psi^j({\cal
O}))_{\gamma^j_0}\cap C^1((\gamma_0^j)^\ast TM),\\
&&XN(\psi^j({\cal O}))_{\gamma^j_0}(j\delta):=N(\psi^j({\cal
O}))_{\gamma^j_0}(j\delta)\cap C^1((\gamma_0^j)^\ast TM),\\
&&{\cal N}_{j\tau}:=\Psi_{j\tau}\bigl(N(\psi^j({\cal
O}))_{\gamma^j_0}(j\delta)\bigr),\\
&&{\cal XN}_{j\tau}:=\Psi_{j\tau}\bigl(XN(\psi^j({\cal
O}))_{\gamma^j_0}(j\delta)\bigr),\\
&& {\cal F}^{NX}_{j\tau}:={\cal F}^N_{j\tau}|_{XN(\psi^j({\cal
O}))_{\gamma_0^j}(j\delta)}={\cal F}_{j\tau}|_{XN(\psi^j({\cal
O}))_{\gamma_0^j}(j\delta)},\\
&&{\cal L}_{j\tau}^N:={\cal L}_{i\tau}\bigl|_{{\cal
N}_{j\tau}}\quad\hbox{and}\quad {\cal L}_{j\tau}^{NX}:={\cal
L}_{i\tau}\bigl|_{{\cal XN}_{j\tau}}.
\end{eqnarray*}
Then Banach manifolds $XN(\psi^j({\cal O}))_{\gamma^j_0}$,
$XN(\psi^j({\cal O}))_{\gamma^j_0}(j\delta)$ and ${\cal XN}_{j\tau}$
are dense in Hilbert manifolds in $N(\psi^j({\cal
O}))_{\gamma^j_0}$,  $N(\psi^j({\cal O}))_{\gamma^j_0}(j\delta)$ and
${\cal N}_{j\tau}$ respectively. Moreover, $\Psi_{j\tau}$ restrict
to $C^2$ diffeomorphisms from $N(\psi^j({\cal
O}))_{\gamma^j_0}(j\delta)$ (resp. $XN(\psi^j({\cal
O}))_{\gamma^j_0}(j\delta)$) onto ${\cal N}_{j\tau}$ (resp. ${\cal
XN}_{j\tau}$). ${\cal L}_{j\tau}^N$ and ${\cal L}_{j\tau}^{NX}$ are
$C^{2-0}$ and $C^2$, respectively, and have isolated critical points
$\gamma_0^j$, $j=1,k$. Define
\begin{eqnarray*}
&&\tilde S_{j\tau}:=\Upsilon_{\gamma^j_0}\bigl(N(\psi^j({\cal O
}))_{\gamma^j_0}\bigr),\\
&& \tilde
S_{j\tau}(j\delta):=\Upsilon_{\gamma^j_0}\bigl(N(\psi^j({\cal O
}))(j\delta)_{\gamma^j_0}\bigr),\\
&&X\tilde S_{j\tau}:=\Upsilon_{\gamma^j_0}\bigl(XN(\psi^j({\cal O
}))_{\gamma^j_0}\bigr),\\
&& X\tilde
S_{j\tau}(j\delta):=\Upsilon_{\gamma^j_0}\bigl(XN(\psi^j({\cal O
}))(j\delta)_{\gamma^j_0}\bigr),\\
&&\tilde{\cal L}^S_{j\tau}:=\tilde{\cal L}_{j\tau}\bigl|_{\tilde
S_{j\tau}}\quad\hbox{and}\quad \tilde{\cal
L}^{SX}_{j\tau}:=\tilde{\cal L}_{j\tau}\bigl|_{X\tilde S_{j\tau}}.
\end{eqnarray*}
Then $\tilde S_{j\tau}$ is a Hilbert subspace of $\tilde H_{j\tau}$
of codimension $1$, and $\tilde S_{j\tau}(j\delta)$  is an open
neighborhood of the origin of $\tilde S_{j\tau}$. $X\tilde
S_{j\tau}=\tilde S_{j\tau}\cap \tilde X_{j\tau}$ is a Banach
subspace of $\tilde X_{j\tau}$ of codimension $1$, and $X\tilde
S_{j\tau}(j\delta)$ is an open neighborhood of the origin of
$X\tilde S_{j\tau}$. Moreover, $\tilde{\cal L}_{j\tau}^S$ and
$\tilde{\cal L}_{j\tau}^{SX}$ are $C^{2-0}$ and $C^2$, respectively,
and have isolated critical points $0$, $j=1,k$.

By \cite[(3.12)-(3.13)]{Lu} and (\ref{e:5.8}),
$$
\phi_{j\tau}(\tilde\alpha)(t)=\exp_{\gamma^j_0(t)}(\Phi(t)\tilde\alpha(t))=\Psi_{j\tau}(\gamma_0,
\Phi\tilde\alpha)(t)=(\Psi_{j\tau})_{\gamma^j_0}(\Phi\tilde\alpha)(t)\;\forall
t.
$$
From this, (\ref{e:5.10}) and (\ref{e:5.18}) and \cite[(3.16)]{Lu}
it follows that
\begin{equation}\label{e:5.22}
\left.\begin{array}{ll}
 &{\cal F}^N_{j\tau}(\Phi\tilde\alpha)={\cal
L}_{j\tau}\circ\Psi_{j\tau}({\gamma^j_0}, \Phi\tilde\alpha))={\cal
L}_{j\tau}(\phi_{j\tau}(\tilde\alpha))=\tilde{\cal
L}_{j\tau}(\tilde\alpha)\\
& \hbox{or}\;  {\cal F}_{j\tau}^N(v)=\tilde{\cal
L}_{j\tau}(\Upsilon_{\gamma^j_0}v)
\end{array}\right\}
\end{equation}
for all $\tilde\alpha\in\tilde S_{j\tau}(j\delta)$  or $v\in
N(\psi^j({\cal O}))(j\delta)_{\gamma^j_0}$. Moreover, (\ref{e:5.21})
and (\ref{e:5.19}) imply the commutative diagram
$$
\begin{CD}
N({\cal O})(\delta)_{\gamma_0} @>\Upsilon_{\gamma_0}
>> \tilde S_\tau(\delta)
 \\
@V \psi^k VV @VV \psi^k V \\
N(\psi^k({\cal O}))(k\delta)_{\gamma_0^k}@>\Upsilon_{\gamma_0^k}
>>\tilde S_{k\tau}(k\delta)
\end{CD}
$$
and thus
\begin{equation}\label{e:5.23}
\psi^k(\tilde S_\tau)\subset\tilde S_{k\tau}.
\end{equation}
These show  that Claim~\ref{cl:5.5} is equivalent to

\begin{claim}\label{cl:5.6}
 There exist small open neighborhoods
$V^{(j)}$ of $0\in \tilde S_{j\tau}$
 with $\psi^k(V^{(1)})\subset V^{(k)}$, such that $\psi^k$ induces isomorphisms {\footnotesize
$$
(\psi^k)_\ast: H_\ast\bigl((\tilde{\cal L}^S_{\tau})_{c}\cap
V^{(1)}, ((\tilde{\cal L}^S_{\tau})_{c}\setminus\{0\})\cap
V^{(1)};\K\bigr)\to H_\ast\bigl(\tilde{\cal L}^S_{k\tau})_{kc}\cap
V^{(k)}, ((\tilde{\cal L}^S_{k\tau})_{kc}\setminus\{0\})\cap
V^{(k)};\K\bigr).
$$}
\end{claim}

We shall prove it with Corollary~\ref{cor:2.8}. By the definitions
below (\ref{e:5.21}) we have
\begin{equation}\label{e:5.24}
T_{\gamma_0^j}{\cal N}_{j\tau}=N(\psi^j({\cal
O}))_{\gamma^j_0}\quad\hbox{and}\quad T_{\gamma_0^j}{\cal
XN}_{j\tau}=XN(\psi^j({\cal O}))_{\gamma^j_0}
\end{equation}
 and
therefore decompositions
\begin{equation}\label{e:5.25}
T_{\gamma_0^j}H_{j\tau}(\alpha^j)=T_{\gamma_0^j}{\cal
N}_{j\tau}\oplus\R(\gamma_0^j)^\cdot\quad\hbox{and}\quad
T_{\gamma_0^j}XH_{j\tau}(\alpha^j)=T_{\gamma_0^j}{\cal
XN}_{j\tau}\dot{+}\R(\gamma_0^j)^\cdot.
\end{equation}
Here $XH_{j\tau}(\alpha^j)=X_{j\alpha}\cap H_{j\tau}(\alpha^j)$. Let
${\bf H}^-(b)$, ${\bf H}^0(b)$ and ${\bf H}^+(b)$  denote the
negative, null, and negative space of a continuous symmetric
bilinear form $b$ on a Hilbert space, respectively. Since
$\tilde{\cal L}^X_{j\tau}={\cal L}^X_{j\tau}\circ\phi_{j\tau}$
implies
$$
(B_{j\tau}(0)\xi,\eta)_{W^{1,2}}=d^2\tilde{\cal
L}^X_{j\tau}(0)(\xi,\eta)=d^2{\cal
L}^X_{j\tau}(\gamma_0^j)(d\phi_{j\tau}(0)\xi, d\phi_{j\tau}(0)\eta)
$$
for any $\xi,\eta\in T_{\gamma_0^j}XH_{j\tau}(\alpha^j)$,
$d^2\tilde{\cal L}^X_{j\tau}(0)$ and $d^2{\cal
L}^X_{j\tau}(\gamma_0^j)$ may be extended into continuous symmetric
bilinear forms on Hilbert spaces
$T_{\gamma_0^j}XH_{j\tau}(\alpha^j)$ and $\tilde
H_{j\tau}=W^{1,2}(S_{j\tau},\R^n)$, respectively. We also use
$d^2\tilde{\cal L}^X_{j\tau}(0)$ and $d^2{\cal
L}^X_{j\tau}(\gamma_0^j)$ to denote these extensions without
occurring of confusions below. Then we  infer that
\begin{equation}\label{e:5.26}
{\bf H}^\ast(d^2{\cal
L}^X_{j\tau}(\gamma_0^j))=d\phi_{j\tau}(0)\bigl({\bf
H}^\ast(d^2\tilde{\cal L}^X_{j\tau}(0))\bigr)\subset
C^1((\gamma_0^j)^\ast TM)
\end{equation}
and have dimensions $m^\ast({\cal L}_{j\tau}, \gamma_0^j)$ for
$\ast=0,-$ and $j=1,k$. Note that the second differentials $d^2{\cal
L}_{j\tau}^{NX}(\gamma_0^j)$ are the restrictions of $d^2{\cal
L}_{j\tau}^{X}$ to $T_{\gamma_0^j}{\cal XN}_{j\tau}$ as symmetric
bilinear forms, $j=1,k$. It is clear that
$\R(\gamma_0^j)^\cdot\subset {\bf H}^0(d^2{\cal
L}^X_{j\tau}(\gamma_0^j))$ and
$$
{\bf H}^\ast(d^2{\cal L}^{NX}_{j\tau}(\gamma_0^j))\subset{\bf
H}^\ast(d^2{\cal L}^X_{j\tau}(\gamma_0^j)),\quad\ast=0,-.
$$
(Here we understand $d^2{\cal L}^{NX}_{j\tau}(\gamma_0^j)$ in ${\bf
H}^\ast(d^2{\cal L}^{NX}_{j\tau}(\gamma_0^j))$ as the extension of
it on the corresponding Hilbert space).  The finiteness of
dimensions of these spaces and (\ref{e:5.25}) imply
\begin{equation}\label{e:5.27}
\left.\begin{array}{ll}
 & {\bf H}^0(d^2{\cal
L}^X_{j\tau}(\gamma_0^j))={\bf H}^0(d^2{\cal
L}^{NX}_{j\tau}(\gamma_0^j))\oplus\R(\gamma_0^j)^\cdot,\\
&{\bf H}^-(d^2{\cal L}^{NX}_{j\tau}(\gamma_0^j))={\bf H}^-(d^2{\cal
L}^X_{j\tau}(\gamma_0^j)),\;j=1,k.
\end{array}\right\}
\end{equation}
Since  the restriction of $\Psi_{j\tau}$ to $XN(\psi^j({\cal
O}))_{\gamma^j_0}$ is a $C^2$ diffeomorphism onto ${\cal
XN}_{j\tau}$, whose differential at $0=(\gamma_0^j,0)$ is the
identity on $XN(\psi^j({\cal O}))_{\gamma^j_0}$ it follows from
(\ref{e:5.24}) that $d^2{\cal L}^{NX}_{j\tau}(\gamma_0^j))$ and
$d^2{\cal F}^{NX}_{j\tau}(0))$ (both defined on $XN(\psi^j({\cal
O}))_{\gamma^j_0}$) are same and hence have the same extensions. The
final claim leads to
\begin{equation}\label{e:5.28}
{\bf H}^\ast(d^2{\cal L}^{NX}_{j\tau}(\gamma_0^j))={\bf
H}^\ast(d^2{\cal F}^{NX}_{j\tau}(0)),\quad\ast=0,-.
\end{equation}
By (\ref{e:5.22}) we have also
$${\cal F}^{NX}_{j\tau}=\tilde{\cal
L}^{XS}_{j\tau}\circ\bigl(\Upsilon_{\gamma_0^j}\bigm|_{XN(\psi^j({\cal
O}))_{\gamma^j_0}}\bigr)
$$
and so
\begin{equation}\label{e:5.29}
\Upsilon_{\gamma_0^j}\bigl({\bf H}^\ast(d^2{\cal
F}^{NX}_{j\tau}(0))\bigr)={\bf H}^\ast(d^2\tilde{\cal
L}^{XS}_{j\tau}(0)),\quad\ast=0,-
\end{equation}
because $\Upsilon_{\gamma_0^j}\bigm|_{XN(\psi^j({\cal
O}))_{\gamma^j_0}}$ are Banach space isomorphisms onto $X\tilde
S_{j\tau}$, $j=1,k$. From (\ref{e:5.26})-(\ref{e:5.29}) we get
\begin{equation}\label{e:5.30}
\left.\begin{array}{ll} &\dim {\bf H}^0(d^2\tilde{\cal
L}^{X}_{j\tau}(0))=\dim {\bf H}^0(d^2\tilde{\cal
L}^{XS}_{j\tau}(0))+1,\vspace{1mm} \\
&\dim {\bf H}^-(d^2\tilde{\cal L}^{X}_{j\tau}(0))=\dim {\bf
H}^-(d^2\tilde{\cal L}^{XS}_{j\tau}(0)).
\end{array}\right\}
\end{equation}
(Here ${\bf H}^\ast(d^2\tilde{\cal L}^{XS}_{j\tau}(0))$ is
understand as before). These imply that
\begin{equation}\label{e:5.31}
{\bf H}^-(d^2\tilde{\cal L}^{X}_{j\tau}(0))= {\bf
H}^-(d^2\tilde{\cal L}^{XS}_{j\tau}(0)),\;j=1,k
\end{equation}
and that ${\bf H}^0(d^2\tilde{\cal L}^{XS}_{j\tau}(0)$ have
codimension one in ${\bf H}^0(d^2\tilde{\cal L}^{X}_{j\tau}(0))$,
$j=1,k$. Let $e$ be a nonzero element in the orthogonal
complementary of ${\bf H}^0(d^2\tilde{\cal L}^{XS}_{\tau}(0)$  in
${\bf H}^0(d^2\tilde{\cal L}^{X}_{\tau}(0))$ with respect to the
inner product of $\tilde H_{\tau}$. Then
$$
{\bf H}^0(d^2\tilde{\cal
L}^{X}_{j\tau}(0))= {\bf H}^0(d^2\tilde{\cal
L}^{XS}_{j\tau}(0))\oplus\R\psi^j(e),\;j=1,k
$$
Recall the orthogonal decompositions $\tilde
H_{j\tau}=M^0(\tilde\gamma_0^j)\oplus M(\tilde\gamma_0^j)^-\oplus
M(\tilde\gamma_0^j)^+$ according to the null, negative, and positive
definiteness of the operator $B_{j\tau}(0)$, where
$B_{j\tau}(0)=B_{j\tau}(\tilde\gamma_0^j)$ is given as above
(\ref{e:3.7}), $j=1,k$. Because ${\bf H}^\ast(d^2\tilde{\cal
L}^{X}_{j\tau}(0))=M^\ast(\tilde\gamma_0^j)$, $\ast=+, 0, -$, we may
obtain orthogonal decompositions
\begin{equation}\label{e:5.32}
\tilde H_{j\tau}= \tilde S_{j\tau}\oplus\R\psi^j(e),\quad j=1,k.
\end{equation}
 Let $P^{(j)}$ be
the orthogonal projections from $\tilde H_{j\tau}$ onto $\tilde
S_{j\tau}$ in the decompositions. Since $\psi^j(e)\in\tilde
X_{j\tau}$, using the Banach inverse operator theorem one easily
prove that $P^{(j)}|_{\tilde X_{j\tau}}$ are continuous linear
operators from $\tilde X_{j\tau}$ onto $X\tilde S_{j\tau}$. It
follows that the map
\begin{equation}\label{e:5.33}
  A^S_{j\tau}:
\tilde S_{j\tau}(j\delta)\cap\tilde X_{j\tau}\to \tilde
X_{j\tau},\quad x\mapsto P^{(j)} A_{j\tau}(x)
\end{equation}
is $C^1$ and that the map
\begin{equation}\label{e:5.34}
 B^S_{j\tau}: \tilde S_{j\tau}(j\delta)\cap\tilde X_{j\tau}\to L_s\bigl(\tilde S_{j\tau}, \tilde S_{j\tau})
\end{equation}
given by $\tilde B^S_{j\tau}(x)=P^{(j)} B_{j\tau}(x)|_{\tilde
S_{j\tau}}$ is  continuous. Here as before the topology on $\tilde
S_{j\tau}(j\delta)\cap\tilde X_{j\tau}$ is one induced by $\tilde
X_{j\tau}$. It is not hard to check that the tuples
$$
\bigl(\tilde S_{j\tau}, X\tilde S_{j\tau},  \tilde{\cal
L}^S_{j\tau},  A^S_{j\tau},  B^S_{j\tau}\bigr)
$$
satisfy the assumptions in Theorem~\ref{th:1.1}, and specially
$$
d^2 \tilde{\cal L}^{XS}_{j\tau}
  (0)(\xi,\eta)=\bigl(B^S_{j\tau}(0)\xi,
  \eta\bigr)_{W^{1,2}}
$$
for any $\xi, \eta\in \tilde S_{k\tau}$. These, (\ref{e:5.17}) and
(\ref{e:5.30}) lead to
$$
m^0(\tilde{\cal L}^S_{\tau}, 0)=m^0(\tilde{\cal L}^S_{k\tau},
0)\quad\hbox{and}\quad m^-(\tilde{\cal L}^S_{\tau},
0)=m^-(\tilde{\cal L}^S_{k\tau}, 0).
$$
Clearly,   (\ref{e:5.23}) implies that $\psi^k(X\tilde
S_{\tau})\subset X\tilde S_{k\tau}$. From these and (4.13) and
(\ref{e:5.32})-(\ref{e:5.35}) it follows that
$$
 \psi^k\bigl( A^S_{k\tau}(x)\bigr)=
 A^S_{k\tau}\bigl(\psi^k(x)\bigr)\quad{\rm
and}\quad\psi^k\bigl( B^S_{k\tau}(x)\xi\bigr)=
B^S_{k\tau}\bigl(\psi^k(x)\bigr)\psi^k(\xi)
$$
for any $x\in \tilde S_\tau(\delta)\cap \tilde X_{\tau}$ and $\xi\in
\tilde S_\tau$. Now Claim~\ref{cl:5.6} follows from
Corollary~\ref{cor:2.8}. Theorem~\ref{th:5.1} is proved.
$\Box$\vspace{2mm}

Define $C_\ast({\cal F}^N_\tau, 0; \K):=H_\ast(W({\cal
O})_{\gamma_0}, W({\cal O})^-_{\gamma_0}; \K)$ and
\begin{eqnarray*}
&&C_\ast({\cal L}_\tau, {\cal O}; \K):=H_\ast\bigl(\widehat W({\cal
O})), \widehat W({\cal O})^-; \K\bigr),\\
&&C_\ast({\cal F}_\tau, {\cal O}; \K):=H_\ast(W({\cal O}), W({\cal
O})^-; \K)
\end{eqnarray*}
via the relative singular homology. By (\ref{e:5.10}) and
(\ref{e:5.16}), $\Psi_\tau$ induces obvious isomorphisms
\begin{equation}\label{e:5.35}
(\Psi_\tau)_\ast: C_\ast({\cal L}_\tau, {\cal O}; \K)
 \cong C_\ast({\cal F}_\tau, {\cal O};\K).
\end{equation}
The bundle trivializations under Claim~\ref{cl:5.2} and
\cite[(2.13), (2.14)]{Wa} lead to
\begin{eqnarray}\label{e:5.36}
  C_q({\cal F}_{\tau}, {\cal O}; \K)
&\cong& \oplus^q_{j=0}\left[C_{q-j}
  \left({\cal F}^N_{\tau}, 0; \K\right) \otimes
     H_j(S_{\tau}; \K)\right]\nonumber  \\
&\cong& C_{q-1}
  \left({\cal F}^N_{\tau}, 0; \K\right)\nonumber \\
&\cong& C_{q-1}
  \left(\tilde{\cal L}^S_{\tau}, 0; \K\right)
 \end{eqnarray}
for any $q\in \{0\}\cup\N$, where the third ``$\cong$" is due to
$$
{\cal F}^{N}_{\tau}=\tilde{\cal
L}^{S}_{\tau}\circ\bigl(\Upsilon_{\gamma_0}\bigm|_{N({\cal
O})_{\gamma_0}}\bigr)
$$
by (\ref{e:5.22}).  Recall that $m^-(\tilde{\cal L}^S_{\tau},
0)=m^-(\tilde{\cal L}_{\tau}, 0)=m^-({\cal O})$ by (\ref{e:5.30})
and (\ref{e:5.5}). Applying Corollary~\ref{cor:1.2} to $\tilde{\cal
L}^S_{\tau}$, and $\tilde{\cal L}_{\tau}$ (if ${\cal O}$ is a
constant orbit) we obtain Lemma~4.12 in \cite{Lu}, i.e.

\begin{lemma}\label{lem:5.7}
Suppose that $C_q({\cal L}_{\tau}, {\cal O};\K)\ne 0$ for ${\cal
O}=S_{\tau}\cdot \gamma$. Then
$$
q-2n\le q-1- m^0_{\tau}({\cal O})\le m^-_{\tau}({\cal O})\le q-1
$$
if ${\cal O}$ is not a single point critical orbit, i.e. $\gamma$ is
not constant,  and
$$
q-2n\le q- m^0_{\tau}({\cal O})\le m^-_{\tau}({\cal O})\le q
$$
otherwise.
\end{lemma}

Here is Lemma~4.13 in \cite{Lu}.

\begin{lemma}\label{lem:5.8}
Suppose that $C_q({\cal L}_{\tau}, {\cal O};\K)\ne 0$ for ${\cal
O}=S_{\tau}\cdot \gamma$. If either ${\cal O}$ is not a single point
critical orbit and $q>1$, or ${\cal O}$ is a single point critical
orbit and $q>0$, then each point in ${\cal O}$ is non-minimal saddle
point.
\end{lemma}

\noindent{\bf Proof}.\quad  When ${\cal O}$ is a single point
critical orbit and $q>0$, the conclusion follows from
Corollary~\ref{cor:1.2}. Now assume that ${\cal
O}=S_\tau\cdot\gamma$ is not a single point critical orbit and
$q>1$. By (\ref{e:5.35}) and (\ref{e:5.36})  we have
\begin{eqnarray*}
  0\ne C_q({\cal L}_{\tau}, {\cal O}; \K)\cong C_{q-1}
  \left(\tilde{\cal L}^S_{\tau}, 0; \K\right).
 \end{eqnarray*}
By Corollary~\ref{cor:1.2}, $\gamma$ and hence every point of ${\cal
O}$ is a non-minimal saddle point of ${\cal L}_{\tau}$.
$\Box$\vspace{2mm}

\begin{remark}\label{rm:5.9}
{\rm Let us outline how our method above can be used to give the
shifting theorem of critical groups of the energy functional of a
Finsler metric on a compact manifold at a nonconstant critical
orbit. For a regular Finsler metric $F$ on $M$, by \cite{Me} the
energy functional
$$
{\cal L}:
H_1,\;\gamma\mapsto\int^1_0L(\gamma(t),\dot\gamma(t))dt=\int^1_0[F(\gamma(t),\dot\gamma(t))]^2dt,
$$
is $C^{2-0}$, and satisfies the (PS) condition. For a nonconstant
critical orbit ${\cal O}=S_1\cdot\gamma_0$ there exists a constant
$c>0$ such that $F(\gamma(t),\dot\gamma(t)\equiv c>0$ for any
$\gamma\in{\cal O}$. Note that $\gamma_0$ is at least $C^2$. Suppose
that $\gamma_0^\ast TM\to S_1$ is trivial. As usual (cf. Section~3
and \cite[\S2]{Lu}) we may assign its Maslov-type index
$i_1(\gamma_0)$ and $\nu_1(\gamma_0)$.

By (\ref{e:5.35}) and (\ref{e:5.36}) we have
\begin{equation}\label{e:5.37}
 C_q({\cal L}, {\cal O}; \K)\cong C_{q-1}
  ({\cal F}^N, 0; \K)\quad\forall q\in\N\cup\{0\}.
\end{equation}
Now we modify $L$ near the zero section to $\hat L$ so that
$L(x,v)=\hat L(x,v)$ if $F(x,v)>\frac{1}{2}c$ and that $\hat L$
satisfying the conditions (L1)-(L3) in \cite{Lu}. We also choose it
to ensure that one may use the stability theorem of critical groups
(cf. \cite[Th.5.6]{Ch} and \cite[Th.3.6]{CiDe}) to prove
\begin{equation}\label{e:5.38}
C_{\ast}({\cal F}^N, 0; \K)\cong C_{\ast}
  (\hat{\cal F}^N, 0; \K).
\end{equation}
Then we can apply Corollary~\ref{cor:1.2} to get a $\delta>0$ and a
(unique) $C^1$-map
$$
h:{\bf H}^0(d^2{\cal L}^{NX}(\gamma_0))\cap {\bf B}_\delta(XN({\cal
O})_{\gamma_0})\to {\bf H}^-(d^2{\cal L}^{NX}(\gamma_0))+ {\bf
H}^+(d^2{\cal L}^{NX}(\gamma_0))\cap XN({\cal O})_{\gamma_0}
$$
with $h(0)=0$, such that
\begin{equation}\label{e:5.39}
 C_{j}(\hat{\cal F}^N, 0; \K)=C_{j-i_1({\cal O})}(\hat{\cal F}^{N\circ}, 0;
 \K)\;\forall j,
  \end{equation}
 where $i_1({\cal O})=i_1(\gamma_0)=\dim {\bf H}^-(d^2{\cal L}^{NX}(\gamma_0))$ and
$\nu_1(\gamma_0)-1=\dim{\bf H}^0(d^2{\cal L}^{NX}(\gamma_0))$,
 ${\bf B}_\delta(XN({\cal
O})_{\gamma_0})$ is a ball of radius $\delta$ and centrad at $0$ in
$XN({\cal O})_{\gamma_0}$, and
$$
\hat{\cal F}^{N\circ}:{\bf H}^0(d^2{\cal L}^{NX}(\gamma_0))\cap {\bf
B}_\delta(XN({\cal O})_{\gamma_0})\to\R
$$
is given by $\hat{\cal F}^{N\circ}(\xi)=\hat{\cal F}^{N}(\xi+
h(\xi))=\hat{\cal L}\bigl( \exp_{\gamma_0}(\xi+ h(\xi))\bigr)$. Note
that $\xi+ h(\xi)\in C^1((\gamma_0)^\ast TM)$ sits in a small
neighborhood of the zero section if $\delta>0$ is small enough.
Hence we may require
$$
F(\gamma(t), \dot\gamma(t))>\frac{3}{4}c\quad\forall t
$$
for any $\gamma=\exp_{\gamma_0}(\xi+ h(\xi))$ with $\xi\in {\bf
H}^0(d^2{\cal L}^{NX}(\gamma_0))\cap {\bf B}_\delta(XN({\cal
O})_{\gamma_0})$. This implies
$$
{\cal F}^{N\circ}(\xi):={\cal F}^{N}(\xi+ h(\xi))={\cal L}\bigl(
\exp_{\gamma_0}(\xi+ h(\xi))\bigr)=\hat {\cal L}\bigl(
\exp_{\gamma_0}(\xi+ h(\xi))\bigr)
$$
for any $\xi\in {\bf H}^0(d^2{\cal L}^{NX}(\gamma_0))\cap {\bf
B}_\delta(XN({\cal O})_{\gamma_0})$. From these and
(\ref{e:5.37})-(\ref{e:5.39}) we get the following \textsf{shifting
theorem}
\begin{equation}\label{e:5.40}
 C_q({\cal L}, {\cal O}; \K)\cong C_{q-1-i_1({\cal O})}
  ({\cal F}^{N\circ}, 0; \K)\quad\forall q\in\N\cup\{0\}.
\end{equation}
The detailed proof (including the case of Finsler-geodesics on $M$
joining orthogonally two submanifolds $M_1$ and $M_2$ of $M$) will
be given in other place. (It seems that this result was proved by
the finite dimensional approximations of the loop space by broken
geodesics in Finsler geometry). }
\end{remark}

\section{The corrections of  Sections~5, 6,7 in \cite{Lu}}\label{sec:6}

Since we  only make corrections for the proofs of Theorems~4.4, 4.7
in \cite{Lu}, the proofs in Sections~5,6 of \cite{Lu} are correct
except that ``the generalized Morse lemma'' in line 10 and `` the
shifting theorem ([14] and [7,p.50])'' in lines 4-5 from bottom
should be changed into ``Theorem~1.1'' and ``Corollary~1.2'' in this
paper, respectively.

For Section~7 in \cite{Lu} we need to make a few of replacements as
follows:

\noindent{``Lemma~4.12'' in line 5 on the page 3021 of \cite{Lu},\\
``(4.67)'' in line 12 on the page 3021 of \cite{Lu},\\
``Lemma~4.12'' in line 2 from bottom on the page 3021 of \cite{Lu},\\
``(4.53) '' in line 1 on the page 3022 of \cite{Lu},\\
``Theorem~4.11'' in line 5 from bottom on the page 3022 of \cite{Lu},\\
``Lemma~4.13'' in line 11 on the page 3024 of \cite{Lu},\\
\textsf{are respectively changed into}:  `` Lemma~\ref{lem:5.7}'',
``(\ref{e:5.35}) and (\ref{e:5.36})'', `` Lemma~\ref{lem:5.7}'',
``(\ref{e:5.5})'', ``Theorem~5.1'' and `` Lemma~\ref{lem:5.8}'' in
this paper.

\begin{remark}\label{rm:6.1}
{\rm For a Tonelli Lagrangian $L\in C^2(S_\tau\times TM,\R)$  and a
$\tau$-periodic solution $\gamma$ of the corresponding Lagrangian
system (\ref{e:3.2}), assume that $\gamma^\ast TM\to S_\tau$ is
trivial one can still assign two sequences of integers (Maslov-type
index) $\{i_{k\tau}(\gamma^k)\,|\,k\in\N\}$ and
$\{\nu_{k\tau}(\gamma^k)\,|\, k\in\N\}$ (cf. Remark~\ref{rm:5.9}),
and the mean index
$$
\hat i_\tau(\gamma):=\lim_{k\to\infty}\frac{i_{k\tau}(\gamma^k)}{k}.
$$
If $\gamma$ is isolated as a critical point of ${\cal L}_{\tau}$ in
$H_\tau(\alpha)$ then we may define the critical group $C_\ast({\cal
L}_\tau,\gamma;\K)=H_\ast(({\cal L}_\tau)_c\cap U, ({\cal
L}_\tau)_c\cap(U\setminus\{\gamma\});\K)$ as usual, where $c={\cal
L}_\tau(\gamma)$. Furthermore, suppose that $L$ has  global
Euler-Lagrange flow (cf.\cite{AbF, Fa}). By Lemma~ 5.2 of \cite{AbF}
there exists a number $R(A)$ for every $A>0$ such that for any
$R>R(A)$ and for any Lagrangian $L^R$ which is a convex quadratic
$R$-modification of $L$ one has: if $\alpha$ is a critical point of
${\cal L}_\tau^R$ such that ${\cal L}_\tau^R(\alpha)\le A$, then
$\|\dot\alpha\|_\infty\le R(A)$. In particular, such a $\alpha$ is
an extremal curve of ${\cal L}_\tau$ and ${\cal
L}_\tau(\alpha)={\cal L}^R_\tau(\alpha)$. Let us take $A>c$ and
$R>R(A)+ \|\dot\gamma\|_\infty$. Clearly, $\gamma$ is a critical
point of ${\cal L}_\tau|_{X_\tau(\alpha)}$ and hence that of ${\cal
L}^R_\tau|_{X_\tau(\alpha)}$. This implies that $\gamma$ is also
 a critical point of ${\cal L}^R_\tau$ because of the density of
 $X_\tau(\alpha)$ in $H_\tau(\alpha)$. If $\{\gamma_n\}$ is a
 sequence of critical points of ${\cal L}^R_\tau$ converging to
 $\gamma$, we may assume ${\cal L}^R_\tau(\gamma_n)<A\;\forall n$.
Hence each $\gamma_n$ is a critical point of ${\cal L}_\tau$. This
shows that $\gamma$ is an isolated critical point of ${\cal
L}^R_\tau$, and thus $C_\ast({\cal L}_\tau,\gamma;\K)\cong
C_\ast({\cal L}^R_\tau,\gamma;\K)$. According to Section~\ref{sec:1}
the Morse index $m^-({\cal L}_\tau^R,\gamma)$ and nullity
 $m^0({\cal L}_\tau^R,\gamma)$ are well-defined. By
 (\ref{e:3.16})-(\ref{e:3.19}) we have
$m^-({\cal L}_\tau^R,\gamma)=i_{\tau}(\gamma)$ and $m^-({\cal
L}_\tau^R,\gamma)=\nu_{\tau}(\gamma)$. Using these we may derive
that for a Tonelli Lagrangian $L\in C^2(S_\tau\times TM,\R)$ with
global Euler-Lagrange flow under Assumption $F(\alpha)$ in
\cite[page 3010]{Lu}, Claim~5.1 of \cite{Lu} still holds and
Lemma~5.2 of \cite{Lu} becomes:

\noindent{\it Claim 1}.  For each $k\in\N$ there exists
$\gamma'_k\in{\cal K}({\cal L}_{k\tau}, \alpha^{k\tau})$
 such that
$$
C_r({\cal L}_{k\tau}, \gamma'_k; \K)\ne 0\quad{\rm and}\quad r-2n\le
r-\nu_{k\tau}(\gamma'_k)\le i_{k\tau}(\gamma'_k)\le r.
$$

Similarly, for such a system Lemma~5.3 of \cite{Lu} is also true if
$\hat m^-_\tau(\gamma)$ is replaced by $\hat i_\tau(\gamma)$. In
particular, Corollary~5.4 of \cite{Lu} holds if $\hat
m^-_\tau(\gamma_j)$ is replaced by $\hat i_\tau(\gamma_j)$. We have
also the corresponding conclusions in Sections~6,7 of \cite{Lu}. }
\end{remark}

\section{Postscripts}\label{sec:7}

The first draft of this new correction version was completed in
February 2010. The key is to find a new splitting lemma which is
very suitable for our question. We reported the content of the
splitting lemma on conference of Symplectic Geometry and Physics
held at Chern Mathematics Institute on May 17-23, 2010. Then I
concentrate my efforts on developing new theory \cite{Lu2, Lu3}. I
feel most apologetic for submitting this correction version late.

There exists a gap in the original correction \cite{Lu1} with
Jiang's splitting lemma \cite{JM}. That is, we need to prove that
for an isolated critical point $\gamma_0\in X_\tau\subset H_\tau$
the critical groups $C_\ast({\cal L}_\tau, \gamma_0)$ and
$C_\ast({\cal L}_\tau|_{X_\tau}, \gamma_0)$ are isomorphic.  This
may be proved with the finite dimensional approximations of the loop
space by broken Euler-Lagrangian loop space \cite{Maz}, which is an
analogue of the finite dimensional approximations of the loop space
by broken geodesics in Riemannian and Finsler geometry.   Actually
we can also use Theorem~2.10 and Theorem~3.11 in \cite{Lu2} to prove
it.

Using the methods of this paper it seems more easily to generalize
the results in \cite{Lu} to a larger class of Lagrangians as done in
\cite{AbF, Maz} by modifying them to ones satisfying (L1)-(L3) in
\cite{Lu}.

Recently, with Floer homological methods from Ginzburg's proof of
the Conley conjecture \cite{Gi} Doris Hein \cite{He} proved the
existence of infinitely many periodic orbits  for cotangent bundles
of oriented, closed manifolds, and Hamiltonians, which are quadratic
at infinity.

\end{document}